\documentclass[12pt]{amsart}
\usepackage{amssymb}  
\usepackage{latexsym} 
\usepackage[all]{xy}
\usepackage{comment}
\usepackage{url}

\newcommand{\Abar}{\overline{A}}
\newcommand{\Bbar}{\overline{B}}

\newcommand{\C}{{\mathbb C}}

\newcommand{\CC}{{\mathcal C}}
\newcommand{\Char}{\operatorname{char}}

\newcommand{\Diag}{\operatorname{Diag}}

\newcommand{\dd}{\mathfrak{d}}

\newcommand{\dv}{{\,\mid\,}}

\newcommand{\F}{{\mathbb F}}

\newcommand{\G}{{\mathcal G}}
\newcommand{\HH}{{\mathcal H}}
\newcommand{\GL}{\operatorname{GL}}
\newcommand{\Jac}{\operatorname{Jac}}

\newcommand{\Kbar}{\overline{K}}
\newcommand{\kbar}{\overline{k}}

\newcommand{\Magma}{{\sf {MAGMA} }}
\newcommand{\ndv}{\nmid}
\newcommand{\nr}{{\text{\rm nr}}}
\newcommand{\sh}{{\text{\rm sh}}}

\newcommand{\OK}{{\mathcal{O}_K}}
\newcommand{\pp}{\mathfrak{p}}
\newcommand{\pf}{\operatorname{pf}}
\newcommand{\Div}{\operatorname{Div}}
\newcommand{\PGL}{\operatorname{PGL}}
\newcommand{\PP}{{\mathbb P}}
\newcommand{\Q}{{\mathbb Q}}
\newcommand{\Qbar}{\overline{Q}}
\newcommand{\Rbar}{\overline{R}}
\newcommand{\rdet}{\operatorname{rd}}
\newcommand{\R}{{\mathbb R}}
\newcommand{\Real}{{\operatorname{Re}}}
\newcommand{\rank}{\operatorname{rank}}
\newcommand{\ra}{{\longrightarrow}}
\newcommand{\Sel}{\operatorname{Sel}}
\newcommand{\SL}{\operatorname{SL}}
\newcommand{\T}{\mathcal{T}}
\newcommand{\tr}{{\text{\rm tr}}}
\newcommand{\Y}{{\mathcal{Y}}}
\newcommand{\LL}{{\mathcal{L}}}
\newcommand{\oh}{{\mathcal{O}_E}}
\newcommand{\x}{{\bf x}}

\newcommand{\Z}{{\mathbb Z}}
\newcommand{\Aut}{\operatorname{Aut}}
\newcommand{\adj}{\operatorname{adj}}

\newenvironment{spmatrix}{\left(\begin{smallmatrix}}{\end{smallmatrix}\right)}

\newenvironment{Proof}{\par\noindent{\sc Proof:}}%
                      {\hspace*{\fill}\nobreak$\Box$\par\medskip}
\newenvironment{ProofOf}[1]{\par\noindent{\sc Proof of #1:}}%
                       {\hspace*{\fill}\nobreak$\Box$\par\medskip}

\newtheorem{Proposition}{Proposition}[section]
\newtheorem{Theorem}[Proposition]{Theorem}
\newtheorem{Lemma}[Proposition]{Lemma}
\newtheorem{Corollary}[Proposition]{Corollary}

\theoremstyle{definition}
\newtheorem{Definition}[Proposition]{Definition}
\newtheorem{Remark}[Proposition]{Remark}
\newtheorem{Example}[Proposition]{Example}

\numberwithin{equation}{section}

\newcounter{nootje}
\begin{document}

\title[Minimisation and reduction]%
{Minimisation and reduction of 2-, 3- and 4-coverings of elliptic curves}

\author{J.E.~Cremona} 
\address{Mathematics Institute,
         University of Warwick,
         Coventry CV4 7AL, UK}
\email{J.E.Cremona@warwick.ac.uk}

\author{T.A.~Fisher}
\address{University of Cambridge,
         DPMMS, Centre for Mathematical Sciences,
         Wilberforce Road, Cambridge CB3 0WB, UK}
\email{T.A.Fisher@dpmms.cam.ac.uk}

\author{M.~Stoll}
\address{Mathematisches Institut, Universit\"at Bayreuth, 95440
  Bayreuth, Germany}
\email{Michael.Stoll@uni-bayreuth.de}

\date{12th August 2009}  

\begin{abstract}
In this paper we consider models for genus one curves of degree~$n$
for $n=2$, $3$ and~$4$, which arise in explicit $n$-descent on
elliptic curves.  We prove theorems on the existence of minimal models
with the same invariants as the minimal model of the Jacobian elliptic
curve and provide simple algorithms for minimising a given model,
valid over general number fields.  Finally, for genus one models
defined over~$\Q$, we develop a theory of reduction and again give
explicit algorithms for $n=2$, $3$ and~$4$.

\end{abstract}

\maketitle

\renewcommand{\baselinestretch}{1.1}
\renewcommand{\arraystretch}{1.3}

\renewcommand{\theenumi}{\roman{enumi}}


\tableofcontents   



\section{Introduction}
\label{intro}

Let $E$ be an elliptic curve defined over a number field~$K$.  An
$n$-descent on~$E$ computes the $n$-Selmer group
of~$E$, which parametrises the everywhere locally soluble
$n$-coverings of~$E$ up to isomorphism. 
An {\em $n$-covering of~$E$} is a principal
homogeneous space $C$ for~$E$, together with a map $\pi : C \to E$
that fits into a commutative diagram
\[ \xymatrix{ C \ar[dr]^{\pi} \ar@{-->}[d]_{\psi} \\
              E \ar[r]^{\cdot n} & E
            }
\]
where $\psi : C \to E$ is an isomorphism defined over the algebraic
closure~$\overline{K}$, compatible with the structure of~$C$ as a
principal homogeneous space. In a series of papers~\cite{ndescent}, it
is shown how to produce explicit equations of covering curves from a
more abstract representation of the Selmer group.
(The latter is computed, at least for $n$ prime, in \cite{SchaeferStoll}.)

In general, an $n$-covering~$C$ can be realised as a smooth curve of
degree~$n$ inside a Severi-Brauer variety~$S$ of dimension~$n-1$ (when
$n = 2$, we obtain a double cover of a conic instead of an embedding).
If $C$ has points everywhere locally, as will be the case when $C$
represents an element of the $n$-Selmer group of~$E$, then the same
statement is true of~$S$, and hence $S \cong \PP^{n-1}$, so that $C$
has a degree-$n$ model in projective space. Thus, for $n =
2$, we get a double cover of~$\PP^1$ ramified in four points, for $n =
3$, we get a plane cubic curve, and for $n = 4$, we get an
intersection of two quadrics in~$\PP^3$. For larger~$n$, these models
are no longer complete intersections, but can be given by a number of
quadratic equations.

In this paper, we will focus on the problem of how to produce ``nice''
models of the covering curves, i.e., models given by equations with
small integral coefficients, in the cases $n = 2$, $3$ and~$4$. The
advantage of having such a nice model is two-fold.  On the one hand,
rational points on the covering curve can be expected to be of smaller
height on a model with small coefficients, and therefore will be found
more easily. On the other hand, if no rational points are found, one
would like to use the covering curve as the basis for a further
descent, and the necessary computations are greatly facilitated when
the given model is nice.

This problem naturally splits into two parts: {\em Minimisation} and
{\em Reduction}.  Minimisation makes the invariants of the model
smaller by eliminating spurious bad primes and reducing the exponents
of primes of bad reduction, to obtain a ``minimal model''. 
We prove the following theorem.  (See Section~\ref{g1mod} for the
definitions of models for $n$-coverings and their invariants.)
\begin{Theorem}
\label{thmglobal} 
  Let $n=2,3$ or $4$. Let $K$ be a number field of class number one,
  and $E$ an elliptic curve defined over $K$.  If $\CC$ is an
  $n$-covering of~$E$ which is everywhere locally soluble (i.e. $\CC$
  has points over all completions of~$K$) then $\CC$ has a model with
  integral coefficients and the same discriminant as 
  a global minimal Weierstrass equation for~$E$.
\end{Theorem}
By contrast, reduction attempts to reduce the size of the coefficients
by an invertible integral (i.e., unimodular) linear change of
coordinates, which leaves the invariants unchanged. Both processes are
necessary to obtain a nice model: minimisation without reduction
will provide a model with small invariants, but most likely rather
large coefficients, whereas reduction without minimisation will not be
able to make the coefficients really small, since the invariants will
still be large.

After introducing the kinds of models 
we will be using and their invariants in Section~\ref{g1mod}, we state
our main results on minimisation over local fields in 
Section~\ref{min:statements}, and discuss how they relate 
to earlier work. The most important of these results 
(the Minimisation Theorem, Theorem~\ref{minthm}) 
is proved in Section~\ref{sec:unproj}. The proof is short and 
transparent, but is not algorithmic. We remedy this in
Section~\ref{sec:M} where we give practical algorithms for 
computing minimal models, that may be seen as generalising 
Tate's algorithm~\cite{Tate}. In Section~\ref{sec:globalmin} we
deduce Theorem~\ref{thmglobal} from our local results, and explain how
it may be generalised to arbitrary number fields. Moreover, as our local 
minimisation results make no restriction on the characteristic of 
the local field, they have more general global applications; 
in particular, one obtains results over function fields as well 
as number fields.

The algorithms of Section~\ref{sec:M} may be combined with 
the Minimisation Theorem to prove the Strong Minimisation Theorem
(Theorem~\ref{mainthmA} (i)). This states that if an $n$-covering of
$E$ (defined over a local field, and represented by a degree-$n$ model) 
is soluble over the maximal unramified extension, then it has a
model with integral coefficients and the same discriminant as a minimal
Weierstrass equation for $E$. In Section~\ref{sec:critmodels} we prove
the converse (Theorem~\ref{mainthmA} (ii)), thereby showing that the 
Strong Minimisation Theorem is best possible.

In Section~\ref{sec:R} we discuss reduction for general
$n$-coverings, and more specifically for $n = 2$, $3$ and~$4$.
Our results for reduction only cover the case where the ground field
is~$\Q$.  A comparable theory of reduction over a general number field
would be very useful in practice, but has not yet been sufficiently
developed. We end in Section~\ref{sec:E} by giving some examples of
both minimisation and reduction (over $K= \Q$). 
All our algorithms (for $n=2,3,4$ and $K = \Q$) have been implemented in 
(and contributed to) \Magma~(see \cite{magma}).

\medskip

As stated earlier, the main application of our results is in explicit
$n$-descent on elliptic curves over number fields.  Minimisation and
reduction of binary quartics is also used in the invariant theory
method for $2$-descent (see \cite{BSD1} and \cite{Cremona}).  For
$n=3$, Djabri and Smart in their ANTS III article~\cite{DS} consider
the possibility of carrying out $3$-descent using invariant theory in
a similar way; one stumbling-block there was the inability to minimise
plane cubic models for $3$-coverings.


\section{Genus one models}
\label{g1mod}

In this section, we specify the models of the covering curves that we
will use, together with their invariants $c_4$, $c_6$, and~$\Delta$.
For completeness and later reference we include the case $n=1$.  Note
that we use the term ``genus one model'' to include singular models,
which do not define curves of genus one.

\begin{Definition}
A {\em Weierstrass equation}, or {\em genus one model of degree 1}, is an 
equation of the form
\begin{equation}
\label{weqn1}
 y^2 + a_1 x y + a_3 y = x^3 + a_2 x^2 + a_4 x + a_6. 
\end{equation}
The space of all Weierstrass equations with coefficients $a_1, \ldots,
a_6$ in a ring $R$ will be denoted $X_1(R)$. We say that two such
models are {\em $R$-equivalent} if they are related by substitutions
\begin{equation}
\label{subst:weqn}
x \leftarrow u^2 x + r \qquad y \leftarrow u^3 y + u^2 s x + t
\end{equation}
for some $u \in R^\times$ and $r,s,t \in R$.  We write $\G_1(R)$ for
the group of all transformations $[u;r,s,t]$ and define
$\det([u;r,s,t])= u^{-1}$.  The invariants $c_4$, $c_6$ and $\Delta$
are certain primitive polynomials in $a_1, \ldots, a_6$ with integer
coefficients, satisfying $c_4^3-c_6^2=1728\Delta$ (see
e.g. \cite[Chapter~III]{Si1}).
\end{Definition}

\begin{Definition}
\label{def:g1m2}
  A {\em genus one model of degree 2}, or {\em generalised binary
    quartic}, is an equation of the form $$ y^2 + P(x,z) y = Q(x,z) $$
  where $P$ and $Q$ are homogeneous polynomials of degrees 2 and 4.
  We sometimes abbreviate this as $(P,Q)$.  The space of all such
  models with coefficients in a ring $R$ is denoted $X_2(R)$.  Two
  such models are {\em $R$-equivalent} if they are related by
  substitutions $x \leftarrow m_{11} x + m_{21} z$, $z \leftarrow
  m_{12} x + m_{22} z$ and $y \leftarrow \mu^{-1} y + r_0 x^2 + r_1 x
  z + r_2 z^2$ for some $\mu \in R^\times$, $r = (r_0,r_1,r_2) \in
  R^3$ and $M = (m_{ij}) \in \GL_2(R)$. We write $\G_2(R)$ for the
  group of all such transformations~$[\mu,r,M]$, and define
  $\det([\mu,r,M]) = \mu \det (M)$.

  A generalised binary quartic $y^2 + P(x_1,x_2)y = Q(x_1,x_2)$ over a
  field $K$ defines a subscheme $\CC_{(P,Q)} \subset \PP(1,1,2)$, the
  ambient space being a weighted projective space with coordinates
  $x_1$, $x_2$, $y$.  The model $\Phi=(P,Q)$ is {\em $K$-soluble} if
  $\CC_{\Phi}(K) \not= \emptyset$.

  The binary quartic
  $F(x,z) = a x^4 + b x^3 z + c x^2 z^2 + d x z^3 + e z^4$
  has invariants $c_4(F) = 2^4 I$ and $c_6(F) = 2^5 J$, where~$I$
  and~$J$ are given by
  \begin{align*}
  I & = 12 a e - 3 b d + c^2, \\
  J & = 72 a c e - 27 a d^2 - 27 b^2 e + 9 b c d - 2 c^3.
  \end{align*}
  The discriminant $\Delta= (c_4^3- c_6^2)/1728$ is 16 times
  the usual discriminant of a quartic polynomial.
  The invariants of a generalised binary quartic are 
  obtained by completing the square, {\em i.e.}
  $c_4(P,Q)= c_4(\frac{1}{4} P^2+Q)$ and so on.
  We find that $c_4$, $c_6$ and $\Delta$ are primitive integer 
  coefficient polynomials in the coefficients of $P$ and $Q$, again
  satisfying $c_4^3-c_6^2=1728\Delta$.
\end{Definition}

Earlier work on $2$-coverings, including \cite{BSD1} and
\cite{SCmin2}, used the more restrictive binary quartic models
with~$P=0$.  We use generalised binary quartics here, in order to
obtain more uniform local results at places with residue
characteristic~$2$.

\begin{Definition}
  A {\em genus one model of degree 3} is a ternary cubic. We write 
  $X_3(R)$ for the space of all ternary cubics with coefficients
  in a ring $R$. Two such models are {\em $R$-equivalent} if they are related by 
  multiplying by $\mu \in R^\times$ and then 
  substituting $x_j \leftarrow \sum_{i=1}^3 m_{ij} x_i$ 
  for some $M = (m_{ij}) \in \GL_3(R)$. We write $\G_3(R) = R^\times
  \times \GL_3(R)$ 
  for the group of all such transformations~$[\mu,M]$, and define  
  $\det([\mu,M]) = \mu \det (M)$.

  A ternary cubic $F(x,y,z)$ over a field $K$ 
  defines a subscheme $\CC_F \subset \PP^2$.
  The model $F$ is {\em $K$-soluble} if $\CC_F(K) \not= \emptyset$.

  The invariants $c_4$ and~$c_6$ may be defined as follows. Let
  \begin{equation*} 
    H(F) =  \det\left(\begin{array}{ccc}
                       F_{xx} & F_{xy} & F_{xz} \\
                       F_{yx} & F_{yy} & F_{yz} \\
                       F_{zx} & F_{zy} & F_{zz}
                     \end{array}\right)
  \end{equation*}
  be the {\em Hessian} of~$F$, which is again a ternary cubic.  
  Then we have
  \[ H(H(F)) = 48\,c_4(F)^2 F + 16\,c_6(F) H(F) \,; \]
  the sign of~$c_4(F)$ is fixed by requiring that 
  $\Delta = (c_4^3 - c_6^2)/1728$ has integer coefficients.
  Then $c_4$, $c_6$ and $\Delta$ are primitive integer coefficient 
  polynomials in the coefficients of $F$ and
  satisfy $c_4^3-c_6^2=1728\Delta$.
\end{Definition}

\begin{Definition}
\label{def:g1m4}
  A {\em genus one model of degree 4}, or {\em quadric intersection},
  is an ordered pair~$(Q_1,Q_2)$ of quadrics (homogeneous polynomials of
  degree 2) in four variables. The space
  of all such models with coefficients in a ring $R$ is denoted
  $X_4(R)$.  Quadric intersections $(Q_1,Q_2)$ and $(Q_1',Q_2')$ are
  {\em $R$-equivalent} if they are related by putting $Q'_1 = m_{11}
  Q_1 + m_{12}Q_2$ and $Q'_2 = m_{21} Q_1 + m_{22}Q_2$ for some $M =
  (m_{ij}) \in \GL_2(R)$ and then substituting $x_j \leftarrow
  \sum_{i=1}^4 n_{ij} x_i$ for some $N = (n_{ij}) \in \GL_4(R)$.  We
  write $\G_4(R) = \GL_2(R) \times \GL_4(R)$ for the group of all such
  transformations~$[M,N]$, and define $\det([M,N]) = \det(M) \det (N)$.

  A quadric intersection $\Phi=(Q_1,Q_2)$ over a field $K$ defines a
  subscheme $\CC_{\Phi} \subset \PP^3$.  The model $\Phi$ is {\em
    $K$-soluble} if $\CC_{\Phi}(K) \not= \emptyset$.

  The invariants $c_4$ and~$c_6$ may be defined 
  as follows. Let $A$ and $B$ be the matrices of 
  second partial derivatives of $Q_1$ and $Q_2$. 
  Then $F(x,z) = \det(Ax + Bz)$ is a binary quartic. We define 
  $c_4(Q_1,Q_2) = 2^{-4} c_4(F)$, $c_6(Q_1,Q_2) = 2^{-6} c_6(F)$
  and $\Delta(Q_1,Q_2) = 2^{-12} \Delta(F)$. 
  These scalings are chosen so that $c_4$, $c_6$ and $\Delta$ 
  are primitive integer coefficient polynomials in the 
  coefficients of $Q_1$ and~$Q_2$. They satisfy $c_4^3-c_6^2=1728\Delta$.
\end{Definition}

Earlier work on $4$-coverings, including \cite{WomackThesis} and
\cite{SiksekThesis}, used pairs of symmetric matrices rather than
pairs of quadrics.  We use quadrics here, in order to obtain more
uniform local results at places with residue characteristic~$2$.

\begin{Remark}
  There is also a definition of {\em genus one model of degree
  $5$}, see \cite{g1inv}. The minimisation and reduction of
  these models (and possible extensions to larger $n$) 
will be the subject of future investigations. 
\end{Remark}

\begin{Remark}
\label{rem1}
There is a natural way in which we can re-write 
a Weierstrass equation (a genus one model of degree 1)
as a genus one model of degree $n=2,3$ or $4$ 
(see Lemma~\ref{lem1}). We have normalised the invariants
$c_4$, $c_6$ and $\Delta$ so that they agree with the usual formulae 
(see e.g. \cite[Chapter III]{Si1}) when specialised to 
one of these `Weierstrass models'.
\end{Remark}

\begin{Definition}
\label{def:inv}
  Let $K$ be a field and $\Kbar$ its algebraic
  closure. Let $K[X_n]$ be the polynomial ring in the 
  coefficients of a genus one model of degree $n$. 
  A polynomial $F \in K[X_n]$ is an {\em invariant} of {\em weight} $k$
  if $F \circ g = \det(g)^k F$ for all $g \in \G_n(\Kbar)$.
\end{Definition}

For $n=1,2,3,4$ we defined polynomials $c_4, c_6, \Delta \in \Z[X_n]$
with $c_4^3 - c_6^2 = 1728 \Delta$. These have the following
properties.

\begin{Theorem}
\label{thm:invjac}
  Let $n =1,2,3$ or $4$.
  \begin{enumerate}
    \item The polynomials $c_4, c_6, \Delta \in K[X_n]$ are invariants
          of weights $4$, $6$ and $12$.
    \item A genus one model $\Phi \in X_n(K)$ defines a smooth
      curve~$\CC_{\Phi}$ of genus one (over $\Kbar$) if and only if
      $\Delta(\Phi) \not= 0$.
    \item If $\Char(K) \not= 2,3$ then $c_4$ and $c_6$ generate the
      ring of invariants. Moreover if $\Phi \in X_n(K)$ with
      $\Delta(\Phi) \not= 0$ then the Jacobian of the curve
      $\CC_{\Phi}$ has Weierstrass equation
          $$ y^2 = x^3 - 27 c_4(\Phi) x - 54 c_6(\Phi). $$
  \end{enumerate}
\end{Theorem}
\begin{Proof}
The invariants $c_4$, $c_6$ and $\Delta$ were known to the nineteenth
century invariant theorists. The observation that they give a formula
for the Jacobian is due to Weil \cite{Weil1}, \cite{Weil2}. See
\cite{AKM3P} for a brief survey, or \cite{g1inv} for a proof of the
theorem exactly as it is stated here.
\end{Proof}

As was first pointed out to us by Rodriguez-Villegas,
it is possible to work back through Tate's formulaire 
(see e.g. \cite[Chapter III]{Si1}) to write the invariants 
$c_4$ and $c_6$ in terms of polynomials $a_1, \ldots, a_6$.

\begin{Lemma}
  \label{lem:abinv}
There exist $a_1,a_2,a_3,a_4,a_6 \in \Z[X_n]$ 
  and $b_2,b_4,b_6 \in \Z[X_n]$ with
  \begin{align*}
   b_2 & =  a_1^2+4 a_2     & c_4 & = b_2^2-24 b_4 \\ 
   b_4 & =  a_1 a_3 + 2 a_4 & c_6 & =  -b_2^3 + 36 b_2 b_4 - 216 b_6. \\
   b_6 & =  a_3^2 + 4 a_6 
  \end{align*} 
\end{Lemma}
\begin{Proof}
The lemma is proved by splitting into the cases $n=2,3,4$ 
and giving explicit formulae for the $a$-invariants. 
(The case $n=1$ is a tautology.)

\paragraph{Case $n=2$.} 
The $a$-invariants of the generalised binary quartic
$$ y^2 + (l x^2 + m xz + n z^2) y 
     = a x^4 + b x^3 z + c x^2 z^2 + d x z^3 + e z^4 $$
are
\begin{equation}
  \label{jac2gen}
  \begin{aligned}
    a_1 & =  m \\
    a_2 & =  c- ln  \\
    a_3 & =  l d + n b \\
    a_4 & =  - 4 a e + b d  -(l^2 e + ln c + n^2 a)  \\
    a_6 & =  - 4 a c e + a d^2 + b^2 e 
    - (l^2 c e + m^2 a e + n^2 a c + l n b d ) + lm b e  + mn a d.  
  \end{aligned} 
\end{equation}

\paragraph{Case $n=3$.} 
The $a$-invariants of the ternary cubic
$$  a x^3 + b y^3 + c z^3 + a_2 x^2 y + a_3 x^2 z +  b_1 x y^2 
                     + b_3 y^2 z + c_1 x z^2  + c_2 y z^2 + m x y z $$
are
\begin{equation}
  \label{jac3gen}
  \begin{aligned} \smallskip
    a_1 & =  m \\ \smallskip
    a_2 & =  -(a_2 c_2+a_3 b_3+b_1 c_1) \\ \smallskip
    a_3 & =  9 a b c - (a b_3 c_2 + b a_3 c_1 + c a_2 b_1)
                 - (a_2 b_3 c_1 + a_3 b_1 c_2) \\
    a_4 & =  -3 (a b c_1 c_2 + a c b_1 b_3  + b c a_2 a_3) \\
        &    \quad + \, a (b_1 c_2^2 + b_3^2 c_1) 
                 + b (a_2 c_1^2 + a_3^2 c_2) 
                 + c (a_2^2 b_3 + a_3 b_1^2) \\ \smallskip
        &    \quad + \, a_2 c_2 a_3 b_3 + b_1 c_1 a_2 c_2 + a_3 b_3 b_1 c_1 \\
    a_6 & =  -27 a^2 b^2 c^2 
               + 9 a b c (a b_3 c_2 + c a_2 b_1 + b a_3 c_1)
               + \ldots +  a b c m^3. 
  \end{aligned}
\end{equation}
These formulae in the case $n=3$ were first given in \cite{ARVT}. 

\paragraph{Case $n=4$.} 
Let $Q = \sum_{i \le j} c_{ij} x_i x_j$ be a quadric in  4 variables. 
Then  
 \[ \det (\tfrac{\partial^2 Q}{\partial x_i \partial x_j} )
  = \pf(Q)^2 + 4 \rdet(Q) \]
where $\pf(Q) = c_{12} c_{34} + c_{13} c_{24} + c_{14} c_{23}$ and 
$\rdet(Q) \in \Z[c_{11},c_{12}, \ldots, c_{44}]$.
We define the $a$-invariants of the quadric intersection
$(Q_1,Q_2)$ to be the $a$-invariants of the generalised binary quartic
\[  y^2 + \pf (x Q_1 + z Q_2) y = \rdet( x Q_1 + z Q_2). \]

\end{Proof}

The polynomials $a_i$ of Lemma~\ref{lem:abinv} are far from unique.
They can be modified by any transformation in $\G_1(\Z[X_n])$,  
i.e. by any transformation of the form $[\pm 1;r,s,t]$
with $r,s,t \in \Z[X_n]$. The following theorem extends
Theorem~\ref{thm:invjac}(iii) to fields of arbitrary characteristic.
(The reader only interested in applications 
over number fields and their completions, may safely skip this result.)

\begin{Theorem}
\label{thm:invjac-a}
 Let $K$ be any field, and $n =1,2,3$ or $4$. For all $\Phi \in
 X_n(K)$ with $\Delta(\Phi) \not= 0$, the Jacobian of the curve $\CC_{\Phi}$ has
 Weierstrass equation
         \begin{equation}
           \label{weqn(Phi)}
           y^2 + a_1 xy + a_3 y
             = x^3 + a_2 x^2 + a_4 x + a_6
\end{equation}
where $a_i = a_i(\Phi)$.
\end{Theorem}

\begin{Proof}
For $n=3$ this is a special case of a theorem of Artin,
Rodriguez-Villegas and Tate \cite{ARVT}. The cases $n=2,4$ may be
proved using similar techniques. We sketch a simplified form of the
proof, covering the cases $n=2$, $3$, and~$4$. (The case $n=1$ is of
course a tautology.)

Let $C/S$ be the universal family over\footnote{In \cite{ARVT}
the authors work over $S=  \operatorname{Spec}(\Z[X_3]) \setminus \{0\}$.
This gives a more general result, but also makes the proof more difficult.}
$S = \operatorname{Spec} (\Z[X_n][\Delta^{-1}])$.
By Theorem~\ref{thm:invjac}(ii) the fibres are smooth projective
curves of genus one. Let $J/S$ be the Jacobian of $C/S$, in the
sense that $J$ is the $S$-scheme representing the relative
Picard functor $\operatorname{Pic}^0_{C/S}$; see \cite[\S9.3, Theorem 1]{BLR}.
Each fibre of $J/S$ is the Jacobian of the corresponding fibre of
$C/S$ and hence an elliptic curve. By a generalisation of the
usual procedure for putting an elliptic curve in Weierstrass form
(see \cite{De}, or \cite[Theorem 2]{ARVT} for a further generalisation)
$J$ is defined as a subscheme of $\PP^2_S$ by the homogenisation of
\begin{equation}
\label{weqn'}
 y^2 + a'_1 xy + a'_3 y = x^3 + a_2' x^2 + a_4' x + a_6'
\end{equation}
for some $a'_1, \ldots, a'_6 \in \Z[X_n][\Delta^{-1}]$. Thus for
every field $K$, and every $\Phi \in X_n(K)$ with $\Delta(\Phi) \not= 0$,
the Weierstrass equation~(\ref{weqn'}) gives a model for
the Jacobian of~$\CC_{\Phi}$.

It only remains to show that~(\ref{weqn(Phi)}) and~(\ref{weqn'}) are
related by a transformation in $\G_1(R)$ where $R =
\Z[X_n][\Delta^{-1}]$. By Theorem~\ref{thm:invjac}(iii) they are
related by some $[u;r,s,t] \in \G_1(K)$ where $K= \Q(X_n)$.
Since for any genus one model with~$\Delta\not=0$, (\ref{weqn(Phi)})
and~(\ref{weqn'}) both specialise to a non-singular Weierstrass
equation, it follows that $u \in R^\times$.  Then, since $R$ is
integrally closed, a standard argument (see \cite[Chapter VII,
  Proposition 1.3]{Si1}) shows that $r,s,t \in R$.
\end{Proof}

We note that $a_1, \ldots, a_6$ are not invariants in the sense
of Definition~\ref{def:inv}.
The ring of invariants when $\Char(K)=2$ or $3$ is described in
\cite[\S10]{g1inv}. As is noted there, these do not give a
formula for the Jacobian.


\section{Minimisation theorems}
\label{sec:T}

\subsection{Statement of results}
\label{min:statements}

Let $K$ be a field with normalised discrete valuation 
$v : K^\times \to \Z$. We write $\OK$ for the valuation ring
(or ring of integers) of $K$ and fix a uniformiser $\pi \in K$. 
We assume throughout that the residue field $k = \OK/\pi \OK$ is perfect.
A field extension $L/K$ is {\em unramified} if there
is a (normalised) discrete valuation $w:L^\times \to \Z$ extending $v$.
The strict Henselisation $K^{\sh}$ of $K$ is an unramified
extension of $K$ that satisfies the conclusions of 
Hensel's lemma and has residue field $\kbar$, the algebraic 
closure of $k$.
(See \cite[Definition~4.18]{MilneETALE} for the precise definition.) 
If $K$ is complete (with respect to $v$) 
then $K^{\sh}$ is the maximal unramified extension $K^\nr$ of $K$ 
as defined in \cite[Chapter III, \S5]{SerreLF}.

We work with genus one models
of degree $n = 1,2,3$ or $4$.
The invariants $c_4$, $c_6$ and $\Delta$ 
of a genus one model were defined in Section~\ref{g1mod}.

\begin{Definition}
\label{mindefs}\strut
\begin{enumerate}
  \item A genus one model $\Phi \in X_n(K)$ is {\em non-singular} if 
        $\Delta(\Phi) \not= 0$.
  \item A genus one model $\Phi \in X_n(K)$ is {\em integral} if it has 
        coefficients in $\OK$.
  \item A non-singular model $\Phi \in X_n(\OK)$ is {\em minimal} 
        if $v(\Delta(\Phi))$ is minimal among all integral
        models $K$-equivalent to $\Phi$, otherwise $\Phi$ is {\em non-minimal}.
\end{enumerate}
\end{Definition}

Algorithms for computing minimal models in the case $n=1$ have
been given by Tate~\cite{Tate}, \cite[Chapter IV, \S9]{Si2} 
and Laska~\cite{Laska}. 
The latter can be refined using Kraus' conditions \cite{Kraus} 
as described in \cite[Chapter V]{Connell} or \cite[\S3.2]{Cremona}.
(Laska's algorithm and its refinements are simpler 
than Tate's algorithm, but are only applicable when $\Char(K) \not=2,3$.)
In Section~\ref{sec:M} we give algorithms for computing minimal 
models in the cases $n = 2,3,4$. 

In the following lemma we define the {\em level} of a genus one model.

\begin{Lemma}
\label{lem:deflevel}
Let $\Phi \in X_n(K)$ be a non-singular model of degree $n$.
Let $\Delta_E$ be the minimal discriminant of $E = \Jac(\CC_\Phi)$. Then 
\begin{enumerate}
\item $v(\Delta(\Phi)) = v(\Delta_E) + 12 \ell$ 
for some integer $\ell$, called the {\em level} of $\Phi$.
\item If $\Char(k) \not= 2,3$ then $\ell = \min \{ \lfloor v(c_4(\Phi))/4
\rfloor, \lfloor v(c_6(\Phi))/6 \rfloor \}$. 
\item The level of an integral model is always non-negative.
\end{enumerate}
\end{Lemma}
\begin{Proof}
If $\Char(k) \not= 2,3$ then this is clear by Theorem~\ref{thm:invjac}
and the standard formulae for transforming Weierstrass equations.
In general (that is, to prove (iii) when $\Char(k)= 2$ or $3$, 
or even to define the level when $\Char(K)=2$ or $3$) 
we use Lemma~\ref{lem:abinv} and Theorem~\ref{thm:invjac-a} instead.
\end{Proof}

The level of $\Phi \in X_n(K)$ may be computed as 
$v(u)$ where $[u;r,s,t] \in \G_1(K)$ is a transformation that minimises
the Weierstrass equation~(\ref{weqn(Phi)}). 

\begin{Definition}
\label{mindefs2}
The {\em minimal level} of $\Phi \in X_n(K)$ is the minimum of the
levels of all integral models $K$-equivalent to $\Phi$. Thus an integral
model $\Phi$ is minimal (see Definition~\ref{mindefs}) 
if and only if it has level equal to this minimal level.
\end{Definition}

If $n=1$ then the minimal level is $0$, for trivial reasons.  So from
now on we take $n = 2,3$ or $4$. The most important result on minimisation
states that every $K$-soluble model has minimal level~$0$, or in other
words, is $K$-equivalent to an integral model
whose discriminant has the same valuation as the discriminant
of a minimal model for the Jacobian elliptic curve.

\begin{Theorem}[Minimisation theorem]
\label{minthm}
Let $\Phi \in X_n(K)$ be non-singular.
If $\CC_{\Phi}(K) \not= \emptyset$ then $\Phi$ has minimal level $0$.
\end{Theorem}

The following strengthening of the Minimisation Theorem shows
that a non-singular model 
has minimal level $0$ if and only if it is $K^\sh$-soluble.

\begin{Theorem}
\label{mainthmA}
Let $\Phi \in X_n(K)$ be non-singular.
\begin{enumerate}
\item {\em (Strong Minimisation Theorem)}. 
If $\CC_{\Phi}(K^\sh) \not= \emptyset$ then $\Phi$ has minimal level~$0$.
\item {\em (Converse Theorem)}. If $\CC_{\Phi}(K^\sh) = \emptyset$ then
the minimal level is at least $1$, and is equal to $1$ if $\Char(k) \ndv n$. 
\end{enumerate}
\end{Theorem}

Algorithms for minimising $K$-soluble binary quartics over $K= \Q_p$
are sketched by Birch and Swinnerton-Dyer~\cite[Lemmas 3,4,5]{BSD1},
with details in the case~$p\not=2,3$.  Their algorithms give a proof
of the Minimisation Theorem for~$n=2$, except when $p=2$ (in which
case further work is required to handle the ``cross terms''). As
pointed out in \cite{SCmin2} this generalises immediately to any local
field $K$ with $\Char(k) \not= 2,3$. The authors extended
these calculations to the case $n=3$ in conjunction with their work
on $3$-descent \cite{ndescent}. 
The case $n=4$ was treated by Womack
in his PhD thesis \cite[Section 2.5]{WomackThesis}, using a method
that goes via the results for $n=2$.

In each case, the approach taken is to 
start with a $K^\sh$-soluble model $\Phi \in X_n(\OK)$ with 
$v(c_4(\Phi)) \ge 4$ and $v(c_6(\Phi)) \ge 6$, and then 
by a series of substitutions to show that $\Phi$ is 
$K$-equivalent to an integral model of smaller level. This leads to
both a proof of the Strong Minimisation Theorem and a practical algorithm for
minimising.  However, this traditional approach suffers from the
following drawbacks.
\begin{itemize}
\item It is necessary to split into a large number of 
(elementary yet tedious) cases, and the number of
cases grows rapidly with $n$. 
\item The modifications required if $\Char(k) = 2$ or $3$ are
somewhat involved. (The hypothesis that $\Phi$ has 
positive level has to made explicit using either Kraus' conditions
\cite{Kraus} or the ``$a$-invariants'' defined in Lemma~\ref{lem:abinv}.)
\end{itemize}

We take a different approach, in which the 
tasks of proving the Minimisation Theorem and finding a 
practical algorithm for minimising are treated separately.  
A proof of the Minimisation Theorem for $n=2,3$ (in all
residue characteristics) is given in \cite{minbqtc}. 
In Section~\ref{sec:unproj} we simplify the proof and extend to the case $n=4$.
Unfortunately this approach does not lead to any readily 
implementable algorithm, nor does it prove the
Strong Minimisation Theorem. 

In Section~\ref{sec:minbq} (case $n=2$) and Section~\ref{M3} (case $n=3$) we
specify a rather simple-minded procedure and show that, given 
any non-minimal integral model, iterating this procedure will 
eventually decrease the level. This gives an algorithm 
for computing minimal models. In Section~\ref{M4} we give an algorithm in 
the case $n=4$ based on the treatment in Womack's thesis. 
The algorithms for $n=2, 4$ must be modified when
the residue characteristic is $2$, 
as described in Section~\ref{min_char2}. 
These modifications are required since, as noted in Section~\ref{g1mod},
our models for 
$n$-coverings differ slightly from those used previously in the literature.
We have also defined the level, not in an absolute way, 
but by comparison with a minimal model for the Jacobian elliptic curve. 
The combined effect of these changes is that our results are much cleaner 
to state, in particular for residue characteristic~$2$, and can be 
proved uniformly, without assumptions on the ramification index.

\medskip

As is the case for Tate's algorithm, it is clear from
the form of our algorithms (for $n=2,3,4$) that their success 
or otherwise is unchanged by an unramified field extension. 
We deduce the following.
\begin{Theorem}
\label{mainthmB1}
The minimal level of a non-singular genus one model of 
degree $n=2,3$ or $4$ is unchanged by an unramified field extension.
\end{Theorem}
The Strong Minimisation Theorem is then an immediate consequence of 
Theorem~\ref{mainthmB1} and the Minimisation Theorem. 

In Section~\ref{sec:critmodels} we show how to write down examples of
minimal genus one models of positive level. We call the models arising
in our construction {\em critical models}, see
Definition~\ref{defcrit} below.  We show (for $n=2,3$) that any
$K^\sh$-insoluble model is $K$-equivalent to a critical model. There
is a corresponding result for models of degree $n=4$. The proof of the
Converse Theorem (Theorem~\ref{mainthmA}(ii)) is then reduced to a
statement about the possible levels of a critical model (see
Lemma~\ref{lem:critlev}).

\medskip

Theorem~\ref{mainthmA} in the case $n=2$ 
may already be found in \cite[Remarque 21]{Liu}.
We claim that our proof is simpler, and in any case serves 
as a template for our generalisations to $n=3,4$. 
Liu also gives an algorithm for 
minimising~\cite[p.4594, Remarque 11]{Liu} (still for $n=2$), 
which although not made explicit appears to be the same as ours.

We remark that minimisations are not unique, in the sense that
there can be more than one $\OK$-equivalence class of minimal
models $K$-equivalent to a given genus one model. Following on 
from our work and that of Liu, it will be explained in \cite{Sadek}
how to compute the number of such classes.

For a more general, but necessarily less explicit, discussion of the
problem of minimising homogeneous polynomials (of degree $d$ in $n$
variables) see \cite{Kollar}.

\subsection{Proof of the minimisation theorem}
\label{sec:unproj}

In this section only we relax our assumptions on $\OK$ and $K$.
It will only be necessary to assume that $\OK$ is a 
principal ideal domain and $K$ is its field of fractions. 
The definitions of a non-singular model and an integral model
(see Definition~\ref{mindefs}) carry over as before.
We consider models of degree $n=2,3$ or $4$.

Let $E$ be an elliptic curve over $K$, with identity $\oh \in E(K)$,
and let $D$ be a $K$-rational divisor on $E$ of degree $n$.  We write
$[D]$ for the linear equivalence class of $D$.  We pick a basis $f_1,
\ldots,f_n$ for the Riemann-Roch space $\LL(D)$, and let $E \to
\PP^{n-1}$ be the morphism given by $P \mapsto (f_1(P): \ldots :
f_n(P))$. Then according as $n=2,3$ or $4$, we find that $E$ may be written
as either a double cover of $\PP^1$, a plane cubic, or an
intersection of two quadrics in $\PP^3$.  It is therefore defined by
a suitable genus one model~$\Phi\in X_n(K)$. Moreover this model is
uniquely determined up to $K$-equivalence by the pair $(E,[D])$:
replacing~$D$ by an equivalent divisor or changing basis for the
space~$\LL(D)$ only has the effect of a linear change of coordinates
on~$\PP^{n-1}$, so only changes the genus one model by a
$K$-equivalence.
In this situation we say that the genus one model~$\Phi$ 
{\em represents} the pair~$(E,[D])$.

Similarly, we obtain a genus one model~$\Phi\in X_n(K)$, well-defined
up to $K$-equivalence, representing every pair~$(\CC,[D])$ where~$\CC$
is a genus one curve and $D$ a divisor of degree~$n$ on~$\CC$; we have
$\CC\cong\CC_{\Phi}$ (over~$K$), and in particular, $\Phi$ is
$K$-soluble if and only if~$\CC(K)\not=\emptyset$.  Under this
isomorphism, the divisor class~$[D]$ on~$\CC$ maps to a distinguished
divisor class~$[D_{\Phi}]$ of degree~$n$ on~$\CC_{\Phi}$, namely the
class of the fibres of the map $\CC_{\Phi} \to \PP^1$ if $n=2$, or the
hyperplane section if $n=3,4$.  

The following is now immediate.

\begin{Lemma}
\label{whyD}
Every $K$-soluble non-singular genus one model represents a pair
$(E,[D])$ in the manner described above.
\end{Lemma}
\begin{Proof}
Let $\Phi \in X_n(K)$ be a non-singular genus one model. 
It is a tautology that $\Phi$ represents the 
pair~$(\CC_{\Phi},[D_{\Phi}])$. Now if $\Phi$ is $K$-soluble then 
$\CC_{\Phi}$ is a smooth curve of genus one with a rational point,
and hence is an elliptic curve.  
\end{Proof}

The aim of this section is to prove the following theorem.
The Minimisation Theorem (Theorem~\ref{minthm}) is then an immediate
consequence by Lemma~\ref{whyD}.

\begin{Theorem}
\label{thmstar}
Let $E/K$ be an elliptic curve with integral Weierstrass equation
\begin{equation}
\label{weqnstar}
y^2 + a_1 xy + a_3 y = x^3 + a_2 x^2 + a_4 x + a_6 
\end{equation}
and let $D \in \Div_K(E)$ be a divisor on $E$ of degree $n=2,3$ or $4$. 
Then $(E,[D])$ can be represented by an integral genus one
model with the same discriminant as~{\rm(\ref{weqnstar})}.
\end{Theorem}

This theorem states that, in the $K$-equivalence class of genus one
models representing~$(E,[D])$, there is one which is integral and has
the same discriminant as any given integral Weierstrass model for~$E$.
Our strategy for proving this starts with two observations.

Firstly, the claim really does only depend on the divisor class~$[D]$ and
not the given specific divisor~$D$ in that class, since the
$K$-equivalence class of genus one models representing~$(E,[D])$ only
depends on the divisor class.

Secondly, if $\tau_Q : E \to E$ is translation by some point~$Q \in
E(K)$, then the pairs $(E,[D])$ and $(E,[\tau_Q^*D])$ determine
$K$-equivalent genus one models. This follows from the fact that the
map~$E \to \PP^{n-1}$ determined by $[\tau_Q^*D]$ is the composite
of~$\tau_Q$ and the map determined by~$[D]$.

Using the classical facts that every $K$-rational divisor $D$ of
degree $n$ is linearly equivalent to a unique divisor of the
form~$(n-1).\oh + P$ for some $P \in E(K)$, and that divisors on an
elliptic curve are linearly equivalent if and only if they have the
same degree and the same sum, it suffices to prove
Theorem~\ref{thmstar} for such divisors as $P$ runs over a set of
coset representatives for $E(K)/nE(K)$.

In Lemmas~\ref{lem1} and~\ref{lem2} below,
we show by means of explicit formulae that
Theorem~\ref{thmstar} holds in the cases $D = n. \oh$ and $D =
(n-1).\oh + P$ where $P \in E(K)$ is an integral point, that is, a
point with coordinates in $\OK$.  This is already enough to prove
Theorem~\ref{thmstar} in the case $\OK$ is a complete discrete
valuation ring with residue characteristic prime to $n$. Indeed, by
the theory of formal groups, every non-zero element of $E(K)/nE(K)$
may then be represented by an integral point.

In general we rely on the following two lemmas, 
proved later in this section.
\begin{Lemma}[Unprojection lemma]
\label{lem:unproj}
Let $D \in \Div_K(E)$ have degree $2$ or $3$, and let $P \in E(K)$.
If Theorem~{\rm\ref{thmstar}} holds for $D$ then it holds for $D+P$. 
\end{Lemma}
\begin{Lemma}[Projection lemma]
\label{lem:proj}
Let $D \in \Div_K(E)$ have degree $3$ or $4$, and let $P \in E(K)$.
If Theorem~{\rm\ref{thmstar}} holds for $D$ then it holds for $D-P$. 
\end{Lemma}
Theorem~\ref{thmstar} may be deduced from these 
lemmas in more than one way. For example, if $n=3$ or $4$ 
then $D \sim (n-1).\oh +P$ for some $P \in E(K)$. Then we quote 
the result for $D' = (n-1). \oh$ and use the unprojection lemma.
Likewise if $n=2$ or $3$ then $D \sim (n+1).\oh - P$ for some
$P \in E(K)$. Then we quote the result for $D' = (n+1). \oh$ 
and apply the projection lemma to~$D'$.

\medskip

Theorem~\ref{thmstar} in the case $D= n. \oh$ follows from the
formulae we used to normalise the invariants $c_4$, $c_6$ 
and $\Delta$: see Remark~\ref{rem1}.
\begin{Lemma}
\label{lem1}
Let $E$ be an elliptic curve with Weierstrass equation
\begin{equation}
\label{weqnA}
Y^2 + a_1 XY + a_3 Y = X^3 + a_2 X^2 + a_4 X + a_6. 
\end{equation}
Then the pair $(E,[n.\oh])$ determines genus one models as follows:
$$ \begin{array}{l@{\qquad}l}
n= 2 : & y^2 + (a_1 x_1 x_2 + a_3 x_2^2)y = x_1^3 x_2 + a_2 x_1^2 x_2^2 
  + a_4 x_1 x_2^3 + a_6 x_2^4; \\
n= 3: & y^2z + a_1 xyz + a_3 yz^2 - x^3 - a_2 x^2z - a_4 xz^2 - a_6z^3 = 0; \\
n= 4: & \left\{ \begin{aligned} x^2 - zt & = 0 \\ 
y^2 + a_1 xy + a_3 yz - x t - a_2 x^2 - a_4 xz - a_6z^2 & = 0
  \end{aligned}  \right\}. 
\end{array} $$
Moreover, each of these models has the same invariants $c_4$, $c_6$ and 
$\Delta$ as~{\rm(\ref{weqnA})}.
\end{Lemma}
\begin{Proof}
In the case $n=2$ we embed $E$ in $\PP(1,1,2)$ via $(x_1:x_2:y) =
(X:1:Y)$.  In the cases $n=3,4$ we embed $E$ in $\PP^{n-1}$ via
$(z:x:y) = (1: X: Y)$ and $(z:x:y:t) = (1: X : Y : X^2)$ respectively.
The statement about the invariants follows by direct calculation.
\end{Proof}

Next we prove Theorem~\ref{thmstar} in the case $D= (n-1). \oh + P$
where $P \in E(K)$ is an integral point. By a substitution $X
\leftarrow X+X(P)$, $Y \leftarrow Y + Y(P)$ we may assume that $P$ is
the point $(0,0)$.

\begin{Lemma}
\label{lem2}
Let $E$ be an elliptic curve with Weierstrass equation
\begin{equation}
\label{weqnB}
Y^2 + a_1 XY + a_3 Y = X^3 + a_2 X^2 + a_4 X
\end{equation}
and let $P = (0,0)$.  
Then the pair $(E,[(n-1).\oh+P])$ determines genus one models as follows:
$$ \begin{array}{l@{\qquad}l}
n= 2 : & y^2 + (-x_1^2 + a_1 x_1 x_2 + a_2 x_2^2)y 
   = - a_3 x_1x_2^3 - a_4 x_2^4; \\
n= 3: & y^2z - x^2 y + a_1 xyz + a_2 yz^2 + a_3 x z^2 + a_4 z^3 = 0; \\
n= 4: & \left\{ \begin{aligned} zt - xy + a_1 yz + a_3 z^2 & = 0 \\ 
y^2 - xt + a_2 y z + a_4 z^2 & = 0
  \end{aligned}  \right\}. 
\end{array} $$
Moreover, each of these models has the same invariants $c_4$, $c_6$ and 
$\Delta$ as~{\rm(\ref{weqnB})}.
\end{Lemma}
\begin{Proof}
The rational function 
$$ F = \frac{ Y + a_1 X + a_3 }{X} = \frac{X^2 + a_2 X + a_4}{Y} $$
belongs to the Riemann-Roch space $\LL(\oh + P)$. 
In the case $n=2$ we embed $E$ in $\PP(1,1,2)$ via $(x_1:x_2:y) = (F:1:X)$.
In the cases $n=3,4$ we embed $E$ in $\PP^{n-1}$ 
via $(z:x:y) =  (1:F:X)$ and $(z:x:y:t) = (1:F:X:Y)$  respectively.
The statement about the invariants follows by direct calculation.
\end{Proof}

It remains to prove Lemmas~\ref{lem:unproj} and~\ref{lem:proj}.
One observation that we use in the proofs is the following.

\begin{Lemma}
\label{trans}
The group $\SL_n(\OK)$ acts transitively on $\PP^{n-1}(K)$.
\end{Lemma}
\begin{Proof} Since $\OK$ is a principal ideal domain this is standard. 
See for example \cite[Exercise 6 on p.186]{Jacobson}.
\end{Proof} 

The following lemma explains how to pass between results for generalised 
binary quartics (case $n=2$) and ternary cubics (case $n=3$).

\begin{Lemma}
\label{lem:2<->3}
Let $D \in \Div_K(E)$ be a divisor of degree $2$ and let 
$P \in E(K)$. Let $f_1, f_2, f_3$ be binary forms over 
$K$ with $\deg(f_i) = i$. The following statements are equivalent.
\begin{enumerate}
\item The pair $(E,[D])$ is represented by the generalised binary quartic 
\begin{equation}
\label{eqnC2}
y^2 + f_2(x_1, x_2) y = f_1(x_1,x_2) f_3(x_1,x_2) 
\end{equation}
and $P$ is the point defined by $f_1 = y = 0$.
\item The pair $(E,[D+P])$ is represented by the ternary cubic
\begin{equation}
\label{eqnC3}
f_1(X,Z) Y^2 - f_2(X,Z) Y - f_3(X,Z) = 0 
\end{equation}
and $P$ is the point $(X:Y:Z) = (0:1:0)$. 
\end{enumerate}
\end{Lemma}

\begin{Proof}
We first show that the curves $C_2$ and $C_3$ defined by~(\ref{eqnC2})
and~(\ref{eqnC3}) are isomorphic. An isomorphism $\phi: C_2 \to C_3$
is given by
\begin{align*}
\phi: (x_1:x_2:y) \mapsto (X:Y:Z) & = (x_1 f_1(x_1,x_2) : y + f_2(x_1,x_2) :
x_2 f_1(x_1,x_2) ) \\ & = (x_1 y : f_3(x_1,x_2) : x_2 y)
\end{align*}
with inverse
\[ \phi^{-1} : (X:Y:Z) \mapsto (x_1:x_2:y) 
= (X : Z : f_1(X,Z) Y - f_2(X,Z)). \]
The isomorphism identifies the points $\{ f_1 = y = 0 \} \in C_2(K)$
and $(0:1:0) \in C_3(K)$. To prove the equivalence of (i) and (ii)
we note that if $D = P_1 + P_2$ is a fibre of 
the map $C_2 \to \PP^1 \, ; \, (x_1: x_2:y) \mapsto (x_1:x_2)$ 
then the points $\phi(P_1), \phi(P_2)$ and $(0:1:0)$ are collinear
on $C_3 \subset \PP^2$. 
\end{Proof}

There is an entirely analogous result for passing between
ternary cubics (case $n=3$) and quadric intersections (case $n=4$).

\begin{Lemma}
\label{lem:3<->4}
Let $D \in \Div_K(E)$ be a divisor of degree $3$ and let 
$P \in E(K)$. Let $\ell_1, \ell_2, q_1, q_2$ be ternary forms over
$K$ with $\deg(\ell_i)=1$ and $\deg(q_i) = 2$.
The following statements are equivalent.
\begin{enumerate}
\item The pair $(E,[D])$ is represented by the ternary cubic
\begin{equation}
\label{eqnC3a}
\ell_1(x_1,x_2,x_3) q_2 (x_1, x_2,x_3) 
   - \ell_2(x_1,x_2,x_3) q_1 (x_1, x_2, x_3) = 0
\end{equation}
and $P$ is the point defined by $\ell_1 = \ell_2 = 0$.
\item The pair $(E,[D+P])$ is represented by the quadric intersection
\begin{equation}
\label{eqnC4}
\begin{aligned}
\ell_1(x_1,x_2,x_3) x_4 + q_1 (x_1, x_2,x_3) &= 0 \\
\ell_2(x_1,x_2,x_3) x_4 + q_2 (x_1, x_2,x_3) &= 0 
\end{aligned}
\end{equation}
and $P$ is the point $(x_1: x_2 : x_3:x_4) = (0:0:0:1)$.
\end{enumerate}
\end{Lemma}

\begin{Proof}
We first show that the curves $C_3$ and $C_4$ defined by~(\ref{eqnC3a})
and~(\ref{eqnC4}) are isomorphic. An isomorphism $\phi: C_3 \to C_4$
is given by
\[ \phi: (x_1:x_2:x_3) \mapsto  
(x_1 \ell_1 : x_2 \ell_1 : x_3 \ell_1 : - q_1 )
 = (x_1 \ell_2 : x_2 \ell_2 : x_3 \ell_2 : - q_2 ) \]
with inverse
\[ \phi^{-1} : (x_1:x_2:x_3:x_4) \mapsto (x_1:x_2:x_3). \]
This isomorphism identifies the points $\{\ell_1 =\ell_2 = 0\} \in C_3(K)$
and $(0:0:0:1) \in C_4(K)$. To prove the equivalence of (i) and (ii)
we note that if $C_3 \subset \PP^2$ meets some line in the 
divisor $D = P_1 + P_2 +P_3$ then the points $\phi(P_1), \phi(P_2), \phi(P_3)$ 
and $(0:0:0:1)$ are coplanar on $C_4 \subset \PP^3$. 
\end{Proof}

A generic computation shows that the genus one models~(\ref{eqnC2})
and~(\ref{eqnC3}) in Lemma~\ref{lem:2<->3} have the same discriminant.
Likewise the models~(\ref{eqnC3a}) and~(\ref{eqnC4}) in Lemma~\ref{lem:3<->4}
have the same discriminant.

\medskip

\begin{ProofOf}{Lemma~\ref{lem:unproj}} 
(i) Let $D \in \Div_K(E)$ be a divisor of degree 2, and suppose
the pair $(E,[D])$ is represented by an integral generalised
binary quartic of discriminant $\Delta$. By Lemma~\ref{trans}
(with $n=2$) we may assume that $P$ is the point 
$(x_1: x_2 :y )= (1:0:\eta)$ for some $\eta \in K$. 
Since $\OK$ is integrally closed it follows that $\eta \in \OK$. 
By making a substitution $y \leftarrow y + \eta x_1^2$ we may assume 
that $\eta = 0$. Our model is now of the form~(\ref{eqnC2}) with
$f_1(x_1,x_2) = x_2$.
Then the ternary cubic~(\ref{eqnC3}) is an integral model of
discriminant $\Delta$ representing the pair $(E,[D+P])$. 

\noindent
(ii) Let $D \in \Div_K(E)$ be a divisor of degree 3, and suppose
the pair $(E,[D])$ is represented by an integral ternary cubic
of discriminant $\Delta$. By Lemma~\ref{trans} (with $n=3$) 
we may assume that $P$ is the point $(x_1: x_2 :x_3 )= (0:0:1)$. 
Our model is now of the form~(\ref{eqnC3a}) with $\ell_1 =x_1$
and $\ell_2 = x_2$. We may choose the quadratic forms $q_1$ and $q_2$ 
to have coefficients in $\OK$.
Then the quadric intersection~(\ref{eqnC4}) is an integral model of
discriminant $\Delta$ representing the pair $(E,[D+P])$. 
\end{ProofOf}

\begin{ProofOf}{Lemma~\ref{lem:proj}} 
(i) Let $D \in \Div_K(E)$ be a divisor of degree 3, and suppose
the pair $(E,[D])$ is represented by an integral ternary cubic 
of discriminant $\Delta$. By Lemma~\ref{trans}
(with $n=3$) we may assume that $P$ is the point 
$(x_1: x_2 :x_3 )= (0:0:1)$. Our model is now of the form~(\ref{eqnC3}).
Then the generalised binary quartic~(\ref{eqnC2}) is an integral model of
discriminant $\Delta$ representing the pair $(E,[D-P])$. 

\noindent
(ii) Let $D \in \Div_K(E)$ be a divisor of degree 4, and suppose
the pair $(E,[D])$ is represented by an integral quadric intersection
of discriminant $\Delta$. By Lemma~\ref{trans} (with $n=4$) 
we may assume that $P$ is the point $(x_1: x_2 :x_3 :x_4)= (0:0:0:1)$
Our model is now of the form~(\ref{eqnC4}) for some forms
$\ell_1,\ell_2, q_1, q_2$ with coefficients in $\OK$.
Then the ternary cubic~(\ref{eqnC3}) is an integral model of
discriminant $\Delta$ representing the pair $(E,[D-P])$. 
\end{ProofOf}

\begin{Remark} 
In principle these proofs give an algorithm for
minimising $K$-soluble models, but only once a $K$-rational 
point is explicitly known. Although it is easy to decide solubility 
over local fields, such an algorithm would require that we find a local point to
sufficiently high precision.
Hence our comment that this is not a readily implementable algorithm.
\end{Remark}


\section{Minimisation algorithms}
\label{sec:M}

In this section we give algorithms for minimising binary quartics
(case $n=2$), ternary cubics (case $n=3$) and quadric intersections
(case $n=4$). As in Section~\ref{min:statements}
we work over a field $K$ which is the field of
fractions of a discrete valuation ring $\OK$. 
There is no need to assume that $K$ is complete (or even Henselian). 
We fix a uniformiser
$\pi$ and write $k= \OK/\pi \OK$ for the residue field.  
In the cases $n=2,4$ we initially assume that $\Char(k) \not=2$, 
leaving the case $\Char(k) =2$ to Section~\ref{min_char2}.

Our algorithms for $n=2,3$ share some common features which
we now elucidate. In these cases we specify a procedure that takes 
as input an integral genus one model of positive level, and 
returns a $K$-equivalent integral model of the same or smaller level. 
We then show that if the model is non-minimal then the level must decrease
after finitely many iterations, and give a bound $N$ on the number of
iterations required.  This also gives a test for minimality: if $N$
iterations of the procedure fail to decrease the level, then the model
must be minimal.

The proofs are by induction on the {\em slope}, which we define as the
least valuation of the determinant of a matrix $M \in \GL_n(K)$ 
with entries in $\OK$ that can be used to decrease the level. 
The slope of a minimal model is undefined.
The arguments we use are incapable of proving the Minimisation
Theorem, since we assume at the outset that the given model has
a slope, {\em i.e.} is non-minimal.

The following lemma is used to show that our procedure
gives a well-defined map on $\OK$-equivalence
classes. This is useful, since it means we are free to replace our 
model by an $\OK$-equivalent one at any stage of the proof.
We write $I_m$ for the $m$ by $m$ identity matrix.

\begin{Lemma}
\label{welldef}
Let $\GL_n(K)$ act on $\PP^{n-1}$ in the natural way (via left
multiplication of column vectors by matrices).  Let $\alpha =
\Diag(I_r, \pi I_{n-r})$ for some $0 < r < n$.  Then the subgroup of
$\GL_n(\OK)$ consisting of transformations whose reduction mod $\pi$
preserves the subspace $\{ x_{r+1} = \ldots = x_n = 0 \}$ is
$$ \GL_n(\OK) \cap \alpha \GL_n(\OK) \alpha^{-1} $$
\end{Lemma}
\begin{Proof} Identifying $\PP^{n-1}(K)$ with the non-zero elements
  of~$K^n$ modular scalars, $\GL_n(\OK)$ is the subgroup preserving
  $\mathcal{O}_K^n$ and we are interested in the subgroup which also
  preserves $\mathcal{O}_K^r\oplus (\pi\OK)^{n-r} =
  \alpha(\mathcal{O}_K^n)$.  The statement is now clear.
\end{Proof}

This lemma is used as follows. Suppose that $\Phi$ and $\Psi$ 
are $\GL_n(\OK)$-equivalent models, and the matrix relating them 
is one whose reduction mod $\pi$ preserves the subspace
$\{ x_{r+1} = \ldots = x_n = 0 \}$. Then the models
$\Phi'$ and $\Psi'$ obtained by applying $\alpha= \Diag(I_r, \pi I_{n-r})$ 
to both $\Phi$ and $\Psi$, will again be $\GL_n(\OK)$-equivalent.

\subsection{Minimisation of $2$-coverings} 
\label{sec:minbq}

Let $F \in K[x,z]$ be a binary quartic, say
$$F (x,z) = a x^4 + b x^3 z + c x^2 z^2 + d x z^3 + e z^4.$$ Viewing
the set of these as a subset of~$X_2(K)$, the group
of~$K$-equivalences between binary quartics
is~$K^{\times}\times\GL_2(K)$, where $[\mu,M]$ acts 
as $[\mu,(0,0,0),M] \in \G_2(K)$. 
Note that $[\pi^{-2},\pi I_2]$ acts
trivially, so we may if convenient assume that~$M$ has entries
in~$\OK$, not all in~$\pi\OK$.

We say that an integral binary quartic $F$ is {\em minimal} if $v(\Delta(F))$ 
is minimal among all integral binary quartics $K$-equivalent to $F$.
If $\Char(k) = 2$ then this need not be the same as being minimal as
a generalised binary quartic.
We define the valuation $v(F)$ to be the minimum of the valuations of
the coefficients. If $v(F) \ge 2$, then $F$ is not 
minimal, and indeed dividing through by~$\pi^2$ gives a 
$K$-equivalent integral model of smaller level.
The algorithm for minimising binary quartics is described
in the following theorem.

\begin{Theorem}
\label{minbq}
Let $F \in \OK[x,z]$ be a non-singular
binary quartic. Suppose that $v(F) = 0$ or $1$, but $F$ has positive level.
If $\Char(k)=2$ then further assume that $F$ is non-minimal.
Then
\begin{enumerate}
\item The reduction mod $\pi$ of $F_1(x,z) = \pi^{-v(F)} F(x,z)$
has either a triple or quadruple root defined over $k$.
\item The following procedure replaces $F$ by a $K$-equivalent integral model 
of the same level. 
\begin{itemize}
\item Move the repeated root
of $F_1(x,z)$ mod $\pi$ to $(x:z) = (0:1)$. 
\item Replace $F(x,z)$ by $\pi^{-2} F(\pi x,z)$.
\end{itemize}
\item If $F$ is non-minimal then the procedure in (ii) gives $v(F) \ge 2$
after at most $2$ iterations.
\end{enumerate}
\end{Theorem}

\begin{Proof}
We first prove the theorem in the case $F$ is non-minimal.
By hypothesis there exists $[\mu,M] \in K^{\times}\times\GL_2(K)$ with
$v(\mu\det(M))\le-1$ such that
the transform of~$F$ by~$[\mu,M]$ is still integral.
The slope $s$ of~$F$ is 
the least possible valuation of $\det M$, for $M$ such a
matrix with entries in $\OK$. By Lemma~\ref{welldef} we are free to
replace $F$ by any $\OK$-equivalent binary quartic.  So, putting 
$M$ in Smith normal form, we may assume that
\begin{equation*}
F(\pi^sx,z) \equiv 0 \pmod{\pi^{2s+2}} 
\end{equation*}
where $s$ is the slope.  For $s \ge 2$, this condition works out as
$\pi^2 \dv c$, $\pi^{s+2} \dv d$ and $\pi^{2s+2} \dv e$. So the only
possible slopes are $s=0,1,2$ (as if these conditions hold for
some~$s>2$, then they also hold for $s=2$, and $s$ was defined to be
minimal). If $s=0$, then $v(F) \ge 2$ contrary to hypothesis.  If
$s=1$, then the coefficients of $F$ have valuations satisfying
\begin{equation*}
\ge 0 \quad \ge 1 \quad \ge 2 \quad \ge 3 \quad \ge 4. 
\end{equation*}
So either $v(F)=0$ and $F(x,z)$ mod $\pi$ has a quadruple root
at $(x:z)=(0:1)$, or $v(F)=1$ and $\pi^{-1}F(x,z)$ mod $\pi$ has
a triple or quadruple root at $(x:z)=(0:1)$. 
If $s=2$, then the coefficients of $F$ have valuations satisfying
\begin{equation*}
\ge 0 \quad = 0 \quad \ge 2 \quad \ge 4 \quad \ge 6. 
\end{equation*}
Then $F(x,z)$ mod $\pi$ has a triple root at $(0:1)$.
In each of the cases $s=1,2$ statements (i) and (ii) of the 
theorem are now clear. Moreover the procedure in (ii)
returns a $K$-equivalent integral model of smaller slope. 
Hence at most 2 iterations are required to give 
$v(F) \ge 2$, establishing (iii).

It remains to prove (i) and (ii) in the case
$\Char(k) \not=2$ and $F$ has positive level (but could be minimal).
Statement (i) follows from the fact that $F_1$ mod $\pi$ is a null form,
{\em i.e.} both the invariants $I$ and $J$ vanish.
(Since $k$ is perfect the multiple root is defined over $k$.)
For (ii) we must show that
if $v(F)=0$ and the reduction of $F$ mod $\pi$ has a repeated root at 
$(x:z)=(0:1)$ then $\pi^2 \dv e$. But
in this case there are smooth $\kbar$-points on the
reduction of $\CC$ mod $\pi$ where $\CC = \{y^2 = F(x,z) \}$. 
So after an unramified extension we may assume 
that $\CC(K) \not= \emptyset$. Then Theorem~\ref{minthm} shows 
that $F$ is non-minimal, and our earlier argument applies.
\end{Proof}

To give a satisfactory analogue of this algorithm when $\Char(k)=2$
we must work with generalised binary quartics.  
We give details in~Section~\ref{min_char2}.


\subsection{Minimisation of 3-coverings} 
\label{M3}
The valuation $v(F)$ of a ternary cubic
\[  F(x,y,z) = a x^3 + b y^3 + c z^3 + a_2 x^2 y + a_3 x^2 z 
  + b_1 x y^2 + b_3 y^2 z + c_1 x z^2  + c_2 y z^2 + m x y z \] 
is the minimum valuation of a coefficient. If $v(F) \ge 1$ then $F$ 
is non-minimal, and indeed dividing through by $\pi$ gives a $K$-equivalent
integral model of smaller level. The algorithm for minimising ternary 
cubics is described in the following theorem.

\begin{Theorem} 
\label{mintc}
Let $F \in X_3(\OK)$ be a non-singular ternary cubic. Suppose $v(F)=0$, 
but $F$ has positive level. Then
\begin{enumerate}
\item The singular locus of the reduction 
\begin{equation*}
\label{singloc3}
{\mathcal S} = 
\{ (x:y:z) \in \PP^2 \mid F \equiv \tfrac{\partial F}{\partial x}
  \equiv \tfrac{\partial F}{\partial y} \equiv \tfrac{\partial F}{\partial z}
  \equiv 0 \pmod{\pi}  \} \end{equation*}
is either a point or a line, and is defined over $k$.
\item The following procedure replaces $F$ by a $K$-equivalent integral 
ternary cubic of the same level.
\begin{itemize}
\item Make a $\GL_3(\OK)$-transformation to move the singular 
locus ${\mathcal S}$ to the point $(1:0:0)$, respectively 
the line $\{z=0\}$.
\item  Replace $F(x,y,z)$ by $\pi F(\pi^{-1}x,y,z)$, 
respectively $\pi^{-1} F( x,y,\pi z)$.
\end{itemize}
\item If $F$ is non-minimal then the procedure in (ii) gives $v(F) \ge 1$
after at most $4$ iterations.
\end{enumerate}
\end{Theorem} 

\begin{Proof} 
We are given that $F$ has positive level. It follows that its
reduction mod $\pi$ is a null-form, {\em i.e.} the invariants 
$c_4$, $c_6$ and $\Delta$ all vanish.
The classification of singular ternary cubics (up to equivalence
over an algebraically closed field) is well known. See for 
example \cite[\S 10.3]{Dolg} or \cite{PoonenCubics}. 
The possible null-forms are either a cuspidal cubic, 
a line touching a conic, three lines through a common point,
a double line and a line, or a triple line. So over $\kbar$ 
the singular locus of the reduction is either a point or a line.
Since $k$ is perfect, this point or line is already defined over $k$.
This proves (i).

Next we prove (ii) and (iii) in the case 
$F$ is non-minimal. By hypothesis there exists $[\mu,M] \in \G_3(K) 
= K^\times \times \GL_3(K)$ with $v(\mu \det M) \le -1$ such that 
the transform of $F$ by $[\mu,M]$ is still integral. Since 
$[\pi^{-3},\pi I_3]$ acts trivially, we may assume that~$M$ has entries
in~$\OK$. The slope $s$ of $F$ is the least possible valuation 
of $\det M$, for $M$ such a matrix with entries in $\OK$. 
By Lemma~\ref{welldef} we are free to replace $F$ by 
any $\OK$-equivalent ternary cubic. So, putting $M$ in Smith normal form, 
we may assume that
\begin{equation}
\label{admiss}
 F(x,\pi^a y, \pi^b z) \equiv 0 \pmod{\pi^{a+b+1}} 
\end{equation} 
for some $0 \le a \le b$ with $a+b=s$. If $a=b=0$, then $v(F) \ge 1$,
contrary to hypothesis. If $a=0$ and $b \ge 1$, then 
the reduction of $F$ mod $\pi$ only involves the
monomials $x z^2$, $y z^2$ and $z^3$. Hence 
${\mathcal S}$ is the line $\{ z= 0\}$. If $a \ge 1$, then
the coefficients of $x^3$, $x^2y$ and $x^2 z$ all vanish 
mod $\pi$. Hence ${\mathcal S}$ is either the point $(1:0:0)$ 
or a line through this point. In each of these
cases it is clear that the procedure in (ii) returns an integral model 
of the same level and smaller slope. Moreover 
it gives $v(F) \ge 1$ after a finite number of 
iterations (bounded by the initial slope).
The next lemma shows that the only possible slopes are
$0,1,2,3$ and $5$. Hence at most $4$ iterations 
are required, establishing (iii).

It remains to prove (ii) in the case $F$ has positive level (but could 
be minimal). We must show that if $(1:0:0)$ is the only singular
point on the reduction then $F(1,0,0) \equiv 0 \pmod{\pi^2}$. But
in this case there are smooth $\kbar$-points on the
reduction. So after an unramified extension we may assume 
that $\CC_F(K) \not= \emptyset$. Then Theorem~\ref{minthm} shows 
that $F$ is non-minimal, and our earlier argument applies.
\end{Proof}

We say that a pair $(a,b)$ is {\em admissible} for $F$ 
if~(\ref{admiss}) holds.

\begin{Lemma}
If some pair $(a,b)$ with $0 \le a \le b$ is admissible for $F$ 
then at least  one of the pairs $(0,0)$, $(0,1)$, $(1,1)$, 
$(1,2)$ or $(2,3)$ is admissible for $F$.
\end{Lemma}
\begin{Proof}
Suppose $(a,b)$ is admissible for $F$. We make the observations:
\begin{itemize}
\item If $a = 0$ and $b \ge 1$ then $(0,1)$ is admissible. 
\item If $a=b \ge 1$ then $(1,1)$ is admissible.
\item If $a \ge 1$ and $b \ge 2a $ then $(1,2)$ is admissible.
\item If $a \ge 2$ and $b \ge a+1$ then $(2,3)$ is admissible.
\end{itemize}
The only remaining possibility is $(a,b)=(0,0)$.
\end{Proof}

\begin{Example}
We apply our algorithm to a cuspidal cubic. (Although this is singular,
there are $\pi$-adically close smooth ternary cubics that
are treated in the same way by our algorithm.) 
An arrow labelled $(0,a,b)$ 
indicates that we make the transformation
$[\pi^{-a-b},\Diag(1,\pi^a,\pi^b)]$.
$$ \begin{array}{rcl}
x z^2 - y^3
& \stackrel{(0,1,1)}{\ra} & x z^2 - \pi y^3 \\ 
& \stackrel{(0,0,1)}{\ra} & \pi x z^2 - y^3 \\
& \stackrel{(0,1,0)}{\ra} & x z^2 - \pi^2 y^3 \\
& \stackrel{(0,0,1)}{\ra} & \pi(x z^2 - y^3) 
\end{array} $$
So this is an example where our algorithm takes the maximum 
possible of $4$ iterations to give $v(F) \ge 1$.
\end{Example}


\subsection{Minimisation of 4-coverings} 
\label{M4}
In this section we prove 
Theorems~\ref{mainthmA}(i) and~\ref{mainthmB1} in the case
$n=4$, assuming that $\Char(k)\not=2$.  The proofs are constructive
and give an algorithm for minimising quadric intersections.
The modifications required when $\Char(k)=2$ are described in the 
next section.

We define a map
\begin{equation}
\begin{aligned}
\dd: X_4(K) & \to X_2(K) \\
   (Q_1,Q_2) & \mapsto F(x,z) = \det(Ax + Bz) 
\end{aligned}
\end{equation}
where $A$ and $B$ are the matrices of second partial 
derivatives of $Q_1$ and $Q_2$. As noted in Definition~\ref{def:g1m4} we have
$\Delta(Q_1,Q_2) = 2^{-12} \Delta(F)$. 

\begin{Lemma} 
\label{4to2cover}
Let $(Q_1,Q_2) \in X_4(K)$ be a non-singular
quadric intersection. Then $F = \dd (Q_1,Q_2)$ is non-singular,
and there is a morphism of genus one curves 
$\CC_{(Q_1,Q_2)} \to \CC_F$ defined over $K$. 
\end{Lemma}
\begin{Proof}  A formula for this morphism is given by classical
invariant theory, as we now recall from \cite{AKM3P}, \cite{MSS}.
We write the binary quartic $F = \dd(Q_1,Q_2)$ as
$F(x,z) = a x^4 + b x^3z + c x^2z^2 + d x z^3+ e z^4$, and let
$T_1$ and $T_2$ be the quadrics whose matrices of second partial
derivatives $M_1$ and $M_2$ are determined by
\begin{equation}
\label{defnT1T2}
  \adj\bigl(\adj(A) x + \adj(B) z \bigr)
     = a^2 A x^3 + a M_1 x^2 z + e M_2 x z^2 + e^2 B z^3 \,.
\end{equation}
Then $J^2 \equiv F(T_1,-T_2) \mod (Q_1,Q_2)$ where
$J = \frac{1}{4} 
\frac{\partial(Q_1,Q_2,T_1,T_2)}{\partial(x_1,x_2,x_3,x_4)}$.
\end{Proof}

\begin{Lemma} 
\label{commdiag}
If $[M,N] \in \G_4(K)$ then there is a commutative
diagram
\[ \xymatrix@+1.2em{ 
   X_4(K)  \ar[r]^{[M,N]} \ar[d]_{\dd} 
    & X_4(K) \ar[d]^{\dd} \\
   X_2(K)  \ar[r]^{[\det N,M]} & X_2(K). 
}
\]
In particular $\dd$ induces a well-defined map on
$K$-equivalence classes.
\end{Lemma}
\begin{Proof} This is clear.
\end{Proof}

Following the treatment in Womack's thesis \cite{WomackThesis}, we deduce 
the Minimisation Theorem for $n=4$ from the Minimisation 
Theorem for $n=2$. The modifications required to prove 
Theorems~\ref{mainthmA}(i) and~\ref{mainthmB1} are given at the end of 
this section (see Proposition~\ref{strong4} below).

\begin{Proposition}
\label{weakmin}
If $(Q_1,Q_2) \in X_4(K)$ is non-singular and $K$-soluble 
then it is $K$-equivalent to an integral model of level $0$.
\end{Proposition}
\begin{Proof}
Since $(Q_1,Q_2)$ is $K$-soluble, it follows by Lemma~\ref{4to2cover}
that $\dd(Q_1,Q_2)$ is $K$-soluble.
So by the minimisation theorem for $n=2$ we know that
$\dd(Q_1,Q_2)$ is $K$-equivalent to an integral
binary quartic $F(x,z)$ of level $0$. It is clear by 
Lemma~\ref{commdiag} that $(Q_1,Q_2)$ is $K$-equivalent to a 
quadric intersection $(Q'_1,Q'_2)$ with $\dd(Q'_1,Q'_2) = F$. 
The following lemma shows we may take $(Q'_1,Q'_2)$ integral. This is
then the required integral model of level $0$.
\end{Proof}
Notice that the next three lemmas are false when $\Char(k)=2$, 
as we could otherwise use the above proof 
to find integral models of level $-v(2)$.
\begin{Lemma}
\label{iterateMRL}
Let $(Q_1,Q_2) \in X_4(K)$ be a $K$-soluble non-singular 
quadric intersection. If $\dd(Q_1,Q_2)$ 
is integral then $(Q_1,Q_2)$ is $K$-equivalent to an 
integral quadric intersection 
$(Q'_1,Q'_2)$ with $\dd(Q'_1,Q'_2) = \dd(Q_1,Q_2)$.
\end{Lemma}
\begin{Proof}
By a transformation $[\mu I_2,I_4]$ for suitable $\mu \in \OK$
we obtain an integral quadric intersection $(Q'_1,Q'_2)$ with 
$\dd(Q'_1,Q'_2) = \mu^4 \dd(Q_1,Q_2)$. 
We now apply the following lemma, as many times as required, at
each stage preserving the integrality of $(Q'_1,Q'_2)$ while
dividing $\dd(Q'_1,Q'_2)$ by a square in $\pi \OK$.
\end{Proof}
Recall that we write $v(F)$ for the minimum of the valuations
of the coefficients of the binary quartic $F$. The following is
Womack's ``main reduction lemma''.
\begin{Lemma} \label{MRL}
Let $(Q_1,Q_2) \in X_4(\OK)$ be a non-singular $K$-soluble
integral quadric intersection. If $F = \dd(Q_1,Q_2)$ 
satisfies $v(F) \ge 2$ then $(Q_1,Q_2)$ is $K$-equivalent to
an integral quadric intersection of smaller level by means
of a transformation $[\lambda I_2, N] \in \G_4(K)$ 
with $\lambda \in K^\times$ and $N \in \GL_4(K)$.
\end{Lemma}

The following geometric lemma prepares for the proof of
Lemma~\ref{MRL}. We say that two pairs of quadratic forms in $m$
variables are $k$-equivalent if they are in the same orbit
for the natural action of $\GL_2(k) \times \GL_m(k)$.
(This extends our earlier definition in the case $m=4$.)
Over an algebraically closed field, the lemma may alternatively
be deduced from the classification of pairs of quadrics 
using the Segre symbol, as given in \cite[Chapter XIII, \S11]{HodgePedoe}.

\begin{Lemma}
\label{geomlem}
Let $Q_1$ and $Q_2$ be quadratic forms in $m = 3$ or $4$ variables over
a field $k$ with $\Char(k) \not=2$. Let $A$ and $B$ be
the matrices of second partial derivatives of $Q_1$ and $Q_2$.
Assume that
\begin{itemize}
\item $\{Q_1 = Q_2 = 0 \} \subset \PP^{m-1}$ is not
a cone, i.e. $\ker(A) \cap \ker(B) = 0$, and 
\item The binary form $F(x,z) = \det(A x + B z)$ is identically zero.
\end{itemize}
Then the $k$-equivalence class of $(Q_1,Q_2)$ is uniquely
determined:
\begin{enumerate}
\item If $m=3$  then $(Q_1,Q_2)$ is $k$-equivalent to $(x_1 x_2, x_2 x_3)$ 
\item If $m=4$ then $(Q_1,Q_2)$ 
is $k$-equivalent to $(x_1 x_2, x_2 x_3 - x_4^2)$.
\end{enumerate}
\end{Lemma}

\begin{Proof}
(i) We must show that the gcd of $Q_1$ and $Q_2$ is a linear form,
and for this we may assume that $k$ is algebraically closed.
Since some quadric in the pencil has rank 2, we may assume 
that $Q_1 = x_1 x_2$. Then the condition $\det(A x + Bz)=0$ works out as 
$b_{33} = b_{13} b_{23} = \det B = 0$. 
Swapping $x_1$ and $x_2$ if necessary, we may assume
that $b_{13} = b_{33} = 0$. Then $b_{23} \not= 0$ (otherwise we would have
a cone) and the condition $\det B =0$ forces $b_{11} = 0$. 
Making a substitution for $x_3$ now puts
$(Q_1,Q_2)$ in the required form.

\noindent
(ii) Suppose  $\{Q_1 = Q_2 = 0\} \subset \PP^3$ has a singular 
point defined over $k$.
Moving this point to $(1:0:0:0)$, it is easy to reduce to the case
\[ A = \left(\begin{smallmatrix} 
\begin{smallmatrix} 0 & 1 \\ 1 & 0 \end{smallmatrix} &
\begin{smallmatrix} 0 & 0 \\ 0 & 0 \end{smallmatrix} \\
\begin{smallmatrix} 0 & 0 \\ 0 & 0 \end{smallmatrix} & A' 
\end{smallmatrix} \right) ,
 \qquad B = \left(\begin{smallmatrix}
\begin{smallmatrix} 0 & 0 \\ 0 & * \end{smallmatrix} & 
\begin{smallmatrix} 0 & 0 \\ * & * \end{smallmatrix} \\
\begin{smallmatrix} 0 & * \\ 0 & * \end{smallmatrix} & 
B' \end{smallmatrix} \right).
\]
The condition 
$\det(A x + B z)=0$ now becomes $\det(A' x + B' z)=0$. Hence we may
assume that $A'$ and $B'$ are scalar multiples of 
$(\begin{smallmatrix} 0 & 0 \\ 0 & 1 \end{smallmatrix})$. Then $b_{23} \not=0$
(otherwise we have a cone) and a substitution in $x_3$ brings us
to the case 
\[ (Q_1, Q_2) = (x_1 x_2 + \lambda x_4^2, x_2 x_3 + \mu x_4^2) \]
for some $\lambda,\mu \in k$. Replacing one of these quadrics
by a suitable linear combination, and then making a substitution in
$x_1$ and $x_3$ to compensate, we may assume that $\lambda =0$. 
Then $\mu \not= 0$ (otherwise we have a cone) and we rescale to
get $\mu = - 1$.

By Theorem~\ref{thm:invjac}(ii) there is a singular point defined
over $\kbar$. So running the above proof over $\kbar$ shows that 
$\{Q_1 = Q_2 = 0\} \subset \PP^3$ is the union of a conic
and a line, meeting at a unique point. This point of intersection
is a $k$-rational singular point. Our earlier proof now applies.
\end{Proof}

\begin{ProofOf}{Lemma~\ref{MRL}} 
We write $\Qbar_1$ and $\Qbar_2$ for the reductions of $Q_1$ and $Q_2$
mod $\pi$. In the proof we often arrive at one of the following 
three special situations.
\begin{description}
\item[Situation 1]
  The reduction $\CC_{(\Qbar_1,\Qbar_2)}$ contains a 
  plane defined over $k$. \\
  By a $\GL_4(\OK)$-transformation we may move the plane to $\{x_1=0\}$. 
  We apply the transformation $[\pi^{-1}I_2,\Diag(\pi,1,1,1)]$ 
  to give an integral model of smaller level.

\item[Situation 2] 
  The reduction  $\CC_{(\Qbar_1,\Qbar_2)}$ is a cone over a point 
  $\x \in \PP^3(k)$ and moreover 
  $Q_1(\x) \equiv Q_2(\x) \equiv 0 \pmod{\pi^2}$. \\
  By a $\GL_4(\OK)$-transformation we may move the point to
  $(1:0:0:0)$. We apply the transformation $[I_2,\Diag(\pi^{-1},1,1,1)]$
  to give an integral model of smaller level.

\item[Situation 3] 
  The reduction  $\CC_{(\Qbar_1,\Qbar_2)}$ contains a line defined over $k$. \\
  By a $\GL_4(\OK)$-transformation we may move the line to
  $\{ x_1 = x_2 = 0\}$. The ``flip-flop'' transformation 
  $[\pi^{-1} I_2,\Diag(\pi,\pi,1,1)]$ gives an integral model of the
  same level.
\end{description}

Let $A$ and $B$ be the matrices of second partial derivatives 
of $Q_1$ and $Q_2$. Let $\Abar$ and $\Bbar$ be their reductions
mod $\pi$. We split into cases according to the value of the
{\em common nullity}, defined as
$s = \dim (\ker \Abar \cap \ker \Bbar )$.

If $s=0$ then by Lemma~\ref{geomlem}(ii) we are in Situation 3. 
Applying the ``flip-flop'' transformation brings us to the
case $s \ge 1$.

If $s=1$ we may assume that $\Qbar_1$ and $\Qbar_2$ are quadratic
forms in $x_2,x_3,x_4$ only. Let $A'$ and $B'$ be the 3 by 3
matrices of second partial derivatives. Then 
\begin{equation}
\label{eq:s=1}
F(x,z) \equiv (a_{11} x + b_{11} z) \det(A'x + B'z) \pmod{\pi^2}. 
\end{equation}
Since $v(F) \ge 2$ we have either 
$a_{11} \equiv b_{11} \equiv 0 \pmod{\pi^2}$ in which case we are
in Situation 2, or $\det(A' x + B' z)=0$ in which case 
Lemma~\ref{geomlem}(i) shows we are in Situation 1.

If $s \ge 2$ we may assume that $\Qbar_1$ and $\Qbar_2$ 
are binary quadratic forms in $x_1$ and $x_2$.
If $\Qbar_1$ and $\Qbar_2$ simultaneously represent 0 over $k$, then 
we are in Situation 1. Otherwise we apply
the ``flip-flop'' transformation $[\pi^{-1} I_2, \Diag(\pi,\pi,1,1)]$
to give an integral model $(R_1,R_2)$ of the same level. 
Then $\Rbar_1$ and $\Rbar_2$ are binary quadratic forms 
in $x_3$ and $x_4$. If $\Rbar_1$ and $\Rbar_2$ simultaneously 
represent 0 over $k$ then we are in Situation 1.
Otherwise we obtain a contradiction to our hypothesis that $(Q_1,Q_2)$ 
is $K$-soluble. Indeed if $(x_1:x_2:x_3:x_4)$ were a $K$-point with 
$\min\{v(x_i): 1 \le i \le 4 \}=0$ then from $Q_1(\x) \equiv Q_2(\x) 
\equiv 0 \pmod{\pi}$ we deduce $x_1 \equiv x_2 \equiv 0 \pmod{\pi}$ 
and from $Q_1(\x) \equiv Q_2(\x) 
\equiv 0 \pmod{\pi^2}$ we deduce $x_3 \equiv x_4 \equiv 0 \pmod{\pi}$.
\end{ProofOf}

This completes the proof of Proposition~\ref{weakmin}.
We now modify the proof so that we can deduce 
Theorems~\ref{mainthmA}(i) and~\ref{mainthmB1} in the case $n=4$
from the corresponding results for $n=2$.
The situation considered at the end of the last paragraph motivates
the definition of a critical model, 
see Definition~\ref{defcrit}(c) below. 

\begin{Proposition}
\label{strong4}
If $(Q_1,Q_2) \in X_4(K)$ is non-singular then it is $K$-equivalent to either
\begin{enumerate}
\item an integral model $\Phi \in X_4(\OK)$ with $\dd (\Phi)$
minimal (and hence $\Phi$ minimal), or
\item a critical model, as specified in Definition~\ref{defcrit}(c) below.
\end{enumerate}
\end{Proposition}
\begin{Proof}
By Lemma~\ref{commdiag} we may assume that $\dd (Q_1,Q_2)$
is a minimal binary quartic. We then follow the proof of 
Lemma~\ref{iterateMRL}, but without the hypothesis of $K$-solubility.
This hypothesis was only used at 
the end of the proof of Lemma~\ref{MRL}.
We may assume that one of the pairs, say 
$\Qbar_1$ and $\Qbar_2$, simultaneously represents 0 over $\kbar$. 
(Otherwise we would have a critical model.)
If they do not simultaneously represent 0 over $k$, then they
must be linearly dependent. So it is clear we can reduce the level, 
but not necessarily using a transformation of the specified form.
In the proof of Lemma~\ref{iterateMRL} 
we repeatedly applied Lemma~\ref{MRL}.
For the final application it does not matter what transformation we use.
In all earlier applications we have $v(F) \ge 3$. 
If $A_1, B_1$ and $A_2, B_2$ are the 2 by 2 matrices representing 
the pairs of binary quadratic forms $\Qbar_1, \Qbar_2$ 
and $\Rbar_1, \Rbar_2$ then
\[ F(x,z) \equiv \pi^2 \det(A_1 x + B_1 z) \det(A_2 x + B_2 z) 
\pmod{\pi^3}. \]
The hypothesis $v(F) \ge 3$ therefore ensures that one of the pairs 
simultaneously represents 0 over $k$. We are then in Situation 1.
\end{Proof}

In Lemma~\ref{critprops2} (see below) we show that critical models 
are minimal. Hence the proof of Proposition~\ref{strong4} gives an 
algorithm for minimising quadric intersections, even in the 
case they are not $K$-soluble.
Proposition~\ref{strong4} also allows us to deduce the case $n=4$
of Theorems~\ref{mainthmA}(i) and~\ref{mainthmB1} from the
case $n=2$. Here we use the easy facts that critical models 
are $K^\sh$-insoluble, and remain critical after any unramified
field extension.

\subsection{Minimisation in residue characteristic 2}
\label{min_char2}
We describe how to modify our algorithms in 
the cases $n=2,4$ when $\Char(k)=2$.
In the case $n=2$ the issue is that we must work with generalised
binary quartics instead of just binary quartics.
Recall that a generalised binary quartic, 
or genus one model of degree~2, 
is an equation of the form $$y^2 + P(x,z) y = Q(x,z)$$
where $P$ and $Q$
are homogeneous polynomials of degrees 2 and 4.
We label the coefficients of $P$ and $Q$ as $l,m,n$ and $a,b,c,d,e$.
We observe that in characteristic~$2$ the binary quadratic form 
$\partial^2 Q / \partial x \partial z = bx^2 + dz^2$ is a covariant of 
the quartic~$Q$. Moreover this covariant vanishes if and only if $Q$ is 
a square. (Recall that $k$ is perfect, and so every element of $k$ 
is a square.)

We say that two models are {\em $y$-equivalent} if they are
related by a {\em $y$-substitution}, that is, a substitution of the
form $x \leftarrow x$, $z \leftarrow z$, 
$y \leftarrow y + r_0 x^2 + r_1 x z + r_2 z^2$.
The {\em valuation} of $(P,Q) \in X_2(\OK)$ is 
$$ v(P,Q) = \max \{ \min ( 2 v(P'),v(Q')) : (P',Q') 
\text { is $y$-equivalent to } (P,Q) \}. $$
It is easy to check that $v(P,Q)$ only depends 
on the $\OK$-equivalence class of $(P,Q)$. 
If $v(P)=0$, or $v(P) \ge 1$ and $Q(x,z)$ is not a square mod $\pi$, then
$v(P,Q) =0$. Otherwise we can make a $y$-substitution so that 
$v(Q) \ge 1$. Then either $v(Q)=1$ in which case $v(P,Q)=1$, or
$v(Q) \ge 2$ in which case $(P,Q)$ is non-minimal, and indeed
dividing $P$ and $Q$ through by $\pi$ and $\pi^2$ gives a $K$-equivalent 
integral model of smaller level. 
Theorem~\ref{minbq} has the following analogue.

\begin{Theorem}
\label{mingbq2}
Let $(P,Q) \in X_2(\OK)$
be a non-singular generalised binary quartic. 
Suppose that $v(P,Q)=0$ or $1$, but $(P,Q)$ has positive level.
  { 
   \begin{enumerate}
     \item 
The reduction mod $\pi$ of
$$ Q_1(x,z) = \left\{ \begin{array}{ll}
P(x,z) & \text{ if } v(P)=0, \\
\partial^2 Q / \partial x \partial z 
& \text{ if } v(P) \ge 1 \text{ and } v(P,Q)=0, \\
\pi^{-1} Q(x,z) & \text{ if } v(P) \ge 1 \text{ and } v(Q)=1 
\end{array} \right. $$ 
has a unique repeated root defined over $k$. 
     \item
The following procedure replaces $(P,Q)$ by a $K$-equivalent 
integral model of the same level. 
\begin{itemize}
\item If $v(P,Q)=1$ then make a $y$-substitution so that $v(Q) \ge 1$. 
\item 
Move the repeated root of $Q_1(x,z)$ mod $\pi$ to $(x:z)=(0:1)$.
\item Make a $y$-substitution so that $\pi \dv e$. (This is possible
since $\pi \dv n$ and every element of $k$ is a square.)
\item Replace $P(x,z)$ by $\pi^{-1} P(\pi x,z)$ and $Q(x,z)$ 
by $\pi^{-2} Q(\pi x,z)$.
\end{itemize}
\item If $(P,Q)$ is non-minimal then the procedure in (ii) 
gives $v(P,Q) \ge 2$ after at most $2$ iterations.
  \end{enumerate}
  }
\end{Theorem}
\begin{Proof}
We first show that if (i) holds for $(P,Q)$
then it holds for any $\OK$-equivalent model $(P',Q')$. 
We say that forms $f,g \in k[x,z]$ are $k$-equivalent if $f(x,z) = \lambda
g(\alpha x + \beta z,\gamma x + \delta z)$ for some $\lambda,
\alpha,\beta,\gamma,\delta \in k$ with 
$\lambda (\alpha \delta - \beta \gamma) \not=0$. Each of the following 
claims is an easy consequence of the definition of $\OK$-equivalence 
(as given in Section~\ref{g1mod}) and our assumption that $\Char(k)=2$.
\begin{itemize}
\item The reductions mod $\pi$ of $P(x,z)$ and $P'(x,z)$  are
  $k$-equivalent; in particular, $v(P)=0\iff v(P')=0$.
\item If $v(P) \ge 1$ then the reductions mod $\pi$ of 
$\partial^2 Q / \partial x \partial z $ and
$\partial^2 Q' / \partial x \partial z $ 
are $k$-equivalent; note that~$v(P,Q)=v(P',Q')$. 
\item If $v(P) \ge 1$ 
  and $v(Q) =v(Q') = 1$ 
  then the reductions mod $\pi$ of $\pi^{-1} Q(x,z)$ and
  $\pi^{-1}Q'(x,z)$  are $k$-equivalent.  
\end{itemize}
It is now clear that if (i) holds for $(P,Q)$ then it holds for $(P',Q')$. 

Next we show that the procedure in (ii) gives a well defined map on
$\OK$-equivalence classes. This does not automatically follow from
Lemma~\ref{welldef}, since we also have to consider
$y$-substitutions. 
Suppose we start with some model satisfying (i), 
and carry out the first three steps of the procedure in (ii)
in two different ways. The result is a pair of $\OK$-equivalent
models $(P,Q)$ and $(P',Q')$ related by some $[1,r,M] \in \G_2(\OK)$.
Since the reduction of $M$ mod $\pi$ fixes the repeated root $(0:1)$
we have $\pi \dv m_{21}$. Labelling the coefficients
of $(P,Q)$ in the usual way, and likewise for $(P',Q')$, we have
$\pi \dv n,e$ and $\pi \dv n',e'$. Therefore $\pi \dv r_2$.  
It is now routine to check that if (ii) holds for $(P,Q)$, i.e. $\pi \dv n,d$
and $\pi^2 \dv e$, then (ii) holds for $(P',Q')$, i.e. $\pi \dv n',d'$
and $\pi^2 \dv e'$. Moreover the transformed models are related
by $[1,(\pi r_0, r_1, \pi^{-1}r_2), \Diag(\pi,1) M \Diag(\pi^{-1},1)] 
\in \G_2(\OK)$. Thus the procedure gives a well-defined map 
on $\OK$-equivalence classes.

We are now free in the proof to replace $(P,Q)$ by any
$\OK$-equivalent model. So if $(P,Q)$ is non-minimal we may assume that
$P(\pi^sx,z) \equiv 0 \pmod{\pi^{s+1}}$ and 
$Q(\pi^sx,z) \equiv 0 \pmod{\pi^{2s+2}}$ 
for some integer $s \ge 0$. We call the least such integer $s$
the {\em slope}. As happened for binary quartics, the only possible 
slopes are $s=0,1,2$. If $s=0$ then $v(P,Q) \ge 2$ contrary to
hypothesis. 
If $s=1$ then the coefficients of $(P,Q)$ have valuations satisfying
\begin{equation*}
 \ge 0 \quad \ge 1 \quad \ge 2 \quad \quad \ge 0 \quad \ge 1 \quad 
\ge 2 \quad \ge 3 \quad \ge 4. 
\end{equation*}
If $v(P)=0$ then $P(x,z)$ mod $\pi$ has a double root at
$(x:z)=(0:1)$. Otherwise, since every element of $k$ is a square,
we can make a $y$-substitution 
$y \leftarrow y + r_0 x^2$ so that $v(Q) \ge 1$.
Then $\pi^{-1}Q(x,z)$ mod $\pi$ has 
either a triple or quadruple root at $(x:z)=(0:1)$. 
If $s=2$ then the coefficients of $(P,Q)$ have valuations satisfying
\begin{equation*}
 \ge 0 \quad \ge 1 \quad \ge 3 \quad \quad \ge 0 \quad = 0 \quad 
\ge 2 \quad \ge 4 \quad \ge 6. 
\end{equation*}
So in this case $v(P,Q)=0$. If $v(P)=0$ then $P(x,z)$ mod $\pi$ has
a double root at $(x:z)=(0:1)$. Otherwise 
$b x^2 + d z^2$ mod $\pi$ has a double 
root at $(x:z)=(0:1)$. 
In each of the cases $s=1,2$ it is now clear that 
the procedure in (ii) returns a $K$-equivalent integral model of
smaller slope. Hence at most 2 iterations are required to 
give $v(P,Q) \ge 2$, establishing (iii).

It remains to give prove (i) and (ii) in the case
$(P,Q)$ has positive level (but could be minimal).
If $(P,Q)$ is $K^\sh$-soluble then after an unramified extension
$\CC_{(P,Q)}(K) \not= \emptyset$. Then Theorem~\ref{minthm} shows
that $(P,Q)$ is non-minimal, and our earlier argument applies. Otherwise,
we show in Proposition~\ref{convprop2} below, that $(P,Q)$ is 
$\OK$-equivalent to a model whose coefficients have valuations satisfying
$$ \ge 1 \quad \ge 1 \quad \ge 2 \quad \quad =1 \quad \ge 2 \quad 
\ge 2 \quad \ge 3 \quad =3. $$
Statements (i) and (ii) are then clear.
\end{Proof}

Next we modify the algorithm for minimising quadric intersections,
as presented in Section~\ref{M4}.
First we replace $\dd$ by the map
\begin{equation}
\begin{aligned}
\dd': X_4(K) & \to X_2(K) \\
   (Q_1,Q_2) & \mapsto (P,Q) = (\pf (x Q_1 + z Q_2), \rdet(x Q_1 + z Q_2))
\end{aligned}
\end{equation}
where $\pf$ and $\rdet$ were defined in the proof of Lemma~\ref{lem:abinv}.
Then $\Delta(Q_1,Q_2) = \Delta(P,Q)$. We call $(P,Q)$ the {\em doubling} 
of $(Q_1,Q_2)$. (The reason for this name is that $\dd'$ 
acts as multiplication-by-2 on the Weil-Chatelet group.)
The analogue of Lemma~\ref{4to2cover}
(using $\dd'$ instead of $\dd$) is immediate if
$\Char(K) \not= 2$. Indeed 
the covering map $\CC_{(Q_1,Q_2)} \to \CC_{(P,Q)}$
is given by $(x_1:x_2:x_3:x_4) \mapsto (T_1:-T_2:J')$
where $J' = \frac{1}{2} ( J - l T_1^2 + m T_1 T_2 - n T_2^2 )$,
and $l,m,n$ are the coefficients of~$P$.
If $\Char(K)=2$ then the  role of $J'$ is taken by
\begin{align*} 
  J'' &= \tfrac{1}{2} \left( J - l T_1^2 + m T_1 T_2 - n T_2^2 
         + m n (l T_1 + m T_2) Q_1  + l m (n T_2 + m T_1) Q_2 \right. \\
&~ \left. \hspace{10em}  
   l^2 n^3 Q_1^2 +  l m n (l n + m^2) Q_1 Q_2 + l^3 n^2 Q_2^2 
\right).
\end{align*}
It may be verified by direct calculation that $T_1$, $T_2$ and $J''$ have 
coefficients in $\Z[X_4]$. Moreover $T_1$ and~$T_2$ 
cannot both vanish identically on $\CC_{(Q_1,Q_2)}$. 
(We checked this for the models specified in Lemma~\ref{lem1}, and then 
used the covariance of $T_1$ and~$T_2$.) Hence in all characteristics 
there is a morphism $\CC_{(Q_1,Q_2)} \to \CC_{(P,Q)}$ given by
$(x_1:x_2:x_3:x_4) \mapsto (T_1:-T_2:J'')$

The diagram in Lemma~\ref{commdiag}
(using $\dd'$ instead of $\dd$)
no longer commutes, but it does commute up to $y$-equivalence, and this
is sufficient for our purposes.

\begin{Definition}
\label{Qdefs}
  Let $Q \in k[x_1, \ldots, x_m]$ be a quadratic form in $m$ variables.  
  \begin{enumerate}    
\item The {\em kernel} $\ker(Q)$ of $Q$ is the subspace of $k^m$ defined by the 
vanishing of $Q$ and all its partial derivatives.
(Recall that $k$ is perfect, so the restriction of $Q$ to the subspace 
where all the partial derivatives vanish is the square of a linear form.)
The {\em rank} of $Q$ is $m - \dim \ker(Q)$.
\item The {\em discriminant} of $Q$ is
          \[ \Delta_m(Q)
                = \left\{ 
                    \begin{array}{ll}
                      \det ( \frac{\partial^2 Q}{\partial x_i \partial x_j}) 
                        & \text{ if $m$ is even} \\
                      \frac{1}{2}
                        \det(\frac{\partial^2 Q}{\partial x_i \partial x_j}) 
                        & \text{ if $m$ is odd.}
                    \end{array}
                  \right. \]
  \end{enumerate}
\end{Definition}
The discriminant $\Delta_m$ is a polynomial in the coefficients
of $Q$ with integer coefficients. Therefore Definition~\ref{Qdefs}(ii) 
is valid in all characteristics. 
Recall that we defined $\pf$ and $\rdet$ so that 
$\Delta_4(Q) = \pf(Q)^2 + 4 \rdet(Q)$.

\begin{Lemma}
\label{geomlem2}
Let $Q_1$ and $Q_2$ be quadratic forms in $m=3$ or $4$ variables over
a field $k$ with $\Char(k) = 2$. Assume that
\begin{itemize}
\item $\{Q_1 = Q_2 = 0 \} \subset \PP^{m-1}$ is not
a cone, i.e. $\ker(Q_1) \cap \ker(Q_2) = 0$, and 
\item if $m=3$ then $\Delta_3(x Q_1 + z Q_2)=0$, whereas if 
$m=4$ then $\pf(x Q_1 + z Q_2) = 0$ and $\rdet(x Q_1 + z Q_2)$ is a square.
\end{itemize}
Then the $k$-equivalence class of $(Q_1,Q_2)$ is uniquely determined, and
is as given in Lemma~\ref{geomlem}.
\end{Lemma}
\begin{Proof} This is similar to the proof of Lemma~\ref{geomlem}.
\end{Proof}

In Lemma~\ref{MRL} we made the hypothesis that $v(F) \ge 2$ 
where $F = \dd(Q_1,Q_2)$. This should now be replaced
by the hypothesis that $\dd'(Q_1,Q_2)$ is $y$-equivalent
to a model $(P,Q)$ with $v(P) \ge 1$ and $v(Q) \ge 2$. Then
\begin{equation}
\label{doubling}
\begin{aligned}
P(x,z) & = \pf(xQ_1 + zQ_2) + 2 h(x,z) \\
Q(x,z) & = \rdet(xQ_1 + zQ_2) -  \pf(xQ_1 + zQ_2) h(x,z) - h(x,z)^2 
\end{aligned} 
\end{equation}
for some $h \in K[x,z]$. Since $(Q_1,Q_2)$ is integral it follows
that $h \in \OK[x,z]$. Then $\pf(x\Qbar_1 + z\Qbar_2)=0$
and $\rdet(x\Qbar_1 + z\Qbar_2)$ is a square. Moreover
if $\rdet(xQ_1 + zQ_2)$ vanishes mod $\pi$ then it vanishes mod $\pi^2$.

The common nullity is $s = \dim (\ker \Qbar_1 \cap \ker \Qbar_2)$.
In the case $s=1$ we may assume that $Q_1$ and $Q_2$ reduce to
quadratic forms in $x_2,x_3,x_4$ only. 
Call these $Q'_1$ and $Q'_2$. The analogue of~(\ref{eq:s=1}) is
\[ \rdet(xQ_1 + z Q_2) \equiv 
  (\alpha x + \beta z) \Delta_3(x Q'_1 + z Q'_2) \pmod{\pi^2} \]
where $\alpha$ and $\beta$ are the coefficients
of $x_1^2$ in $Q_1$ and $Q_2$. In all other respects, the proof
of the Lemma~\ref{MRL} goes through as before. By repeated application 
of this lemma we obtain the following analogue of Lemma~\ref{iterateMRL}.
\begin{Lemma}
\label{iterateMRL2}
Let $(Q_1,Q_2) \in X_4(K)$ be a $K$-soluble non-singular 
quadric intersection. If $\dd'(Q_1,Q_2)$ 
is $y$-equivalent to an integral generalised binary quartic
then $(Q_1,Q_2)$ is $K$-equivalent to an integral quadric 
intersection $(Q'_1,Q'_2)$ such that $\dd'(Q'_1,Q'_2)$ 
is $y$-equivalent to $\dd'(Q_1,Q_2)$.
\end{Lemma}
The Minimisation Theorem for $n=4$ now follows from the Minimisation
Theorem for $n=2$ exactly as before.

The proof of Proposition~\ref{strong4} (with $\dd$ replaced
by $\dd'$) is modified as follows. 
We follow the proof of Lemma~\ref{iterateMRL2} but without the hypothesis
of $K$-solubility. This hypothesis is only used when $s \ge 2$.
In this case 
\[ (\Qbar_1, \Qbar_2) = (\alpha_{11} x_1^2 + \alpha_{12} x_1 x_2 + \alpha_{22} x_2^2,
\beta_{11} x_1^2 + \beta_{12} x_1 x_2 + \beta_{22} x_2^2)  \]
and applying the transformation $[\pi^{-1} I_2,\Diag(\pi,\pi,1,1)]$ 
gives $(R_1,R_2)$ with
\[ (\Rbar_1, \Rbar_2) = (\gamma_{33} x_3^2 + \gamma_{34} x_3 x_4 + \gamma_{44} x_4^2,
\delta_{33} x_3^2 + \delta_{34} x_3 x_4 + \delta_{44} x_4^2). \]
We must show that if $\Qbar_1$ and $\Qbar_2$ are linearly dependent 
and $\dd'(Q_1,Q_2)$ is $y$-equivalent to a model $(P,Q)$ with
$v(P) \ge 2$ and $v(Q) \ge 3$ then one of the pairs $\Qbar_1$, $\Qbar_2$
or $\Rbar_1$, $\Rbar_2$
simultaneously represents $0$ over $k$. Since $s \ge 2$ we already
know that  $\pf(xQ_1 + zQ_2)$ vanishes mod $\pi$ and $\rdet(xQ_1 + zQ_2)$ 
vanishes mod $\pi^2$. It follows by~(\ref{doubling}) that 
$\pf(xQ_1 + zQ_2)$ vanishes mod $\pi^2$ and $\pi^{-2} \rdet(xQ_1 + zQ_2)$
is a square mod $\pi$. Hence
\[ \alpha_{12} \gamma_{34} = \beta_{12} \delta_{34} = \alpha_{12} \delta_{34} + \beta_{12} \gamma_{34} = 0 \]
and
\begin{align*}
\alpha_{12}^2 (\gamma_{33} \delta_{44} + \gamma_{44} \delta_{33}) + \gamma_{34}^2 (\alpha_{11} \beta_{22} 
+ \alpha_{22} \beta_{11}) &= 0 \\
\beta_{12}^2 (\gamma_{33} \delta_{44} + \gamma_{44} \delta_{33}) + \delta_{34}^2 (\alpha_{11} \beta_{22} 
+ \alpha_{22} \beta_{11}) &= 0 {\rlap .}
\end{align*}
Since $\Qbar_1$ and $\Qbar_2$ are linearly dependent we have
$\alpha_{11} \beta_{22} + \alpha_{22} \beta_{11} = 0$. So either $\alpha_{12}=\beta_{12}=0$,
in which case $\Qbar_1$ and $\Qbar_2$ simultaneously represent $0$
over $k$, or $\gamma_{34} = \delta_{34} = \gamma_{33} \delta_{44} + \gamma_{44} \delta_{33} = 0$
in which case  $\Rbar_1$ and $\Rbar_2$ simultaneously represent $0$
over $k$.

\subsection{Minimisation over global fields}
\label{sec:globalmin}

We have so far presented theorems and algorithms 
for minimising genus one models defined over local fields. We now
discuss the global situation, and in particular prove
Theorem~\ref{thmglobal}.  
The following is a more precise version of that theorem.
A genus one model defined over a number 
field $K$ is called {\em integral} if its coefficients belong to the 
ring of integers~$\OK$. 

\begin{Theorem}
\label{thmglobal2}
  Let $n=2,3$ or $4$.
  Let $K$ be a number field of class number one. Let $\Phi \in X_n(K)$ 
  be a non-singular genus one model. If $\CC_\Phi$ is locally soluble at all 
  finite places of $K$ then $\Phi$ is $K$-equivalent to an integral genus 
  one model with the same discriminant as a global minimal model
  for the Jacobian $E$ of $\CC_\Phi$.
\end{Theorem}

\begin{Proof} 
To deduce this result directly from the statement of the Minimisation
Theorem (Theorem~\ref{minthm}) one is naturally led to use a version of
strong approximation. See \cite{minbqtc} for details in the cases
$n=2,3$.  The case $n=4$ is similar. Although these proofs are not
difficult, it is a notable advantage of the algorithmic approach taken
in this section that the passage from local to global becomes a
triviality. 

Indeed, suppose $K$ is a number field with class number one.  Let
$\pp = \pi \OK$ be a prime of~$K$ and put $k = \OK/\pp$. Then for any pair of
$m$-dimensional subspaces $U,V \subset k^n$ there exists $M \in
\SL_n(\OK)$ whose reduction mod $\pp$ takes $U$ to $V$.  (Indeed, the
case $\dim U = \dim V = 1$ is Lemma~\ref{trans}, and the general case
is similar.) 
We can therefore follow the algorithms for minimising at $\pp$, 
using $\pi$ as the uniformiser, 
without changing the level (or integrality) at other primes.

After first scaling the given model to be integral at all primes, we
apply this procedure to the finite number of primes at which the
resulting model has positive level.  This gives an integral model
which has level zero at all primes of~$K$.  By definition of level,
this model has the same discriminant as a global minimal model 
for~$E$, up to a unit factor.  Since this unit must be a $12$th
power, a final scaling by a suitable global unit gives the result.
\end{Proof}

Theorem~\ref{thmglobal} is an immediate corollary since, as recalled
in the introduction, every $n$-covering which is locally soluble 
at all places of~$K$, has a degree-$n$ model.

To extend this theorem to a general number field~$K$, we may replace
integrality by $S$-integrality, where $S$ is a (finite) set of primes
generating the class group, so that the ring of $S$-integers is a
principal ideal domain.  The minimal model may then only be
$S$-integral rather than integral.  Just as with Weierstrass models
for elliptic curves, there may be no global minimal model when the
class number is greater than~$1$.  In practice, we can alternatively
find models which are simultaneously minimal at all primes
in any given finite set, while being at least integral at all other primes.

Similar results may be deduced from our local results in the case
where $K$ is a function field, {\em i.e.}, a finite extension of~$\F_q(t)$.

\section{Minimisation of insoluble genus one models}
\label{sec:critmodels}

We return to working over a discrete valuation field $K$ as
specified in Section~\ref{min:statements}.
In this section we prove the Converse Theorem
(Theorem~\ref{mainthmA}(ii)). This shows that the Strong Minimisation
Theorem (Theorem~\ref{mainthmA}(i)) is best possible.

\begin{Definition} 
\label{defcrit} \strut
\begin{itemize}
  \item[(a)] A generalised binary quartic $(P,Q) \in X_2(\OK)$ 
    is {\em critical} if the valuations of its coefficients 
    $l,m,n,a,b,c,d,e$ satisfy
    $$ \ge 1 \quad \ge 1 \quad \ge 2 \quad \quad =1 \quad \ge 2 \quad 
    \ge 2 \quad \ge 3 \quad =3. $$
  \item[(b)] A ternary cubic $F \in X_3(\OK)$ is {\em critical} 
    if the valuations of its coefficients satisfy the inequalities 
    indicated in the following diagram.
      \begin{equation*}
        \begin{array}{cccccccc}
          & & & \multicolumn{2}{c}{ \, z^3 \,\, } \\
          & & \multicolumn{2}{c}{\, x z^2 \!\!\!} & \multicolumn{2}{c}{ y z^2 } \\
          & \multicolumn{2}{c}{ x^2 z } & \multicolumn{2}{c}{ xyz } & 
            \multicolumn{2}{c}{ y^2 z } \\
          \multicolumn{2}{c}{ x^3 } & \multicolumn{2}{c}{ x^2 y } & 
            \multicolumn{2}{c}{ x y^2 } & \multicolumn{2}{c}{ \,\, y^3 }
        \end{array}
        \qquad
        \begin{array}{cccccccc}
          & & & \multicolumn{2}{c}{ = 2 } \\
          & & \multicolumn{2}{c}{ \ge 2 } & \multicolumn{2}{c}{ \ge 2 } \\
          & \multicolumn{2}{c}{ \ge 1 } & \multicolumn{2}{c}{ \ge 1 } & 
            \multicolumn{2}{c}{ \ge 2 } \\
          \multicolumn{2}{c}{ =0 } & \multicolumn{2}{c}{ \ge 1 } & 
            \multicolumn{2}{c}{ \ge 1 } & \multicolumn{2}{c}{ = 1 }
        \end{array} 
      \end{equation*}
  \item[(c)] A quadric intersection $(Q_1,Q_2) \in X_4(\OK)$ is {\em critical} if 
    the reductions of $Q_1$ and $Q_2$ mod $\pi$ are quadratic
    forms in $x_1$ and $x_2$ with no common root in $\PP^1(\kbar)$, and
    on putting
\begin{equation*}
(R_1,R_2) = [ \pi^{-1} I_2, \Diag(\pi,\pi,1,1) ] (Q_1,Q_2)
\end{equation*}
the reductions of $R_1$ and $R_2$ mod $\pi$ are quadratic
forms in $x_3$ and $x_4$ with no common root in $\PP^1(\kbar)$.
\end{itemize}
\end{Definition}

We show in the next three lemmas that critical models are 
insoluble, minimal and of positive level. We then show
(for $n=2,3$) that every $K^\sh$-insoluble model is $K$-equivalent 
to a critical model. There is a corresponding result for models of
degree $n=4$.

\begin{Lemma} 
\label{critprops}
Critical models are insoluble over $K$.
\end{Lemma}
\begin{Proof}
We give details in the case $n=2$. 
Suppose $(x,y,z) \in K^3$ is a non-zero solution
of $y^2 + P(x,z)y = Q(x,z)$. Clearing denominators we may assume that
$\min \{ v(x),v(z) \} = 0$. It follows that $y \in \OK$.
Then reducing the equation mod $\pi^i$ for $i=1,2,3,4$ we successively 
deduce $\pi \dv y$, $\pi \dv x$, $\pi^2 \dv y$ and $\pi \dv z$.
In particular $\min \{ v(x),v(z) \} > 0$. 
This is the required contradiction. The cases $n=3,4$ are similar.
\end{Proof}
Since the definition of a critical model is unchanged by an
unramified field extension, it follows immediately that
critical models are insoluble over $K^\sh$. 

\begin{Lemma} 
\label{critprops2}
Critical models are minimal.
\end{Lemma}
\begin{Proof}
In the cases $n=2,3$ we give a very quick proof.
Indeed, if $\Phi$ were non-minimal, then our algorithms in 
Sections~\ref{sec:minbq},~\ref{M3} and~\ref{min_char2}
would succeed in reducing the level. 
But on the contrary, when given a critical model, these algorithms 
endlessly cycle between two or three $\OK$-equivalence classes. 
(Treating the case $n=4$ in the same way would give a circular
argument, as the current lemma was cited at the end of~Section~\ref{M4}.)

Alternatively we can imitate the proof of Lemma~\ref{critprops}. 
We give details in the case $n=4$.
We define
\[ s(Q_1,Q_2) = \max \{ -v(\det M) : [M,I_4](Q_1,Q_2) \in X_4(\OK) \}. \]
Suppose $[M,N] \in \G_4(K)$ is a transformation taking the
critical model $\Phi = (Q_1,Q_2)$ to an integral model of smaller
level. We may assume that $N$ has entries in $\OK$, not all
in $\pi \OK$. Let $\xi_j(x_1, \ldots, x_4) = \sum_{i=1}^4 n_{ij} x_i$. 
For $i=1,2$ we put 
\[ Q_i \circ N = Q_i(\xi_1, \ldots, \xi_4) \in \OK[x_1 ,\ldots, x_4]. \]
Our hypothesis is that $s(Q_1 \circ N, Q_2 \circ N) > v(\det N)$. 

If $v(Q_1 \circ N) = 0$ then replacing $Q_2$ by $Q_2 + \lambda Q_1$ 
for suitable $\lambda \in \OK$ we may assume that $v(Q_2 \circ N) > v(\det N)$.
To understand this last condition, we put $N$ in Smith normal form.
Explicitly we write $ N = U \Diag(\pi^a, \pi^b, \pi^c,1) V $
for some $U,V \in \GL_4(\OK)$ and $a \ge b \ge c \ge 0$.
Since $v(Q_2)=0$ we must have $2a > v(\det N) = a+ b + c$ and 
therefore $a-b+c \ge 1$.
It follows that $Q_2 \circ U \equiv x_1 ( \sum_{i=1}^4 \epsilon_i x_i) 
\pmod{\pi^2}$ for some $\epsilon_i \in \OK$ with $\epsilon_2 \equiv
\epsilon_3 \equiv \epsilon_4 \equiv 0 \pmod{\pi}$. In other words,
$Q_2 \equiv \mu \ell_1 \ell_2 \pmod{\pi^2}$ for some $\mu \in \OK$ 
and linear forms $\ell_1, \ell_2 \in \OK[x_1, \ldots, x_4]$ with
$\ell_1 \equiv \ell_2 \pmod{\pi}$. This contradicts the definition
of a critical model (as it would follow that $R_2$ vanishes mod $\pi$). 
Hence $v(Q_1 \circ N) \ge 1$. Similarly 
$v(Q_2 \circ N) \ge 1$. Since $\Qbar_1$ and $\Qbar_2$ are binary
quadratic forms with no common root 
we deduce $\xi_1 \equiv \xi_2 \equiv 0 \pmod{\pi}$. 
Let $\xi'_i = \pi^{-1} \xi_i$ for $i=1,2$. We put
\begin{equation*}
(R_1,R_2) = [ \pi^{-1} I_2, \Diag(\pi,\pi,1,1) ] (Q_1,Q_2).
\end{equation*}
Let $N'$ be the matrix with columns the coefficients of 
$\xi_3, \xi_4, \xi'_1, \xi'_2$. Then $(R_1,R_2)$ is a critical model
and  $s(R_1 \circ N', R_2 \circ N') > v(\det N')$. Repeating
the same arguments we deduce $\xi_3 \equiv \xi_4 \equiv 0 \pmod{\pi}$. 
This contradicts our scaling of the matrix $N$.
\end{Proof}

The next lemma describes the possible levels of a critical model.
For this we need to work explicitly with the ``$a$-invariants''
defined in the proof of Lemma~\ref{lem:abinv}. 
Although $a_1, \ldots, a_6$ are not invariants 
(in the sense of Definition~\ref{def:inv}), they are isobaric in 
the sense that 
\begin{align*}
n=2: \quad & & a_i \circ [\mu,0,\Diag(\xi_1,\xi_2)] 
            & = (\mu \xi_1 \xi_2 )^i a_i \\
n=3: \quad & & a_i \circ [\mu,\Diag(\xi_1,\xi_2,\xi_3)]
            & = (\mu \xi_1 \xi_2 \xi_3)^i a_i \\
n=4: \quad & & a_i \circ [\Diag(\mu_1,\mu_2),\Diag(\xi_1,\xi_2,\xi_3,\xi_4)]
            &  = (\mu_1 \mu_2 \xi_1 \xi_2 \xi_3 \xi_4)^i a_i
\end{align*}
for all $i$. (We use the notation for transformations of genus
one models introduced in Section~\ref{g1mod}.)
In the following we write $t^{(n)}$ as a short-hand for $\pi^{-n} t$.

\begin{Lemma}
\label{lem:critlev}
The level of a critical model is at least $1$ and equal to $1$
if $\Char(k) \ndv n$.
\end{Lemma}

\begin{Proof}
Case $n=2$. By~(\ref{jac2gen}) we have $\pi^i \dv a_i$ for all $i$. 
A convenient way to check this is to note that
 $\pi^{-3/2} P ( \pi^{1/2} x , z)$ and $\pi^{-3} Q( \pi^{1/2} x, z)$ 
have coefficients in $\OK[\pi^{1/2}]$, and then to use the
isobaric property.
It follows that $(P,Q)$ has positive level.
Now suppose that $\Char(k) \not=2$ and $(P,Q)$ has level greater
than $1$. Completing the square we may assume that $l=m=n=0$. Then
$a_1=a_3=0$ and $y^2 = x^3 + a_2^{(2)} x^2 + a_4^{(4)} x + a_6^{(6)}$
is an integral Weierstrass equation of positive level. 
According to Tate's algorithm the cubic polynomial
$$ x^3 + a_2^{(2)} x^2 + a_4^{(4)} x + a_6^{(6)} \equiv
(x + c^{(2)})(x^2- 4 a^{(1)} e^{(3)}) \pmod{\pi} $$
has a triple root defined over $k$. This contradicts the definition of a 
critical model.

\paragraph{Case $n=3$.} 
By~(\ref{jac3gen}) we have $\pi^i \dv a_i$ for all $i$.
  A convenient way to check this is to note that 
  $\pi^{-2} F( \pi^{2/3} x, \pi^{1/3} y, z)$ has coefficients in 
  $\OK[\pi^{1/3}]$,
  and then to use the isobaric property.
  It follows that $F$ has positive level. 
Now suppose that $\Char(k) \not =3$ and $F$ has level greater than 1.
 Then
  $$ y^2 + a_1^{(1)} xy + a_3^{(3)} y
      = x^3 + a_2^{(2)} x^2 + a_4^{(4)} x + a_6^{(6)} $$
  is an integral Weierstrass equation of positive level. 
  By~(\ref{jac3gen}) we find
  $a_2^{(2)} \equiv a_4^{(4)} \equiv 0 \pmod{\pi}$ and 
  \begin{align*}
       a_1^{(1)} & \equiv  m^{(1)} \pmod{\pi} \\
       a_3^{(3)} & \equiv 9 a b^{(1)} c^{(2)} \pmod{\pi} \\
       a_6^{(6)} & \equiv -27(ab^{(1)}c^{(2)})^2
                             + ab^{(1)}c^{(2)}(m^{(1)})^3 \pmod{\pi}. 
  \end{align*}
  So it suffices to show that if there is a Weierstrass equation over $k$
  of the form 
  $$ y^2 + \alpha x y + 9 \beta y = x^3 + (\alpha^3 - 27 \beta) \beta $$
  with $c_4= \Delta=0$, then $\beta=0$. We compute
  $c_4 = \alpha(\alpha^3 - 216 \beta)$ and 
  $\Delta = -\beta(\alpha^3 + 27 \beta)^3$. Since $216+27= 3^5$ is non-zero
  in $k$, it follows that $\beta = 0$ as required.

\paragraph{Case $n=4$.} 
The quadric intersection
$ [ \pi^{-1} I_2, \Diag(\pi^{1/2},\pi^{1/2},1,1)](Q_1,Q_2) $
has coefficients in $\OK[\pi^{1/2}]$. It follows by the isobaric
property of the $a$-invariants that $\pi^i \dv a_i$ for all $i$ 
and hence that $(Q_1,Q_2)$ has positive level.
Now suppose that $\Char(k) \not= 2$.
Then $F = {\dd}(Q_1,Q_2)$ 
satisfies $F(x,z) \equiv \pi^2 f_1(x,z) f_2(x,z) \pmod{\pi^3}$
where $f_1,f_2 \in \OK[x,z]$ are binary quadratic forms, neither
having a repeated root mod $\pi$. (So their product cannot 
have a triple or quadruple root.) 
It follows by Theorem~\ref{minbq}(i) that $F$ and hence
$(Q_1,Q_2)$ has level 1. 
\end{Proof}

\begin{Example}
The following examples of critical models, all of level 2, show
that the hypothesis $\Char(k) \ndv n$ cannot be removed from 
Lemma~\ref{lem:critlev}. 
$$ \begin{array}{ccl}
K = \Q_2 & &  y^2 = 2 x^4 + 24 x^2 z^2 + 8 z^4 \\
K = \Q_3 & &    x^3 + 3 y^3 + 9 z^3 + 18 x y z = 0 \\
K = \Q_2 & & x_1^2 + 2 x_3^2 + 4 x_2 x_4 = x_2^2 + 2 x_4^2 + 4 x_1 x_3 = 0 
\end{array} $$
\end{Example}

The following proposition completes the proof of 
Theorem~\ref{mainthmA}(ii). The doubling map ${\dd}'$ was defined in 
Section~\ref{min_char2}. (If $\Char(k) \not=2$ then we can work with
${\dd}$ instead.)

\begin{Proposition}
\label{convprop2}
Let $\Phi \in X_n(\OK)$ be a $K^\sh$-insoluble minimal genus one model.
\begin{enumerate}
\item If $n=2$ or $3$ then $\Phi$ is $\OK$-equivalent to a critical model.
\item If $n=4$ then $\Phi$ is $K$-equivalent to either a critical model
or an integral model $(Q_1,Q_2)$ with ${\dd}'(Q_1,Q_2)$ critical.
\end{enumerate}
\end{Proposition}

First we need three lemmas. 
\begin{Lemma}
\label{findsmooth}
Let $k$ be an algebraically closed field. Suppose that either
{ \renewcommand{\theenumi}{\alph{enumi}}
\begin{enumerate}
\item $\Phi = (P,Q) \in X_2(k)$ and $P^2 + 4 Q$ is not identically zero,
\item $\Phi = (F) \in X_3(k)$ is non-zero and is not 
the cube of a linear form,
\item $\Phi = (Q_1,Q_2) \in X_4(k)$ and every quadric in the pencil
spanned by $Q_1$ and $Q_2$ has rank at least $2$.
\end{enumerate} }
Then $\CC_{\Phi}$ has a smooth $k$-point (on some $1$-dimensional component).
\end{Lemma}
\begin{Proof}
For $n=2,3$ this is clear. In the case $n=4$ we are looking for
a transverse point of intersection of $Q_1$ and $Q_2$, {\em i.e.}
a point where the Jacobian matrix has rank $2$. We prove the result
more generally for intersections of two quadrics in $m$ variables.
This enables us to reduce to the case $\ker(Q_1) \cap \ker(Q_2) = 0$.
Now let $P$ be a singular point on the quadric intersection. 
(If there is no such point there is nothing to prove.) Then moving this
point to $(1:0: \ldots :0)$ we may assume that
$Q_1 = x_1 x_2 + g_1(x_2, \ldots, x_m)$ and 
$Q_2 = g_2(x_2, \ldots, x_m)$ for some $g_1$ and $g_2$.
Since $\rank(Q_2) \ge 2$ we can pick a smooth point $(x_2: \ldots: x_m)$ 
on $\{Q_2 = 0\} \subset \PP^{m-2}$ with $x_2 \not= 0$. Then solving
the equation $Q_1= 0$ for $x_1$ gives the required transverse point 
of intersection on $\{Q_1= Q_2 = 0\}$.
\end{Proof}

\begin{Lemma}
\label{hlift}
Let $\Phi \in X_n(\OK)$ be a $K^\sh$-insoluble minimal genus one model.
{ \renewcommand{\theenumi}{\alph{enumi}}
\begin{enumerate}
\item If $n=2$ then $\Phi=(P,Q)$ with $v(P,Q)=1$. 
Moreover if $v(Q)=1$ then the reduction of $\pi^{-1}Q(x,z)$ mod $\pi$ 
has either two double roots or a quadruple root (over $\kbar$).
\item If $n=3$ then $\Phi$ is a ternary cubic whose reduction mod $\pi$
is (a constant times) the cube of a linear form.
\item If $n=4$ then there is a rank 1 quadric in the reduced pencil,
{\em i.e.} if
$\Phi=(Q_1,Q_2)$ then $\rank( \lambda \Qbar_1 + \mu \Qbar_2) =1$ for some 
$(\lambda:\mu) \in \PP^1(\kbar)$.
\end{enumerate} }
\end{Lemma}

\begin{Proof}
We recall that $K^\sh$ has residue field $\kbar$. 
The idea of the proof is that if $\Phi$ 
is not of the form listed, then we can use Lemma~\ref{findsmooth} 
to find a smooth $\kbar$-point on the reduction, and 
use the Henselian property to lift it to a $K^\sh$-point,
thereby obtaining a contradiction.

A little more needs to be said in the case $n=2$. 
If $\Char(k) \not = 2$ then completing the 
square gives $v(P) \ge 1$ and Lemma~\ref{findsmooth} shows 
that $v(Q) \ge 1$. If $\Char(k) = 2$ then Lemma~\ref{findsmooth} 
shows that $v(P) \ge 1$. If $Q(x,z)$ mod $\pi$ had a simple root 
over $\kbar$ then we could lift to a $K^\sh$-point on $\CC_{(P,Q)}$ 
with $y=0$. It follows that $Q(x,z)$ is a square mod $\pi$. 
So by a $y$-substitution we may suppose $v(Q) \ge 1$. 
In all residue characteristics we now have $v(P) \ge 1$ and 
$v(Q) \ge 1$. We cannot have $v(Q) \ge 2$ since $(P,Q)$ is minimal.
If $\pi^{-1} Q(x,z)$ mod $\pi$ had a simple root over $\kbar$
then we could lift to a $K^\sh$-point on $\CC_{(P,Q)}$ with $y=0$.
It follows that this polynomial has either two double roots 
or a quadruple root.
\end{Proof}

\begin{Lemma}
\label{fliplimit}
Suppose $(P,Q), (P',Q') \in X_2(\OK)$ are $K$-equivalent 
models of the same level
related by a substitution $[\mu,r,M] \in \G_2(K)$ where $M \in \GL_2(K)$
has Smith normal form $\Diag(1,\pi^s)$. Then $v(\Delta(P,Q)) \ge 2s$.
\end{Lemma}
\begin{Proof}
Let $(P,Q)$ have coefficients $l,m,n,a,b,c,d,e$.
Replacing our models by $\OK$-equivalent ones we may assume
$\mu = \pi^{-s}$ and $M = \Diag(\pi^s,1)$. If we assume for simplicity
that $r=0$, then we have $\pi^s \dv n,d$ and $\pi^{2s} \dv e$. 
Since the discriminant $\Delta \in \Z[X_2]$ belongs to the 
ideal $(n^2,nd,d^2,e)$ it follows that $v(\Delta(P,Q)) \ge 2s$.

For general $r$ we can write the
transformation $[\pi^{-s},r,\Diag(\pi^s,1)]$ either as 
\[ y \leftarrow \pi^s y + r_0 x^2 + r_1 x z + r_2 z^2  
\quad  \text{ followed by }  \quad x \leftarrow \pi^s x \]
or as
\[ x \leftarrow \pi^s x \quad \text{ followed by }  \quad
y \leftarrow \pi^s (y + \pi^s r_0 x^2 + r_1 x z + \pi^{-s} r_2 z^2).  \]
Since $Q'$ has coefficients in $\OK$ we have
$v(r_0^2 + r_0 l - a) \ge -2s$ and
$v(r_2^2 + r_2 n - e) \ge 2s$. Hence $\pi^s r_0, r_2 \in \OK$. 
So replacing our models by $\OK$-equivalent ones we may
assume that $r_0=r_2=0$. Then the middle coefficient of $Q'$ gives 
$v(r_1^2 + r_1 m - c) \ge 0$ and hence
$r_1 \in \OK$. Once more replacing $(P,Q)$ by an $\OK$-equivalent model
we may assume that $r_0=r_1=r_2 =0$. Our earlier proof now applies.
\end{Proof}

\begin{ProofOf}{Proposition~\ref{convprop2}}
We split into the cases $n=2,3,4$.

\medskip

\paragraph{Case $n=2$} 
Applying Lemma~\ref{hlift} to $\Phi = (P,Q)$ we may assume that
$v(P) \ge 1$, $v(Q) =1$, and $\pi^{-1}Q(x,z)$
mod $\pi$ has either two double roots or a quadruple root.

We first rule out the possibility of two double roots. After
an unramified field extension we may assume that these roots are
defined over $k$. So without loss of generality 
$Q(x,z) \equiv \pi x^2 z^2 \pmod{\pi^2}$. 
We replace $P(x,z)$ by $\pi^{-1} P(\pi x,z)$ 
and $Q(x,z)$ by $\pi^{-2} Q(\pi x,z)$. By Lemma~\ref{hlift}
we again have $v(P,Q) \ge 1$. We make a substitution 
$y \leftarrow y + r_2 z^2$ so that $v(P) \ge 1$ and $v(Q) \ge 1$. 
Now $\pi^{-1}Q(x,z)$ mod $\pi$ has a double root at $(x:z)=(1:0)$.
By Lemma~\ref{hlift} it has a second double root, say at $(\lambda:1)$. 
We make the substitution $x \leftarrow x + \lambda z$.
Then $Q(x,z) \equiv \pi x^2 z^2 \pmod{\pi^2}$. We can now 
repeat this process indefinitely. It follows by Lemma~\ref{fliplimit}
that $\Delta(P,Q)=0$. This is the required contradiction.

It remains to consider the case of
a quadruple root, say $Q(x,z) \equiv
\pi x^4 \pmod{\pi^2}$. Let $l_1,m_1,n_1,a_1,b_1,c_1,d_1,e_1$
be the coefficients of $P_1(x,z) = \pi^{-1} P(\pi x,z)$ 
and $Q_1(x,z) = \pi^{-2} Q(\pi x,z)$. 
By Lemma~\ref{hlift} we can make a substitution $y \leftarrow y + r_2 z^2$
so that $\pi \dv n_1,e_1$. Then $\pi^{-1}Q_1(x,z)$ mod $\pi$ has at least
a triple root at $(x:z)=(1:0)$. So by Lemma~\ref{hlift} we have
$\pi^2 \dv d_1$ and $v(e_1)=1$. The coefficients of $(P,Q)$ 
now satisfy the definition of a critical model.

\medskip

\paragraph{Case $n=3$.} 
  By Lemma~\ref{hlift} our ternary cubic $F$ 
  must reduce mod~$\pi$ to the cube of a linear form.
  So without loss of generality, we have
  \[ F = \pi \,f_3(y,z) + \pi \,f_2(y,z)\,x + \pi \,f_1(y,z)\,x^2 
          + a\,x^3 \,.
  \]
  with $\pi \ndv a$. Then $F_1(x, y, z) = \pi^{-1} F(\pi x, y, z)$
  is a minimal  
  ternary cubic and by Lemma~\ref{hlift} its reduction
  mod $\pi$ is the cube of a linear form in $y$ and~$z$. After
  a suitable transformation of $y$ and~$z$, we may assume that
  $f_3(y,z) \equiv b y^3 \pmod{\pi}$ with $\pi \ndv b$ (otherwise $F$ 
  would not be minimal).
  Now $F_2(x, y, z) = \pi^{-1} F_1(x, \pi y, z)$ is again a minimal
  ternary cubic, and its reduction mod $\pi$ is $(c' x + c z)z^2$.
  Again this must be a non-zero cube. So $c'=0$ and $c$ is a unit.
  The coefficients of $F$ now satisfy the definition of a critical model.

\medskip

\paragraph{Case $n=4$.}
We divide the proof into the following two lemmas.

\begin{Lemma} 
\label{newlem1}
Let $(Q_1,Q_2) \in X_4(\OK)$ be a $K^\sh$-insoluble
minimal quadric intersection. Let 
$s = \dim (\ker(\Qbar_1) \cap \ker(\Qbar_2))$ be the common nullity of
the reduced pencil. 
\begin{enumerate}
\item If $s \le 1$ then the reduced pencil contains a unique rank $1$
quadric, and the following procedure replaces $(Q_1,Q_2)$ by a $K$-equivalent
minimal quadric intersection with $s \ge 1$. 
\begin{itemize}
\item Make a $\GL_2(\OK) \times \GL_4(\OK)$-transformation so 
that $Q_2 \equiv x_1^2 \pmod{\pi}$. 
\item Apply the transformation $[\Diag(1,\pi^{-1}),\Diag(\pi,1,1,1)]$.
\end{itemize}
\item If $s \ge 2$ then $(Q_1,Q_2)$ is $\OK$-equivalent to a critical model.
\end{enumerate}
\end{Lemma}

\begin{Proof}
(i) By Lemma~\ref{hlift} there is a rank~1 quadric in the reduced 
pencil. It is unique (and therefore defined over $k$) 
as we would otherwise have $s \ge 2$. The remaining statements are clear. \\
(ii) We may assume that $\Qbar_1$ and $\Qbar_2$ are binary quadratic 
forms in $x_1$ and $x_2$. Since the model is minimal, these forms have 
no common root in $\PP^1(\kbar)$. We put
\[ (R_1,R_2) = [ \pi^{-1} I_2, \Diag(\pi,\pi,1,1) ] (Q_1,Q_2). \]
Then $R_1$ and $R_2$ reduce to binary quadratic forms in 
$x_3$ and $x_4$. Again, since the model is minimal, these forms 
have no common root in $\PP^1(\kbar)$. Hence $(Q_1,Q_2)$ is critical.
\end{Proof}

\begin{Lemma} 
\label{newlem2}
Let $\Phi \in X_4(\OK)$ satisfy the hypotheses 
of Lemma~\ref{newlem1} with $s=1$. If the procedure in Lemma~\ref{newlem1}(i)
may be iterated indefinitely, then $\Phi$ is $\OK$-equivalent
to a quadric intersection $(Q_1,Q_2)$ where the valuations 
of the coefficients of $Q_1$ and $Q_2$ satisfy the inequalities 
indicated in the following diagram:
\begin{equation*}
\begin{array}{cccc} 
x_1^2 & x_1 x_2 & x_1 x_3 & x_1 x_4 \\
      &  x_2^2  & x_2 x_3 & x_2 x_4 \\
      &         &  x_3^2  & x_3 x_4 \\
      &         &         &  x_4^2  
\end{array} 
\quad 
\begin{array}{cccc} 
\ge 0 & \ge 0   &  = 0    & \ge 1   \\
      &  = 0    &  \ge 1  & \ge 1   \\
      &         &  \ge 1  & \ge 1   \\
      &         &         &  = 1    
\end{array} 
\quad 
\begin{array}{cccc} 
 = 0  &  \ge 1  &  \ge 1  &  \ge 1  \\
      &  \ge 1  &  \ge 1  &  = 1    \\
      &         &   = 1   &  \ge 2  \\
      &         &         &  \ge 2 {\rlap .}  
\end{array} 
\end{equation*}
\end{Lemma}

\begin{Proof} We may assume that $\Phi = (Q_1, Q_2)$ has reduction
\begin{equation}
\label{format4}
 (\Qbar_1,\Qbar_2) = ( x_1 \ell(x_2,x_3) + f(x_2,x_3) , c x_1^2 ) 
\end{equation}
for some $c \in k$ and $\ell,f \in k[x_2,x_3]$. 
Since $(Q_1,Q_2)$ is minimal
we have $c f \not= 0$. So the reduction is (set-theoretically) either
a line or a pair of lines. We show in the case of a pair of 
lines that the procedure in Lemma~\ref{newlem1}(i) must give $s \ge 2$
after a finite number of iterations (bounded in terms of the
valuation of the discriminant). The first iteration gives
$(R_1,R_2)$ with
\[ (\overline{R}_1,\overline{R}_2) 
  = ( f(x_2,x_3), g(x_2,x_3,x_4) ) \]
for some $g \in k[x_2,x_3,x_4]$. Since $f$ has rank $2$ we may
assume on replacing $R_2$ by $R_2 + \lambda R_1$ for suitable
$\lambda \in \OK$ that $g$ has rank $1$. 
If $g$ has no coefficient of $x_4^2$ then $s \ge 2$.
Otherwise a $\GL_4(\OK)$-transformation puts
$(\overline{R}_1,\overline{R}_2)$ in the form~(\ref{format4})
with $\ell = 0$ (and the same $f$ as before). 
The process is then repeated. By considering the effect on the doubling 
it follows by Lemma~\ref{fliplimit} that only 
finitely many iterations are possible.

It remains to consider the case where the reduction is (set-theoretically)
a line. We may assume that $\Phi = (Q_1,Q_2)$ and its transforms 
\begin{align*}
(R_1,R_2) & = [\Diag(1,\pi^{-1}), \Diag(\pi,1,1,1)](Q_1,Q_2) \\
(S_1,S_2) & = [ \Diag(\pi^{-1},1) ,\Diag(1,\pi,1,1)](R_1,R_2)
\end{align*}  
under the first two iterations have reductions
\begin{align}
\label{eqn1}
(\overline{Q}_1,\overline{Q}_2) &= 
(x_1(\alpha_1 x_1 + \alpha_2 x_2 + \alpha_3 x_3) + x_2^2 , \,\, x_1^2) \\
\label{eqn2}
(\overline{R}_1,\overline{R}_2) &= 
(x_2^2, \,\, x_2(\beta_2 x_2 + \beta_3 x_3 + \beta_4 x_4) + g(x_3,x_4)) \\
\label{eqn3}
(\overline{S}_1,\overline{S}_2) &=
(\alpha_3 x_1 x_3 + \lambda x_3^3 + \mu x_3 x_4 + \nu x_4^2, \,\, g(x_3,x_4))
\end{align}
for some $\alpha_i,\beta_i,\lambda,\mu,\nu \in k$ and $g \in k[x_3,x_4]$.
By~(\ref{eqn1}) we have $\alpha_3 \not= 0$ (otherwise $s \ge 2$).
Since the reduction cannot be a pair of lines, we see first
by~(\ref{eqn2}) that $g$ has rank $1$, and then by~(\ref{eqn3}) that
$g = \gamma x_3^2$ for some $\gamma \not= 0$.
Finally~(\ref{eqn2}) and~(\ref{eqn3}) show that $\beta_4 \not= 0$ and 
$\nu \not= 0$ (otherwise $s \ge 2$).
The valuations of the coefficients of $Q_1$ and $Q_2$ now satisfy the 
inequalities indicated in the statement of the lemma.
\end{Proof}

Proposition~\ref{convprop2}(ii) follows from the last two lemmas
and the observation that if $(Q_1,Q_2)$ satisfies the conclusions
of Lemma~\ref{newlem2} then its doubling is critical.
\end{ProofOf}


\section{Reduction} 
\label{sec:R}

In this section, we assume that the ground field is~$\Q$. The main
reason for this is that a comparable theory of reduction over a
general number field has not yet been sufficiently developed.

Let $\CC \subset \PP^{n-1}$ be a genus one normal curve defined
over~$\Q$ of degree~$n$ (or, if $n = 2$, let $\CC\to\PP^1$ be a double
cover) with points everywhere locally, so that $\CC$ represents an
element of the $n$-Selmer group of its Jacobian elliptic curve~$E$. If
$n \in \{2,3,4\}$, we can, by the results and algorithms of the
previous sections, assume that $\CC=\CC_{\Phi}$ where $\Phi$ is a
genus one model which is both integral and minimal, so that its
invariants~$c_4$, $c_6$ and~$\Delta$ coincide with those of a minimal
model of~$E$.  This means that the invariants are as small as possible
(in absolute value). However, it does not necessarily mean that the
equations defining~$\CC$ will have small coefficients. To achieve
this, we will employ {\em reduction}.  Leaving aside the aesthetic
value of equations with small coefficients, the main benefit of a
reduced model is that further computations like searching for rational
points on~$\CC$ or performing further descents on~$\CC$ are greatly
facilitated.

The idea of reduction is to find a unimodular transformation ({\em
  i.e.}, an invertible integral linear change of coordinates
on~$\PP^{n-1}$) that makes the equations defining~$\CC$
smaller. Unimodular transformations have the property of preserving
the integrality and invariants of the model, so they will not destroy
its minimality.  In the language of Section~\ref{g1mod}, a unimodular
transformation is just a~$\Z$-equivalence.

If we were allowed to make a coordinate change from $\SL_n(\C)$ instead,
then we could always bring our model into one of the following standard forms,
where in general $a, b \in \C$ (see for example~\cite{Hulek}). When $n = 3$,
we can achieve this normal form even by a transformation from $\SL_3(\R)$,
so in this case, we can take $a, b \in \R$. We will call these forms
{\em Hesse forms}, generalising the classical terminology for $n = 3$.
They are as follows.
\begin{align*}
  n = 2: & \qquad y^2 = a (x_0^4 + x_1^4) + b\, x_0^2 x_1^2 \\
  n = 3: & \qquad a (x_0^3 + x_1^3 + x_2^3) + b\, x_0 x_1 x_2 = 0 \\
  n = 4: & \qquad \left\{\begin{array}{r@{\;=\; 0}}
                           a (x_0^2 + x_2^2) + b\, x_1 x_3 \\
                           a (x_1^2 + x_3^2) + b\, x_0 x_2
                         \end{array}\right.
\end{align*}
In these forms, the coefficients $a$ and~$b$ are bounded in terms of
the invariants, so we can expect them to be small. Therefore, we would
like to come close to a model of this kind, but using a unimodular
transformation.

We need some way of measuring how close two models are.  On the
standard Hesse models, the action of the $n$-torsion of the Jacobian,
$E[n]$, is given by the ``standard representation'' where one
generator multiplies each $x_j$ by~$\zeta_n^j$ and the other generator does
a cyclic shift of the coordinates. (Here $\zeta_n$ denotes a primitive
$n$th root of unity.)  To this representation, we can associate an invariant 
inner product on~$\C^n$, which is unique up to scaling. It
is easy to check that this invariant inner product is just the
standard one on~$\C^n$.  Now our approach is to associate an inner
product to a given model~$\CC$, and consider the model to be close to
a standard model when the associated inner product is close to the
standard one, which means that it is reduced in an appropriate
sense. This is explained in some detail in the following section.

\subsection{The reduction covariant}

Let $K = \R$ or $\C$.  We write $\Y_n(K)$ for the set of all genus one
normal curves of degree $n$ defined over $K$, inside a fixed copy of
$\PP^{n-1}$. (If $n=2$ we consider double covers of $\PP^1$ instead.)
The difference between~$\Y_n(K)$ and $X_n(K)$ is that we now consider
actual curves in~$\PP^{n-1}$ (or the set of ramification points of
$\CC \to \PP^1$ when $n = 2$), instead of defining equations.

Let $\HH^+_n(\C)$ be the space of positive definite Hermitian $n
\times n$ matrices, and $\HH^+_n(\R)$ the space of positive definite
symmetric real $n \times n$ matrices. We can identify these spaces
with the spaces of positive definite Hermitian and real quadratic
forms in $n$~variables, respectively.  There are natural and
compatible (left) actions of $\SL_n(K)$ on $\Y_n(K)$ and $\HH^+_n(K)$
given by the canonical map $\SL_n(K) \to \PGL_n(K) =
\Aut(\PP^{n-1}_K)$ on the one hand and by $g \cdot M = \bar{g}^{-t} M
g^{-1}$ on the other hand (where $\gamma^{-t}$ denotes the inverse
transpose of the matrix~$\gamma$).  If we identify the matrix $M \in
\HH^+_n(K)$ with the quadratic or Hermitian form $Q(x) = \bar{x}^t M
x$, then the compatibility of the actions means that $(g \cdot Q)(g x)
= Q(x)$.

\begin{Theorem}
  \label{redthm}
  For each $n \ge 2$ there is a unique $\SL_n(\C)$-covariant map
  $$ \varphi_\C : \Y_n(\C) \to \HH^+_n(\C)/\R^\times_{>0} \,. $$ 
  This map is compatible with complex conjugation, and hence restricts
  to an $\SL_n(\R)$-covariant map
  $$ \varphi_\R : \Y_n(\R) \to \HH^+_n(\R)/\R^\times_{>0} \,. $$ 
\end{Theorem}

\begin{Proof}
  Let $\CC \to \PP^{n-1}$ be a genus one normal curve defined 
  over $\C$, with Jacobian~$E$. The action of $E[n]$ 
  on $\CC$ extends to $\PP^{n-1}$ and hence defines a group 
  homomorphism $\chi: E[n](\C) \to \PGL_n(\C)$. 
  Lifting to $\SL_n(\C)$ we obtain a diagram
  $$ \xymatrix{ 
       0 \ar[r] & \mu_n \ar[r] \ar@{=}[d] & \ar[d] \ar[r] H_n 
            & \ar[d]^\chi \ar[r] E[n](\C) &  0 \\
       0 \ar[r] & \mu_n \ar[r] & \ar[r] \SL_n(\C) & \ar[r] \PGL_n(\C)
            &  0 \enspace . \!\!\!\!\! }
  $$ The Heisenberg group $H_n$ is a non-abelian group of order $n^3$.
  It comes with a natural $n$-dimensional representation, called the
  Schr\"odinger representation, which is known to be irreducible
  (since it is equivalent to the standard representation mentioned
  above).  Now by the Weyl unitary trick, every irreducible complex
  representation of a finite group has a unique invariant inner
  product. (Recall that existence is proved by averaging over the
  group, and uniqueness (up to $\R^\times_{>0}$) using Schur's lemma.)

  We define $\varphi_\C(\CC)$ 
  to be the (matrix of the) Heisenberg invariant inner product,
  {\em i.e.}, $\varphi_\C(\CC)$ is uniquely determined up to
  positive real scalars by the property that 
    \[  \bar{h}^{-t} \varphi_{\C} ( \CC) h^{-1} = \varphi_{\C} ( \CC) \]
  for all $h \in H_n$. 
  If $g \in \SL_n(\C)$, then the Heisenberg groups $H_n$ and $H_n'$ 
  of $\CC$ and $g \cdot \CC$ are related by $H_n' = g H_n g^{-1}$.
  Then $g \cdot \varphi_\C(\CC) = \bar{g}^{-t} \varphi_\C(\CC) g^{-1}$ 
  is an $H_n'$-invariant 
  inner product, and so must be equal to $\varphi_\C(g \cdot \CC)$. 
  Hence $\varphi_\C$ is $\SL_n(\C)$-covariant. Moreover, since
  $H_n \subset \SL_n(\C)$, this choice of covariant is forced on us.
  The compatibility with complex conjugation is seen in the same way.
\end{Proof}

\begin{Remark} In general $\varphi_\R$ is not the only
$\SL_n(\R)$-covariant. However, it is if the points of $E[n]$ are
defined over $\R$, as happens in the case $n=2$ and $\Delta>0$,
cf.~\cite[Lemma 3.2]{SCred}.
\end{Remark}

In practical terms, we have the following corollary.

\begin{Corollary}
  \label{WeylAverage}
  Let $M_T \in \GL_n(\C)$ describe the action of $T \in E[n](\C)$ on 
  $\CC \to \PP^{n-1}$.
  Then the reduction covariant $\varphi_\C(\CC)$ is
  $$ \sum_{T \in E[n](\C)} \frac{1}{|\det M_T|^{2/n}} 
                           \overline{M}_T^{\,\,t} M_T.
  $$ 
\end{Corollary}

\begin{Proof}
  To get an invariant inner product, we can take any inner product
  and average over its orbit under the action of~$H_n$. Applying this
  to the standard inner product, we find that we can take, up to scaling,
  \begin{equation} \label{eqncov}
    \varphi_{\C}(\CC) = \sum_{h \in H_n} \bar{h}^{-t} h^{-1}
                       = \sum_{h \in H_n} \bar{h}^t h \,.
  \end{equation}
  In the statement of the corollary, $M_T \in \GL_n(\C)$ is any lift
  of the element $\tau_T \in \PGL_n(\C)$ describing the action of~$T$
  on~$\PP^{n-1}(\C)$. The various pre-images of~$\tau_T$ in~$H_n$
  are given by $h = \alpha^{-1} M_T$ where $\alpha \in \C$ with
  $\alpha^n = \det M_T$. We then have
  \[ \bar{h}^t h 
       = \bar{\alpha}^{-1} \alpha^{-1} \overline{M}_T^{\,\,t} M_T
       = \frac{1}{|\det M_T|^{2/n}} \overline{M}_T^{\,\,t} M_T.
  \]
  Since this only depends on~$T$, it is sufficient to take the
  sum in~\eqref{eqncov} just over $T \in E[n](\C)$, 
  instead of over~$h \in H_n$.
\end{Proof}

We can now define what we mean by a reduced genus one normal curve.

\begin{Definition}
  A genus one normal curve $\CC \to \PP^{n-1}$ defined over $\R$ 
  is {\em Minkowski} (respectively {\em LLL}) {\em reduced}
  if $\varphi_\R(\CC)$ is the 
  Gram matrix of a Minkowski (respectively LLL) reduced lattice basis.
\end{Definition}

Note that a lattice basis is (Minkowski or LLL) reduced if it is close
to the standard basis of the standard lattice in the sense that the
basis vectors are (short and) nearly orthogonal. The notion of a
Minkowski reduced model has nice theoretical properties (it is optimal
and essentially unique), whereas for practical purposes, it is important
to be able to compute a reduced lattice basis efficiently; this is
possible when using LLL reduced models.

If we start with some given (minimal) model $\CC \to \PP^{n-1}$,
then in order to reduce it, we first compute its reduction 
covariant~$\varphi_{\R}(\CC)$. We apply the LLL algorithm~\cite{LLL} 
to this Gram matrix, resulting in a
unimodular transformation~$U$ and an LLL reduced Gram matrix~$M$,
such that $M = U^t \varphi_{\R}(\CC) U = U^{-1} \cdot \varphi_{\R}(\CC)$.
We then apply the transformation~$U^{-1}$ to our model~$\CC$. Since
$\varphi_{\R}(\CC)$ is a covariant, we will have that
$\varphi_{\R}(U^{-1} \cdot \CC) = M$ is LLL reduced.
Therefore $U^{-1} \cdot \CC$ is
the (minimal and) reduced model we are looking for.

\medskip

In the following sections we discuss how to compute~$\varphi_\R$.
There are two basic approaches. One is to find the hyperosculating
points of $\CC(\C)$ numerically and to compute the covariant from them.
If $n = 2$, we are looking for the ramification points of the
covering $\CC \to \PP^1$; if $n = 3$, for the flex points of the
plane cubic curve $\CC \subset \PP^2$. The other approach is to use
the $n$-torsion points in~$E(\C)$ instead and compute their action
on~$\PP^{n-1}$. Generally speaking, the first approach leads to simpler
formulas, whereas the second approach tends to be numerically more
stable.


\subsection{Reduction of 2-coverings}

We identify $\HH^+_2(\R)$ with the space of real positive
definite binary quadratic forms, and $\HH^+_2(\R)/\R^\times_{>0}$ with
the upper half plane. This identification maps a real positive definite
binary quadratic form to its unique root in the upper half plane.

\subsubsection{Using the ramification points} \label{Red2ram}

Let $F(x, z) \in \R[x,z]$ be homogeneous of degree~$4$. 
We assume that $f(X) = F(X, 1)$ has degree~4 as well.
(If the leading coefficient is zero, make a change of coordinates first.) 
Let $\theta_1, \ldots, \theta_4 \in \C$ be the roots of~$f$. It is shown in 
\cite{SCred} that $\varphi_\R$ is given by 
$$ \varphi_\R(F)(x, z)
    = \sum_{i=1}^4 \frac{1}{|f'(\theta_i)|}
                      (x - \theta_i z)(x - \overline{\theta_i} z) \,.
$$
This goes back to Julia's thesis~\cite{Julia}, where three different formulas
are given according to the number of real roots of~$f$; see
also~\cite{CremonaReduction}.

The formula is still valid for $\varphi_\C$, in the form
\[ \varphi_\C(F)(x, z)
     = \sum_{i=1}^4 \frac{1}{|f'(\theta_i)|} |x - \theta_i z|^2 \,.
\]
In practice one should first numerically compute the roots of
the resolvent cubic (which is not changed by reduction)
and then compute the roots of~$f$ from these.

\subsubsection{Using the $2$-torsion of~$E$} \label{Red23}

The binary quartic 
\[ F(x,z)= a x^4 + b x^3 z + c x^2 z^2 + d x z^3 + e z^4 \]
has invariants $I$ and $J$ (see Section~\ref{g1mod}) and 
resolvent cubic $r(X) = X^3 - 3 I X + J$. For $\varphi$ a
root of~$r$ we set
$$ \begin{array}{rcl}
     \alpha_1(\varphi) & = & 4 a \varphi - 8 a c + 3 b^2 \\
     \alpha_2(\varphi) & = & b \varphi - 6 a d + b c \\
     \alpha_3(\varphi) & = & (-2 \varphi^2 + 2 c \varphi - 9 b d + 4 c^2)/3  
   \end{array}
$$
and 
$$ W = \begin{pmatrix} 0 & -1 \\ 1 & 0 \end{pmatrix} \,\, ,  
       \quad \text{  } \quad 
   A_\varphi = \begin{pmatrix} \alpha_1(\varphi) & \alpha_2(\varphi) \\
                               \alpha_2(\varphi) & \alpha_3(\varphi) 
               \end{pmatrix}.
$$

\begin{Lemma} \label{Tmatrix2}
  If $\alpha_1(\varphi) \neq 0$, then the action of
  the corresponding point $T \in E[2]$ on~$\PP^1$ is given by
  \[ M_T = W \! A_\varphi \,. \]
\end{Lemma}

\begin{Proof}
  Let $H(x,z)$ be the Hessian
  of~$F$. The pencil spanned by $F$ and~$H$
  contains three degenerate quartics: for each root~$\varphi'$ of the
  resolvent cubic, we have
  $$ \alpha_1(\varphi') \bigl(4 \varphi' F(x,z) - \tfrac{1}{3} H(x,z)\bigr)
       = \big(\alpha_1(\varphi')x^2 + 2 \alpha_2(\varphi') xz
              + \alpha_3(\varphi') z^2\big)^2 \,.
  $$
  Since the action of~$T$ leaves both $F$ and~$H$ invariant, $M_T$
  must induce an involution on~$\PP^1$ that either fixes or swaps the
  roots of the quadratic on the right hand side; there is exactly one
  root $\varphi'$ such that the roots of the corresponding quadratic are
  fixed. Therefore $\varphi' = \varphi$, and the lemma follows by
  checking that $W\!A_\varphi$ does indeed fix the roots of the relevant
  quadratic.
\end{Proof}

\begin{Lemma} \label{itvanishes}
  If $M_T \in \GL_2$ describes the action of $T \in E[2]$
  on $\CC \to \PP^{1}$ then 
  \begin{equation} \label{eqnitvanishes}
    \sum_{T \in E[2]} \frac{1}{\det M_T} M_T^{t} M_T = 0. 
  \end{equation}
\end{Lemma}

\begin{Proof}
  We can verify this generically using the formula of Lemma~\ref{Tmatrix2}.
\end{Proof}

\begin{Proposition}
  Let $F \in \R[x,z]$ be a non-singular binary quartic, with resolvent
  cubic $r(X) = X^3 - 3 I X + J$. 
  \begin{enumerate}
    \item If $\Delta(F) > 0$ then the reduction covariant is 
          $\pm A_\varphi$ where $\varphi$ is the unique root of~$r$ 
          with $\det (A_\varphi) > 0$ and the sign is that 
          of~$\alpha_1(\varphi)$.
    \item If $\Delta(F) < 0$ then the reduction covariant is 
          $$ \Real \bigg( \frac{1}{|\det A_\varphi|}
                           \overline{A}_\varphi A_\varphi 
                          - \frac{1}{\det A_\varphi} A_\varphi^2 \bigg)
          $$
          where $\varphi$ is a complex root of~$r$. 
  \end{enumerate}
\end{Proposition}

\begin{Proof}
  If $\Delta(F) > 0$, then $r$ has three real roots. Since
  $\det(A_\varphi) = - \alpha_1(\varphi) r'(\varphi)/3$, the
  analysis in~\cite{CremonaReduction} shows that there is a unique root 
  $\varphi$ of~$r$ with $\det (A_\varphi) > 0$ (in particular,
  $\alpha_1(\varphi) \neq 0$). By Lemmas \ref{Tmatrix2} and~\ref{itvanishes}
  the reduction covariant simplifies (up to a factor of~$2$) to 
  $$ \sum_{ T \in E[2], \, \det M_T >0 }
        \frac{1}{\det M_T} M_T^{t} M_T
       = I_2 + \frac{1}{\det A_\varphi} A_\varphi^2
       = \frac{\tr A_\varphi}{\det A_\varphi} A_\varphi \,,
  $$
  by the Cayley-Hamilton theorem. So $\pm A_\varphi$ is the positive
  definite symmetric matrix we are looking for, with the sign 
  that makes the top left entry positive.

  If $\Delta(F) < 0$, then $r$ has a pair of complex conjugate roots, 
  say $\varphi$ and $\overline{\varphi}$. If $E[2] = \{0,S,T,\overline{T}\}$,
  then we can take $M_S = M_T \overline{M}_T$, so
  $\det(M_S) = |\det(M_T)|^2 > 0$. By Lemmas \ref{Tmatrix2} 
  and~\ref{itvanishes} again, the reduction covariant simplifies to 
  $$ \Real \bigg( \frac{1}{|\det M_T|} \overline{M}_T^{\,t} M_T 
                   -  \frac{1}{\det M_T} M_T^{t} M_T \bigg)
      = \Real \bigg( \frac{1}{|\det A_\varphi|}
                        \overline{A}_\varphi A_\varphi 
                       - \frac{1}{\det A_\varphi} A_\varphi^2 \bigg) \,.
  $$ 
Notice that we cannot have $\alpha_1(\varphi) =
\alpha_1(\overline{\varphi}) = 0$, since then the resolvent cubic
would have a repeated root, contradicting the fact that $F$ is
non-singular.
\end{Proof}

\subsubsection{The cross terms}

So far, we have shown how to find a unimodular transformation
of the coordinates on~$\PP^1$ that reduces the 2-covering.
(If we start with a generalised binary quartic $(P, Q)$
then we work with $F = P^2 + 4Q$.)
There is still an ambiguity coming from the possibility of making
a $y$-substitution in the general form of a 2-covering.
The most reasonable convention seems to be to arrange that the
cross term coefficients $l,m,n$ are $0$ or~$1$. 
 

\subsection{Reduction of 3-coverings}

\subsubsection{Using the flex points}

Let $F(x,y,z) \in \R[x,y,z]$ be a nonsingular ternary cubic.
In order to find its reduction covariant (as
a positive definite quadratic form $Q(x,y,z)$), we proceed as follows.
Let $H(x,y,z)$ be the Hessian of~$F$ as defined in Section~\ref{g1mod}. 
Then the intersection of~$F = 0$ and~$H = 0$
consists of nine distinct points, the flex points of~$F$.
Three of them are real, the others come in three complex
conjugate pairs. 

There are twelve lines each containing
three of the flex points, coming in four triples of lines
that do not meet in a flex point. (These triples are the
``syzygetic triangles'' mentioned in Section~\ref{Red32} below.)
One of these triples has all
three lines real, call them $L_{11}$, $L_{12}$, $L_{13}$. Another one
has one line real, call it $L_{21}$, and two complex conjugate
lines, call them $L_{22}$ and $L_{23}$. Then $Q$ spans
the one-dimensional intersection of the spaces spanned by
$L_{11}^2$, $L_{12}^2$ and $L_{13}^2$, and by $L_{21}^2$ and 
$L_{22} L_{23}$, respectively.

In order to see why this recipe works, first observe that it clearly
defines an $\SL_3(\R)$-covariant map. We can always make an
$\SL_3(\R)$-transformation to bring~$F$ into the standard Hesse form
$$ F(x,y,z) = a(x^3 + y^3 + z^3)  +  b\,x y z.$$
Then $L_{11}, L_{22}, L_{33}$ are $x$, $y$, $z$, and $L_{21}$, $L_{22}$,
$L_{23}$ are $x+y+z$, $x+ \zeta_3 y + \zeta_3^2 z$, 
$x+ \zeta_3^2 y + \zeta_3 z$ (where $\zeta_3$ is a primitive cube
root of unity). One then looks at the intersection 
$$ \langle x^2, y^2, z^2 \rangle
     \cap \langle (x+y+z)^2, x^2 + y^2 + z^2 - xy - yz - zx \rangle
$$
and finds it is one-dimensional, spanned by $x^2 + y^2 + z^2$,
which is the reduction covariant of any $F$ in Hesse form.

The only way we know to implement this method in practice is 
by numerically solving for the flex points. If the given model
is far from reduced, then usually several of the flex points 
are very close to one another, which makes the computation of the
lines difficult. Another practical problem is that the two 
spaces of quadrics we compute are only approximate and therefore
will usually not have nontrivial intersection. 

\subsubsection{Using the $3$-torsion on~$E$} \label{Red32}

This is the method described in \cite[\S9.5]{ANTS}.
Let $F(x,y,z)$ be a ternary cubic with invariants $c_4$ and $c_6$ 
and Hessian $H$ as defined in~Section~\ref{g1mod}. 
Let $T = (x_T,y_T)$ be a 3-torsion point on the Jacobian
$$E: \quad y^2 = x^3 - 27 c_4 x - 54 c_6.$$ 
Then the 
cubic $\T(x,y,z) = 2 x_T F - 3 H$ is the product of 3 linear
forms. (In \cite[II.7]{Hilbert} it is 
called a ``syzygetic triangle''.) Making a change of coordinates (if 
necessary) we may suppose $\T(1,0,0) \not = 0$. 
We label the coefficients
$$ \T(x,y,z) = r x^3 + s_1 x^2 y + s_2 x y^2  + s_3 y^3 + t_1 x^2 z  
                + t_2 x z^2 + t_3 z^3 + u x y z + v y^2 z + w y z^2.
$$

The proof of \cite[Theorem 7.1]{ANTS} describes how to compute 
a formula for $M_T$, where the entries are polynomials 
in $r,s_1,s_2, \ldots ,w$ and $y_T$.
Up to a scaling, this works out as $M_T = r A + 2 y_T B$ where
{\small
$$ \begin{array}{rcl}
A_{11} & = &  -12 r s_2 w - 36 r s_3 t_2 + 12 r u v + 4 s_1^2 w 
  + 4 s_1 s_2 t_2  - 8 s_1 t_1 v - s_1 u^2 + 12 s_3 t_1^2 \\
A_{12} & = & -54 r s_3 w + 18 r v^2 + 6 s_1 s_2 w - 3 s_1 u v 
  - 6 s_2 t_1 v + 9 s_3 t_1 u \\
A_{13} & = &    -81 r s_3 t_3 + 9 r v w + 9 s_1 s_2 t_3 - 3 s_1 t_2 v
  - 3 s_2 t_1 w + 9 s_3 t_1 t_2 \\ 
A_{21} & = &    36 r s_2 t_2 - 9 r u^2 - 12 s_1^2 t_2 + 12 s_1 t_1 u 
  - 12 s_2 t_1^2 \\ 
A_{22} & = &    24 r s_2 w + 18 r s_3 t_2 - 15 r u v - 8 s_1^2 w 
  - 2 s_1 s_2 t_2 + 10 s_1 t_1 v + 2 s_1 u^2 - 3 s_2 t_1 u - 6 s_3 t_1^2 \\ 
A_{23} & = &    54 r s_2 t_3 - 9 r u w - 18 s_1^2 t_3 + 6 s_1 t_1 w 
  + 3 s_1 t_2 u - 6 s_2 t_1 t_2 \\ 
A_{31} & = &    0 \\ 
A_{32} & = &    -18 r s_2 v + 27 r s_3 u + 6 s_1^2 v - 3 s_1 s_2 u 
  - 18 s_1 s_3 t_1 + 6 s_2^2 t_1 \\ 
A_{33} & = &    -12 r s_2 w + 18 r s_3 t_2 + 3 r u v + 4 s_1^2 w 
  - 2 s_1 s_2 t_2 - 2 s_1 t_1 v - s_1 u^2 + 3 s_2 t_1 u - 6 s_3 t_1^2 
\end{array} $$
}
and $B = r B_1 + (s_1^2 t_2 - s_1 t_1 u + s_2 t_1^2) E_{13}$ with
$$ B_1 = \begin{pmatrix} 
    s_1 u - 2 s_2 t_1 & 
    s_1 v - 3 s_3 t_1 &
    s_1 w - 4 s_2 t_2 - t_1 v + u^2 \\
    -3 r u + 2 s_1 t_1 &
    -3 r v + s_2 t_1 &
    -3 r w + s_1 t_2 \\
    6 r s_2 - 2 s_1^2 &
    9 r s_3 - s_1 s_2 &
    3 r v - s_1 u + s_2 t_1
\end{pmatrix}. $$
(Notes: $E_{ij}$ is the 3 by 3 matrix with $(i,j)$ entry 1 and all other
entries $0$. Our matrices $A$ and $B$ would be called $r^3 (\det P)A$ and 
$r^3 B$ in the notation of \cite{ANTS}.)
This formula comes with the caveat (see \cite[Remark 7.2]{ANTS}) 
that it may give zero.
However, this will never happen for both $T$ and $-T$, so we 
get round the problem by computing $M_T$ as $(M_{-T})^{-1}$. 

Once we have computed $M_T$ for all $T \in E[3]$
the reduction covariant is computed using Corollary~\ref{WeylAverage}.


\subsection{Reduction of 4-coverings}
\label{sec:R4}

We could again try to find the reduction covariant starting from
the 16~hyperosculating points on~$\CC$ and the quadruples of planes
containing four of them, which are the analogue of the syzygetic
triangles. However, this approach does not seem to be very promising.

Instead, we use the fact that below the given 4-covering~$\CC$, there
is a 2-covering $\CC_2$; let $\pi : \CC \to \CC_2$ be the covering map.
If $A$ and~$B$ 
are the symmetric matrices corresponding to the
quadrics defining $\CC \subset \PP^3$, then~$\CC_2$
has equation $y^2 =F(x,z)$ where
\[ F(x,z) := \det(A x +  Bz) \,. \]
Applying reduction to the quartic on the right hand side, we find a
good basis of the pencil of quadrics. It remains to find the reduction
covariant of~$\CC$.

Let $\theta_j \in \C$ ($j = 1,2,3,4$) be the ramification points of 
$\CC_2 \to \PP^1$, {\em i.e.}, the roots of~$f(X) = F(X,1)$. 
Let $G_j$ ($j = 1,2,3,4$) be a linear form (unique up to scaling)
describing the preimage
of~$\theta_j$ on~$\CC \subset \PP^3$. Then (fixing the polynomials
giving the covering map~$\pi : \CC \to \CC_2$) there are 
$\alpha_j \in \C^\times$ such that
\[ (x - \theta_j\,z) \circ \pi = \alpha_j G_j^2 \,. \]
Now the action of $T \in E[4]$ on~$\CC$ induces the action of $2T \in E[2]$
on~$\CC_2$. Therefore the action of $T \in E[2]$ on~$\CC$ will be trivial
on~$\CC_2$, hence the corresponding matrix $M_T \in \SL_4$ will fix
the~$G_j$ up to sign. In fact, it can be checked that the action
of~$E[2]$ on~$\PP^3$ lifts to a representation on~$\C^4$, which is
isomorphic to the regular representation, and the $G_j$ span the four
eigenspaces. So any Hermitian form that is invariant under~$H_4$ must
be invariant under~$E[2]$ and thus be of the form
\[ \sum_{j=1}^4 \lambda_j |G_j|^2 \,. \]
It remains to determine the coefficients~$\lambda_j$.

\begin{Lemma} \label{Red4}
  Keep the notation introduced so far, and let $f(X) = F(X,1)$.
  Then the reduction covariant of~$\CC$ is the positive definite
  Hermitian form
  \[ \varphi_{\C}(\CC) =
     \sum_{j=1}^4 \frac{|\alpha_j|}{|f'(\theta_j)|^{1/2}} |G_j|^2 \,. \]
\end{Lemma}

If $\CC$ is defined over~$\R$, then the restriction of this Hermitian
form to~$\R^4$ 
will be the positive definite quadratic form~$\varphi_{\R}(\CC)$.

\medskip

\begin{Proof}
  We first check that the given form is invariant under $\SL_2(\C)$
  acting on~$\PP^1$ ({\em i.e.}, does not depend on the choice of
  basis of the pencil of quadrics). We know (see Section~\ref{Red2ram} above)
  that $\sum_{j=1}^4 |f'(\theta_j)|^{-1} |x - \theta_j z|^2$
  is an $\SL_2(\C)$-covariant; the same computation (which deals with
  each summand separately) shows that 
  $\sum_{j=1}^4 |f'(\theta_j)|^{-1/2} |x - \theta_j z|$
  is a covariant as well. But $|x - \theta_j z| = |\alpha_j G_j^2|$,
  and the coordinates in~$G_j$ are not affected by the $\SL_2(\C)$-action,
  so the expression given in the statement is invariant.
  
  Now we check that the given form is covariant with respect to the
  action of~$\SL_4(\C)$. But this is clear since every $\alpha_j G_j^2$
  is covariant.
  
  Since we can move any~$\CC$ into standard form by the action of
  $\SL_2(\C) \times \SL_4(\C)$, it now suffices to verify that our
  formula gives the correct result when $\CC$ is in standard form
  \[ a (x_0^2 + x_2^2) + 2b\,x_1 x_3
      = a (x_1^2 + x_3^2) + 2b\,x_0 x_2 = 0 \,.
  \]
  In this case, the 2-covering~$\CC_2$ is given by
  \[ y^2 = (a^4 + b^4) x^2 z^2 - a^2 b^2 (x^4 + z^4) \]
  and the map $\pi$ (see Lemma~\ref{4to2cover} for formulae), 
followed by the map $\CC_2 \to \PP^1$, is given by
  \[ (x : z) = \bigl(b^3(x_1^2 + x_3^2) + 2a^3\,x_0 x_2
                     : -b^3(x_0^2 + x_2^2) - 2a^3\,x_1 x_3\bigr) \,.
  \]
  The roots $\theta_j$ of $f(X) = -a^2 b^2 X^4 + (a^4 + b^4) X^2 - a^2 b^2$
  are $a/b$, $-a/b$, $b/a$, $-b/a$, and up to a common factor $b^4-a^4$,
  we can take $\alpha_j = 1/b$, $1/b$, $1/a$, $-1/a$ and
  $G_j = x_1-x_3$, $x_1+x_3$, $x_0-x_2$, $x_0+x_2$. Also,
  $|f'(\theta_j)| = c |\theta_j|$ for some constant~$c$. 
  Since $|\alpha_j|/|\theta_j|^{1/2}$ has the
  same value $|ab|^{-1/2}$ for all~$j$, our expression gives,
  up to a constant factor again,
  \[ |x_1 - x_3|^2 + |x_1 + x_3|^2 + |x_0 - x_2|^2 + |x_0 + x_2|^2
       = 2\bigl(|x_0|^2 + |x_1|^2 + |x_2|^2 + |x_3|^2) \,,
  \]
  which is the correct result for a 4-covering in standard form.
\end{Proof}

In order to find the $\alpha_j$ and~$G_j$, we can make use of a
result from~\cite{ANTS8}, where it is observed that $\alpha_j G_j^2$
is the quadratic form corresponding to the matrix
\[ e \theta_j^{-1} A + M_1 + \theta_j M_2 + a \theta_j^2 B \,; \]
here $F(x,z) = \det(A x + B z) = a x^4 + b x^3z + c x^2z^2 + d x z^3+ e z^4$
and $M_1$, $M_2$ are obtained from the relation~(\ref{defnT1T2})
in the proof of Lemma~\ref{4to2cover}.


\section{Examples}
\label{sec:E}

In this section we illustrate minimisation and reduction for two
explicit examples over $\Q$ (one a $3$-covering and the other a
$4$-covering).  We then give references to further examples.

\subsection{Minimisation and reduction of a 3-covering} \label{Sec:Ex1}

We consider the elliptic curve 105630d1 in \cite{CrTables} with
Weierstrass equation
\[ E: \qquad y^2 + x y = x^3 + x^2 - 114848533 x - 472424007827. \]
Computing the $3$-Selmer group (see \cite{SchaeferStoll}) we find 
$\Sel^{(3)}(\Q,E) \cong \Z/3\Z$. In~\cite{ndescent} we show how to 
write down elements of the $3$-Selmer group explicitly as $3$-coverings of $E$.
In this case our \Magma programs find (before minimisation and
reduction) that a generator is represented by the $3$-covering 
$\CC \subset \PP^2$ with equation
\begin{align*} F_1(x,y,z) &= 27089 x^3 + 2142 y^3 + 291938 z^3 + 10008 x^2 y - 127341 x^2 z  \\ &~ \qquad  
  + 92937 x y^2 + 104736 y^2 z + 21093 x z^2 - 71172 y z^2 - 2655 x y z.
\end{align*}
(Random choices in the programs mean it need not return the same cubic
every time. However, the answer will always
be $\Q$-equivalent to $F_1$, and this can be checked using the
algorithm in \cite{ANTS}.)
The discriminant of this ternary cubic is $\Delta(F_1) = 
3^{12} \cdot 503^{12}  \cdot \Delta_E$ where 
$\Delta_E = 2^{39} \cdot 3  \cdot 5^{9} \cdot 7^{3}  \cdot 503$
is the minimal discriminant of $E$. So $F_1$ has level
$1$ at the primes $3$ and $503$. Reducing mod $3$ we find
$F_1(x,y,z) = 2(x+z)^3 \pmod{3}$. The level is decreased by the first
iteration of our algorithm (see Theorem~\ref{mintc}). Explicitly we put
 \[ F_2(x,y,z) = \frac{1}{3^2} F_1(3x-y,z,y). \]
Likewise we find  $F_2(x,y,z) \equiv 284(x + 329 y + 33 z)^3 \pmod{503}$ 
and our algorithm puts
\begin{align*} F_3(x,y,z) &= \frac{1}{503^2} F_2(503 x - 33 y+z,z,y - 10 z) \\
  &= 40877301 x^3 - 11504 y^3 + 12 z^3 - 8035425 x^2 y - 64887 x^2 z \\ &~ \qquad \quad + 526580 x y^2  - 200 y^2 z  + 5803 x z^2 
  - 383 y z^2  + 7307 x y z. \end{align*}
The $3$-torsion of $y^2 = x^3 - 27 c_4 x - 54 c_6$ over $\C$ is generated by
\begin{align*} 
S & = (667989.968057, 420236746.168), &
T & = (-264330.994609, 34120617.5970 i). 
\end{align*}
The formulae in Section~\ref{Red32} show that $S$ and $T$ act
on $\{F_3=0\}$ via
\[ M_S = \begin{pmatrix}
  285.46 & -19.022 & 3.4264 \\ 
  4352.6 & -290.04 & 52.341 \\ 
  509.05 & -33.785 & 4.5806 
\end{pmatrix} \]
and 
\[ M_T = \begin{pmatrix}
-50.656 + 47.060 i & 3.2758 - 3.3464 i & 0.11909 + 2.2683 i \\ -786.55 + 
717.15 i & 50.871 - 51.000 i & 1.8675 + 34.587 i \\ -119.84 + 93.073 i & 
7.8268 - 6.5354 i & -0.21547 + 3.9405 i 
\end{pmatrix}. \]
We have scaled these matrices to have determinant $1$.
By Corollary~\ref{WeylAverage} the reduction covariant has
matrix
\[ A = \begin{pmatrix}
 176413988.185 & -11560848.1174 & 3471.84429193 \\ 
-11560848.1174 & 757736.524016 & -1499.92503970 \\
 3471.84429193 & -1499.92503970 & 13237.5156939 
\end{pmatrix}.
\]
Running the LLL algorithm on the lattice with Gram matrix $A$ results
in the unimodular transformation.
\[
U = \begin{pmatrix}
  0  & 0 & 1 \\
  4  & 61 &  6 \\
 -3 & -46 & -4 
\end{pmatrix}.
\]
Accordingly we put $F_4(x,y,z) = F_3(4 y - 3 z,61 y - 46 z,x + 6 y - 4 z)$
and find
\begin{align*} F_4(x,y,z)  &=
12 x^3 + 12 y^3 + 171 z^3 + 65 x^2 y + 65 x^2 z \\
 & \qquad {}         - 94 y^2 z + 87 x z^2 + 101 y z^2 + 7 x y z.
\end{align*}

This ternary cubic has solution 
\[(x:y:z) = ( 345420 : -1638959 : -373029 )\]
which by the formulae in \cite{AKM3P} maps down to a point
\begin{align*}
x & = \frac{-74872620773608422623058757914981065217}{109435039457696221^2} \\[2mm]
y &=  \frac{51043047025320389176098494307847798722958228061916407587}%
           {109435039457696221^3}
\end{align*}
on $E(\Q)$ of canonical height $86.5313\ldots$.
Since the torsion subgroup of $E(\Q)$ is trivial, 
it follows that $\rank E(\Q)=1$. It is equally convenient 
to find this generator using Heegner points.

Note that the \Magma implementation of 3-descent does the minimisation
and reduction automatically. 
To extract the intermediate model $F_1(x,y,z) = 0$, 
one should first specify that 3-descent prints out some of its working, 
using the command {\tt SetVerbose("ThreeDescent",1);}

\subsection{Minimisation and reduction of a 4-covering}

In~\cite[\S8.1]{Skoro}, an example is given of a 4-covering~$\CC$ of
the elliptic curve $E : y^2 = x^3 - 1221$ that represents an element
of exact order~4 in the Shafarevich-Tate group of~$E$. The symmetric
matrices corresponding to the two quadrics defining $\CC \subset \PP^3$
are given as (to keep with our convention, we multiply by~$2$ so that
entries are the second partial derivatives)
\[ A = 2 \begin{spmatrix}
           -1 & 11 & -66 & 396 \\
           11 & -66 & 396 & -2520 \\
           -66 & 396 & -2520 & 16\,335 \\
           396 & -2520 & 16\,335 & -105\,786
         \end{spmatrix}
         \quad\text{and}\quad
   B = 2 \begin{spmatrix}
           -1 & -3 & 33 & -198 \\
           -3 & 33 & -198 & 1188 \\
           33 & -198 & 1188 & -7560 \\
           -198 & 1188 & -7560 & 49\,005
         \end{spmatrix}
   \,.
\]
We will use $x_1, \dots, x_4$ as the coordinates on~$\PP^3$.
We find that
\[ \det(Ax + Bz) = 2^4 \cdot 3^8 (-9 x^4 + 13 x^3 z - 18 x^2 z^2 + 3 z^4) \,, \]
which makes it clear that the model is non-minimal at~$p = 2$ and $p = 3$.
We compute that the discriminant of our quadric intersection is $(2 \cdot 3^4)^{12}$
times the (minimal) discriminant $-2^4 \, 3^5 \, 11^2 \, 37^2$ of~$E$, which shows that
the level at~$2$ is~$1$ and the level at~$3$ is~$4$; the model is already
minimal at all other primes.

We first minimise at~$p = 3$. According to our algorithm (see Section~\ref{M4}),
we have to look at the reductions of $A$ and~$B$ mod~$3$, which are
\[ \bar{A} = \begin{spmatrix}
               1 & 1 & 0 & 0 \\ 1 & 0 & 0 & 0 \\ 0 & 0 & 0 & 0 \\ 0 & 0 & 0 & 0
             \end{spmatrix}
             \quad\text{and}\quad
   \bar{B} = \begin{spmatrix}
               1 & 0 & 0 & 0 \\ 0 & 0 & 0 & 0 \\ 0 & 0 & 0 & 0 \\ 0 & 0 & 0 & 0
             \end{spmatrix}
\]
The common nullity is $s = 2$, and the reduced quadratic forms already
involve only the first two variables. They represent zero simultaneously
over~$\F_3$; the plane $x_1 = 0$ is contained in the reduction of the curve.
So we apply the transformation $[\frac{1}{3} I_2, \Diag(3,1,1,1)]$, resulting
in the new pair of matrices (which we will again denote $A$ and~$B$)
\[ A = \begin{spmatrix}
            -6 &    22 &  -132 &   792 \\
            22 &   -44 &   264 & -1680 \\
          -132 &   264 & -1680 & 10890 \\
           792 & -1680 & 10890 &-70524 \\
       \end{spmatrix}
       \quad\text{and}\quad
   B = \begin{spmatrix}
           -6 &   -6 &   66 & -396 \\
           -6 &   22 & -132 &  792 \\
           66 & -132 &  792 & -5040 \\
         -396 &  792 &-5040 & 32670 \\
       \end{spmatrix}
\]
The level at~$p = 3$ of the new model is~$3$. Reducing mod~3, we have now
\[ \bar{A} = \begin{spmatrix}
               0 & 1 & 0 & 0 \\ 1 & 1 & 0 & 0 \\ 0 & 0 & 0 & 0 \\ 0 & 0 & 0 & 0
             \end{spmatrix}
             \quad\text{and}\quad
   \bar{B} = \begin{spmatrix}
               0 & 0 & 0 & 0 \\ 0 & 1 & 0 & 0 \\ 0 & 0 & 0 & 0 \\ 0 & 0 & 0 & 0
             \end{spmatrix}
\]
The common nullity is again $s = 2$, and there is a plane contained in the
reduction. This time, the plane is $x_2 = 0$, so we swap $x_1$ and~$x_2$
before we apply $[\frac{1}{3} I_2, \Diag(3,1,1,1)]$. The result is a model
of level~2:
\[ A = \begin{spmatrix}
          -132 &    22 &   264 & -1680 \\
            22 &    -2 &   -44 &   264 \\
           264 &   -44 &  -560 &  3630 \\
         -1680 &   264 &  3630 &-23508
       \end{spmatrix}
       \quad\text{and}\quad
   B = \begin{spmatrix}
           66 &   -6 & -132 &   792 \\
           -6 &   -2 &   22 &  -132 \\
         -132 &   22 &  264 & -1680 \\
          792 & -132 &-1680 & 10890
       \end{spmatrix}
\]
Now we get a different situation mod~3:
\[ \bar{A} = \begin{spmatrix}
               0 & 1 & 0 & 0 \\ 1 & 1 & 1 & 0 \\ 0 & 1 & 1 & 0 \\ 0 & 0 & 0 & 0
             \end{spmatrix}
             \quad\text{and}\quad
   \bar{B} = \begin{spmatrix}
               0 & 0 & 0 & 0 \\ 0 & 1 & 1 & 0 \\ 0 & 1 & 0 & 0 \\ 0 & 0 & 0 & 0
             \end{spmatrix}
\]
The common nullity is $s = 1$. We swap $x_1$ and~$x_4$ so that the reduced forms
only involve the last three variables. Then we see that we are in `Situation~2',
so we apply the transformation $[I_2, \Diag(\frac{1}{3}, 1, 1, 1)]$. This
results in a model of level~$1$, given by
\[ A = \begin{spmatrix}
         -2612 &   88 & 1210 & -560 \\
            88 &   -2 &  -44 &   22 \\
          1210 &  -44 & -560 &  264 \\
          -560 &   22 &  264 & -132
       \end{spmatrix}
       \quad\text{and}\quad
   B = \begin{spmatrix}
         1210 & -44 & -560 &  264 \\
          -44 &  -2 &   22 &   -6 \\
         -560 &  22 &  264 & -132 \\
          264 &  -6 & -132 &   66
       \end{spmatrix}
\]
In the last minimisation step at~$p = 3$, the reductions are now
\[ \bar{A} = \begin{spmatrix}
               1 & 1 & 1 & 1 \\ 1 & 1 & 1 & 1 \\ 1 & 1 & 1 & 0 \\ 1 & 1 & 0 & 0
             \end{spmatrix}
             \quad\text{and}\quad
   \bar{B} = \begin{spmatrix}
               1 & 1 & 1 & 0 \\ 1 & 1 & 1 & 0 \\ 1 & 1 & 0 & 0 \\ 0 & 0 & 0 & 0
             \end{spmatrix}
\]
The common nullity is again $s = 1$, and the common kernel is spanned
by $(1, -1, 0, 0)$. We move it to $(1, 0, 0, 0)$ and are in `Situation~2'
again. After applying $[I_2, \Diag(\frac{1}{3}, 1, 1, 1)]$, we obtain
a model that is now minimal at~$p = 3$.
\[ A = \begin{spmatrix}
         -310 &  30 & 418 &-194 \\
           30 &  -2 & -44 &  22 \\
          418 & -44 &-560 & 264 \\
         -194 &  22 & 264 &-132
       \end{spmatrix}
       \quad\text{and}\quad
   B = \begin{spmatrix}
          144 & -14 &-194 &  90 \\
          -14 &  -2 &  22 &  -6 \\
         -194 &  22 & 264 &-132 \\
           90 &  -6 &-132 &  66
       \end{spmatrix}
\]

We still have to minimise at~$p = 2$, using the algorithm described in
Section~\ref{min_char2}. We first find the `double' of our model:
\begin{align*}
  \dd'(A,B) = (P,Q) = \bigl(&2^2 (6413 x^2 - 5665 x z + 1248 z^2), \\
                    &2^2 (41126578 x^4 - 72659303 x^3 z \\
                  &\quad {} + 48099091 x^2 z^2 - 14139840 x z^3 + 1557501 z^4)\bigr)
\end{align*}
We see that we already have $v_2(P) \ge 1$ and $v_2(Q) \ge 2$.
The common kernel of the reductions mod~2 of the two quadratic forms
is spanned by $(1, 1, 0, 1)$ and $(0, 0, 1, 0)$, so the common nullity
is $s = 2$. We change coordinates so that the common kernel is given
by $x_1 = x_2 = 0$. Then the reductions of the quadrics are $x_1^2$ and~$x_2^2$,
so they do not simultaneously represent zero. We apply the `flip-flop'
transformation $[\frac{1}{2} I_2, \Diag(2,2,1,1)]$, after which the
reductions are $x_3 x_4$ and $x_4^2$, so now there is the plane $x_4 = 0$
contained in the reduction of the curve. We swap $x_1$ and~$x_4$ and then
apply $[\frac{1}{2} I_2, \Diag(2,1,1,1)]$ to obtain a pair of matrices
representing a globally minimal model:
\[ A = \begin{spmatrix}
         -728 & -424 &  319 & -474 \\
         -424 & -252 &  187 & -280 \\
          319 &  187 & -140 &  209 \\
         -474 & -280 &  209 & -310
       \end{spmatrix}
       \quad\text{and}\quad
   B = \begin{spmatrix}
          348 &  198 & -152 &  220 \\
          198 &  114 &  -86 &  130 \\
         -152 &  -86 &   66 &  -97 \\
          220 &  130 &  -97 &  144
       \end{spmatrix}
\]

We now apply reduction to this model as described in Section~\ref{sec:R4}.
We have
\[ \det(Ax + Bz) = 4(-9 x^4 + 13 x^3 z - 18 x^2 z^2 + 3 z^4) \,. \]
Following \cite{AKM3P} and~\cite{ANTS8}, we compute the quadratic forms
$T_1$, $T_2$ whose symmetric matrices $M_1$, $M_2$ are given by
\[ \adj\bigl(\adj(A) x + \adj(B) z)\bigr)
      = 4^2 \cdot 81 A x^3 - 4 \cdot 9 M_1 x^2 z
          + 4 \cdot 3 M_2 x z^2 + 4^2 \cdot 9 B z^3 \,.
\]
Then, writing $Q_1$ and~$Q_2$ for the quadratic forms corresponding
to $A$ and~$B$,
\[ \alpha G^2 = 12 \theta^{-1} Q_1 + T_1 + \theta T_2 - 36 \theta^2 Q_2 \]
for $\theta$ a root of $f(X) = \det(X A + B)$. We can for example take
\[ G = (-18 \theta^3 - 28 \theta^2 + 6 \theta + 2) x_1
       + (18 \theta^3 - 26 \theta^2 + 2) x_2
       + (18 \theta^2 + \theta - 3) x_3
       - 2 x_4
\]
and $\alpha = -1395 \theta^3 + 1367 \theta^2 - 2155 \theta - 1001$.
Also, $f'(\theta) = 12(-12\theta^3 + 13\theta^2 - 12\theta$). The matrix
corresponding to $\sqrt{12} \sum_\theta |\alpha| |G|^2/|f'(\theta)|^{1/2}$ is
(to five decimal digits precision)
\[ \begin{pmatrix}
     8857.72019 &  5117.00780 & -3885.97776 &  5665.67630 \\
     5117.00780 &  3080.24124 & -2279.16858 &  3348.18401 \\
    -3885.97776 & -2279.16858 &  1716.07038 & -2498.36286 \\
     5665.67630 &  3348.18401 & -2498.36286 &  3706.96839
   \end{pmatrix} \,.
\]
We apply LLL to this Gram matrix and obtain the reducing transformation
matrix
\[ U = \begin{pmatrix}
          -5 & -2 &  -6 &  0 \\
          -6 & -3 &  -7 & -1 \\
         -15 & -7 & -17 &  0 \\
           3 &  1 &   4 &  1
       \end{pmatrix} \,,
\]
which finally brings the two matrices defining~$\CC$ into the form
\[ U^t A U =
   \begin{spmatrix}
     -2 &  0 & -1 & -2 \\
      0 & -2 & -1 &  0 \\
     -1 & -1 & -2 &  2 \\
     -2 &  0 &  2 & -2
   \end{spmatrix}
   \quad\text{and}\quad
   U^t B U =
   \begin{spmatrix}
      0 &  0 & -1 &  1 \\
      0 &  2 & -1 & -1 \\
     -1 & -1 &  0 & -1 \\
      1 & -1 & -1 & -2
   \end{spmatrix} \,.
\]
These correspond, after a sign change, to the quadratic forms
\begin{align*}
  Q_1 &= x_1^2 + x_1 x_3 + 2 x_1 x_4 + x_2^2 + x_2 x_3 + x_3^2 - 2 x_3 x_4 + x_4^2
    \quad\text{and}\quad \\
  Q_2 &= x_1 x_3 - x_1 x_4 - x_2^2 + x_2 x_3 + x_2 x_4 + x_3 x_4 + x_4^2 \,.
\end{align*}


\subsection{Further examples and applications}

One useful application of the methods described in this paper is to help
find large generators in the Mordell-Weil group of an elliptic curve~$E$.
This has already be demonstrated in Section~\ref{Sec:Ex1}.
Each rational point $P \in E(\Q)$ lifts to one of the $n$-coverings
of~$E$. If we have a nice and small ({\em  i.e.}, minimised and
reduced) model~$\CC$ of this $n$-covering, then the logarithmic
height with respect to $\CC \to \PP^{n-1}$ of the preimage~$Q$ of~$P$ in~$\CC(\Q)$
will be smaller by a factor of about $\frac{1}{2n}$ than the logarithmic
$x$-coordinate height of~$P$ --- standard properties of heights imply that
\[ h(Q) = \frac{1}{2n} h_x(P) + O(1) \]
where the implied constant depends on the equations defining $\CC \to \PP^{n-1}$.
If the equations have small coefficients, this constant should be small
as well. Therefore we can hope to find~$P$ much more easily by searching
for~$Q$ on~$\CC$. In fact, this application was the motivation for the
first tentative steps towards reduction of 4-coverings. The story begins
with~\cite{GPZ}, where the authors determined Mordell-Weil generators for all
Mordell curves $y^2 = x^3 + D$, with $D$ a nonzero integer of absolute
value at most~$10^4$ (in order to determine all the integral points on
these curves), with one exception, $D = 7823$. The analytic rank of
this curve is~$1$, so we know that the Mordell-Weil rank must be
also~$1$; however the Birch and Swinnerton-Dyer Conjecture predicts a
generator of fairly large height. One of us (Stoll) used minimisation
and reduction of 4-coverings in a fairly {\em ad hoc} fashion
to find a good model of the one relevant
4-covering of $E : y^2 = x^3 + 7823$, so that a point search on this
4-covering curve was successful, thus resolving this last open case.
The result was reported in a posting~\cite{StollMailing} to the
{\tt NMBRTHRY} mailing list. We give a short summary of the steps
and the result. By a standard 2-descent, one obtains a 2-covering curve
\[ C : y^2  =  -18 x^4 + 116 x^3 + 48 x^2 - 12 x  + 30 \,. \]
A second 2-descent on~$C$ following~\cite{MSS} produces a 4-covering
of~$E$, whose initial model was given by quadrics with coefficients
of up 15~decimal digits. Using the methods described here, one
finds a model $D \subset \PP^3$ given by
\begin{align*}
  2 x_1 x_2 + x_1 x_3 + x_1 x_4 + x_2 x_4 + x_3^2 - 2 x_4^2 &= 0 \\
  x_1^2 + x_1 x_3 - x_1 x_4 + 2 x_2^2 - x_2 x_3
   + 2 x_2 x_4 - x_3^2 - x_3 x_4 + x_4^2 &= 0
\end{align*}
It is not very difficult to find the point
$P = (116 : 207 : 474 : -332)$ on~$D$. This point then gives rise to the point
\[ Q = \left(\frac{53463613}{32109353}, 
             \frac{23963346820191122}{32109353^2}\right) \]
on~$C$, which in turn finally produces the Mordell-Weil generator on~$E$,
with coordinates
\begin{align*}
  x &= \frac{2263582143321421502100209233517777}{11981673410095561^2} \\[2mm]
  y &= \frac{186398152584623305624837551485596770028144776655756}%
            {11981673410095561^3}
\end{align*}
Note that in the version given in the mailing list posting, the model
was not minimal at~2 (in fact, it had level~2 at~2).

4-descent including minimisation and reduction was also used to find
some of the elliptic curves of high rank and prescribed torsion listed
in~\cite{DujellaTables}, for example the curve with
$E(\Q) \cong \Z/12\Z \times \Z^4$.

Minimised and reduced models of 2-, 3-, and 4-coverings provide the
starting point for the computation of 6- and 12-coverings as described
in~\cite{6and12}. These then allow us to find even larger generators
(of logarithmic canonical height $> 600$). For example, this method
was used to find the last missing generators for curves of prime
conductor and rank at least~2 in the Stein-Watkins database~\cite{StWat}.

A table giving representatives of all elements of order~3
in the Sha\-fa\-re\-vich-Tate groups of all elliptic curves of conductor
$< 130\,000$ can be found at~\cite{Sha3tables}. 
(It is only known that the table is complete if one assumes
the conjecture of Birch and Swinnerton-Dyer.) The final form of
these ternary cubics was obtained by applying the methods described
in this paper to the original models produced by the algorithms described
in~\cite{SchaeferStoll} and~\cite{ndescent}.



\begin{thebibliography}{MKM3P}

\frenchspacing
\renewcommand{\baselinestretch}{1}

\bibitem[AKM$^3$P]{AKM3P}
S.Y. An, S.Y. Kim, D.C. Marshall, S.H. Marshall,
W.G. McCallum and A.R. Perlis,
Jacobians of genus one curves,
{\em J. Number Theory} {\bf{90}} (2001), no. 2, 304--315. 

\bibitem[ARVT]{ARVT}
M. Artin, F. Rodriguez-Villegas and J. Tate, 
On the Jacobians of plane cubics,
{\em Adv. Math.} {\bf{198}} (2005), no. 1, 366--382.

\bibitem[BSD]{BSD1}
B.J. Birch and H.P.F. Swinnerton-Dyer,
Notes on elliptic curves I,  
{\em J. reine angew. Math.}  {\bf{212}} (1963) 7--25. 

\bibitem[BLR]{BLR}
S. Bosch, W. L\"utkebohmert and M. Raynaud,
{\em N\'eron models},
Ergebnisse der Mathematik und ihrer Grenzgebiete (3), {\bf{21}},
Springer-Verlag, Berlin, 1990.

\bibitem[Co]{Connell}
I. Connell, {\em Elliptic Curve Handbook}, (unpublished on-line notes), 
McGill University, 1996. 

\bibitem[Cr1]{Cremona} 
J.E. Cremona, 
{\em Algorithms for modular elliptic curves}, Second edition,
Cambridge University Press, Cambridge, 1997. 

\bibitem[Cr2]{CremonaReduction} 
J.E. Cremona, 
Reduction of binary cubic and quartic forms,
{\em LMS J. Comput. Math.} {\bf 2} (1999), 64--94 (electronic).

\bibitem[Cr3]{CrTables} 
J.E. Cremona, Tables of Elliptic Curves, \\
\url{http://www.warwick.ac.uk/staff/J.E.Cremona/ftp/data/}

\bibitem[CFOSS]{ndescent}
J.E. Cremona, T.A. Fisher, C. O'Neil, D. Simon and M. Stoll, 
Explicit $n$-descent on elliptic curves, 
I Algebra, {\em J. reine angew. Math.}  {\bf 615} (2008) 121-155;
II Geometry,  {\em J. reine angew. Math.}  {\bf 632} (2009) 63-84;
{\em III Algorithms}, in preparation.

\bibitem[De]{De}
P. Deligne,
Courbes elliptiques: formulaire d'apr\`es J. Tate,
in {\em Modular functions of one variable, IV},
(Proc. Internat. Summer School, Univ. Antwerp, Antwerp, 1972),
pp. 53--73. Lecture Notes in Math., Vol. {\bf 476}, Springer, Berlin, 1975.

\bibitem[DS]{DS} 
Z. Djabri and N. P. Smart,
A comparison of direct and indirect methods for computing 
Selmer groups of an elliptic curve, in 
{\em  Algorithmic number theory (ANTS-III),} J. Buhler (ed.), 
Lecture Notes in Comput. Sci. {\bf{1423}}, Springer, Berlin, 1998, 502--513.

\bibitem[Do]{Dolg}
I. Dolgachev, Lectures on invariant theory,
LMS Lecture Note Series, 296, CUP, Cambridge, 2003. 

\bibitem[Du]{DujellaTables}
A. Dujella, High rank elliptic curves with prescribed torsion,
online table at \\
\url{http://web.math.hr/~duje/tors/tors.html}

\bibitem[Fi1]{ANTS}
T.A. Fisher, 
Testing equivalence of ternary cubics, in 
{\em Algorithmic number theory (ANS-VII)}, 
F. Hess, S. Pauli, M. Pohst (eds.), 
Lecture Notes in Comput. Sci. {\bf 4076}, Springer, 2006, 333-345.

\bibitem[Fi2]{minbqtc}
T.A. Fisher,  A new approach to minimising binary quartics and 
ternary cubics, {\em Math. Res. Lett.} {\bf 14} (2007) Issue 4, 597--613.

\bibitem[Fi3]{ANTS8}
T.A. Fisher,
Some improvements to 4-descent on an elliptic curve,
in {\em Algorithmic number theory (ANTS-VIII)},
A.J. van der Poorten, A. Stein (eds.),
Lecture Notes in Comput. Sci. {\bf 5011}, Springer, 2008, 125--138.

\bibitem[Fi4]{g1inv}
T.A. Fisher,
The invariants of a genus one curve, 
{\em Proc. Lond. Math. Soc.} (3) {\bf 97} (2008) 753-782. 

\bibitem[Fi5]{6and12}
T.A. Fisher,
Finding rational points on elliptic curves using 6-descent and 12-descent,  
{\em Journal of Algebra} {\bf{320}} (2008), no. 2, 853-884. 

\bibitem[Fi6]{Sha3tables}
T.A. Fisher, Elements of order 3 in the Tate-Shafarevich group,
online table at \\
\url{http://www.dpmms.cam.ac.uk/~taf1000/g1data/order3.html}

\bibitem[GPZ]{GPZ}
J. Gebel, A. Peth\H{o}  and H.G. Zimmer,
On Mordell's equation,
{\em Compositio Math.} {\bf{110}} (1998), no. 3, 335--367.

\bibitem[Hi]{Hilbert}
D. Hilbert, 
{\em Theory of algebraic invariants},
Cambridge University Press, Cambridge, 1993. 

\bibitem[HP]{HodgePedoe}
W.V.D. Hodge and D. Pedoe, {\em Methods of algebraic geometry},
Volume II, Reprint of the 1952 original, 
Cambridge University Press, Cambridge, 1994. 

\bibitem[Hu]{Hulek}
K. Hulek, {\em Projective geometry of elliptic curves}.  
Ast\'erisque  No. {\bf 137} (1986), 143 pp.

\bibitem[Ja]{Jacobson}
N. Jacobson, 
{\em Basic algebra I}, Second edition, 
W.H. Freeman and Company, New York, 1985. 

\bibitem[Ju]{Julia}
G. Julia, \'Etude sur les formes binaires non quadratiques 
\`a indetermin\'ees r\'eelles ou complexes, {\em Mem. Acad. Sci. l'Inst.
France} {\bf{55}} (1917) 1-293.

\bibitem[Ko]{Kollar}
J. Koll\'ar, 
Polynomials with integral coefficients, equivalent to a given polynomial,
{\em Electron. Res. Announc. Amer. Math. Soc.} {\bf{3}} (1997), 17--27 
(electronic).

\bibitem[Kr]{Kraus}
A. Kraus, 
Quelques remarques \`a propos des invariants $c\sb 4,\;c\sb 6$ et 
$\Delta$ d'une courbe elliptique. 
{\em Acta Arith.} {\bf{54}} (1989), no. 1, 75--80.

\bibitem[La]{Laska}
M. Laska, 
An algorithm for finding a minimal Weierstrass equation for an 
elliptic curve,
{\em Math. Comp.} {\bf{38}} (1982), no. 157, 257--260.

\bibitem[LLL]{LLL}
A.K. Lenstra, H.W. Lenstra and L. Lov\'asz, 
Factoring polynomials with rational coefficients,  
{\em Math. Ann.} {\bf 261} (1982), no. 4, 515--534.

\bibitem[Liu]{Liu}
Q. Liu, 
Mod\`eles entiers des courbes hyperelliptiques sur un corps de 
valuation discr\`ete,
{\em Trans. Amer. Math. Soc.} {\bf{348}} (1996), no. 11, 4577--4610.

\bibitem[M]{magma}
{\sf MAGMA} is described in W. Bosma, J. Cannon and C. Playoust, 
The \Magma algebra system I: The user language, 
{\em J. Symb. Comb.} {\bf{24}},  (1997) 235--265. 
(See also the \Magma home page at
\url{http://magma.maths.usyd.edu.au/magma/}.)

\bibitem[MSS]{MSS}
  J. R. Merriman,  S. Siksek and N. P. Smart,
 Explicit $4$-descents on an elliptic curve,
 {\em Acta Arith.} {\bf 77} (1996), no. 4, 385--404.

\bibitem[Mi]{MilneETALE}
J. S. Milne, {\em Lectures on Etale Cohomology}, Version 2.10, 
available online from \\
\url{http://www.jmilne.org/math/CourseNotes/math732.html}.

\bibitem[Po]{PoonenCubics}
B. Poonen, An explicit algebraic family of genus-one 
curves violating the Hasse principle,
{\em J. Th\'eor. Nombres Bordeaux} {\bf{13}} (2001), no. 1, 263--274. 

\bibitem[Sa]{Sadek}
M. Sadek, {\em Models of genus one curves}, PhD thesis, University
of Cambridge, in preparation.

\bibitem[ScSt]{SchaeferStoll}
E.F. Schaefer and M. Stoll,
How to do a $p$-descent on an elliptic curve,
{\em Trans. Amer. Math. Soc.} {\bf{356}} no. 3 (2004), 1209--1231

\bibitem[Se]{SerreLF}
J.-P. Serre, {\em Local fields}, Graduate Texts in Mathematics 
{\bf 67}, Springer-Verlag, New York-Berlin, 1979.

\bibitem[Sik]{SiksekThesis}
  S. Siksek,
  {\em Descent on curves of genus one},
  PhD thesis, University of Exeter, 1995.
See 
\url{http://www.warwick.ac.uk/staff/S.Siksek/papers/phdnew.pdf}

\bibitem[Sil1]{Si1}
J.H. Silverman, 
{\em The arithmetic of elliptic curves},
Graduate Texts in Mathematics {\bf{106}}, Springer-Verlag, New York, 1986. 

\bibitem[Sil2]{Si2}
J.H. Silverman, 
{\em Advanced topics in the arithmetic of elliptic curves},
Graduate Texts in Mathematics {\bf{151}}, Springer-Verlag, New York, 1994. 

\bibitem[Sk]{Skoro}
  A. Skorobogatov, {\em Torsors and rational points},
  Cambridge University Press, Cambridge, 2001.

\bibitem[SW]{StWat}
  W.A. Stein and M. Watkins,
  A database of elliptic curves---first report,
  {\em Algorithmic number theory} (Sydney 2002), 267--275,
  Lect. Notes Comp. Sci. {\bf 2369}, Springer, Berlin 2002.

\bibitem[SC1]{SCmin2}
  M. Stoll and J.E. Cremona,
Minimal models for 2-coverings of elliptic curves, 
{\em LMS J. Comput. Math.} {\bf{5}} (2002), 220--243.

\bibitem[SC2]{SCred}
  M. Stoll and J.E. Cremona,
On the reduction theory of binary forms,
  {\em J. reine angew. Math.} {\bf 565}, 79--99 (2003).

\bibitem[Sto]{StollMailing}
  M. Stoll, posting to the NMBRTHRY mailing list, 10 January 2002.
  See \\
\hbox{\url{http://listserv.nodak.edu/cgi-bin/wa.exe?A2=ind0201&L=NMBRTHRY&P=R500&I=-3}}.

\bibitem[Ta]{Tate}
J. Tate, Algorithm for determining the type of a singular 
fiber in an elliptic pencil, in {\em Modular functions of one variable IV}, 
B.J. Birch and W. Kuyk (eds.), Lecture Notes in Math.
{\bf{476}}, Springer, Berlin, 1975.
  
\bibitem[We1]{Weil1}
A. Weil, Remarques sur un m\'emoire d'Hermite,
{\em Arch. Math.} (Basel) {\bf{5}}, (1954) 197--202. 

\bibitem[We2]{Weil2}
A. Weil, 
Euler and the Jacobians of elliptic curves, 
{\em Arithmetic and geometry, Vol. I}, 353--359,
Progr. Math., {\bf{35}}, Birkh\"auser, Boston, MA, 1983.

\bibitem[Wo]{WomackThesis}
  T.O. Womack,
  {\em Explicit descent on elliptic curves},
  PhD thesis, University of Nottingham, 2003.
See \url{http://www.warwick.ac.uk/staff/J.E.Cremona/theses/}

\end{thebibliography}
\end{document}